\documentclass[11pt]{article}
\usepackage{amscd}
\usepackage{amssymb,amsfonts,amsmath,psfrag,amscd,cite}
\usepackage{amsmath}
  \usepackage{paralist}
  \usepackage{graphics}
  \usepackage{epsfig}
  \usepackage[colorlinks=true]{hyperref}
  \usepackage[all]{xy}
\usepackage{graphicx}
\newtheorem{theo}{Theorem}
\newtheorem{prop}[theo]{Proposition}

\newtheorem{defi}[theo]{Definition}
\newtheorem{lemm}[theo]{Lemma}
\newtheorem{rema}[theo]{Remark}

\makeatletter \@addtoreset{equation}{section}
\@addtoreset{theo}{section} \makeatother

\begin{document}
\date{}
\title{Symmetric Reduction and Hamilton-Jacobi Equations of\\ the Controlled Underwater Vehicle-Rotor System}
\author{Hong Wang  \\
School of Mathematical Sciences and LPMC,\\
Nankai University, Tianjin 300071, P.R.China\\
E-mail: hongwang@nankai.edu.cn \\\\
\emph{In Memory of Professor Jerrold E. Marsden} \\
June 25, 2020} \maketitle

{\bf Abstract.} In this paper, we first give the regular point
reduction and the two types of Hamilton-Jacobi equation for a regular
controlled Hamiltonian (RCH) system with symmetry and momentum map on the
generalization of a semidirect product Lie group. Next, as an
application of the theoretical results, we consider the underwater
vehicle with two internal rotors as a regular point reducible
RCH system, in the cases of coincident and
non-coincident centers of buoyancy and gravity,
we derive precisely the geometric constraint
conditions of the reduced symplectic form for the
dynamical vector field of the regular point reducible controlled
underwater vehicle-rotor system,
that is, the two types of Hamilton-Jacobi equation for the reduced
controlled underwater vehicle-rotor system by calculation in
detail,  respectively. These researches reveal the deeply internal
relationships of the geometrical structures of phase spaces, the dynamical
vector fields and controls of the system.\\



\tableofcontents

\section{Introduction}

It is well-known that the study of stability and drift of underwater
vehicle dynamics is a famous research work given by Leonard and Marsden
in \cite{lema97}, in which the underwater vehicle system is
considered as a Hamiltonian system with rigid motion symmetry,
the authors give the semidirect product reduction for
the underwater vehicle system and analysis carefully the stability for
the reduced underwater vehicle dynamics, also see
Leonard \cite{le97a, le97b}on the study of underwater vehicle system.\\

It is worthy of noting that the authors in Marsden et al.\cite{mawazh10}
define a regular controlled Hamiltonian (RCH) system, which is a
Hamiltonian system with external force and control.
In general, an RCH system, under the actions of external force and control, is not
Hamiltonian, however, it is a dynamical system closely related to a
Hamiltonian system, and it can be explored and studied by extending
the methods for external force and control in the study of Hamiltonian systems.
Thus, one can emphasize explicitly the impact of external force
and control in the study for the RCH systems.
Now, it is a natural idea that one may consider the
underwater vehicle with internal rotors, and a
control torque acting on the internal rotors,
as a model of a Hamiltonian system with control.
In particular, in Marsden et al.\cite{mawazh10}, the authors give the regular point
reduction and the regular orbit reduction for an RCH system with
symmetry, by analyzing carefully the geometrical and topological
structures of the phase space and the reduced phase space of the
corresponding Hamiltonian system.These research work not only gave a
variety of reduction methods for the RCH systems, but
also showed a variety of relationships of controlled Hamiltonian
equivalence of these systems. Thus, it is a natural problem
to study the underwater vehicle-rotor system
with the control torque acting on the internal rotors
as a regular point reducible RCH system on the generalization
of a semidirect product Lie group. Moreover,
as the application of theoretical result of the RCH system
with symmetry, one can derive the regular point reduced controlled
underwater vehicle-rotor systems by calculation in detail,
in the cases of coincident and
non-coincident centers of buoyancy and gravity,
respectively.\\

On the other hand, we note that
Hamilton-Jacobi theory is an important research subject
in mathematics and analytical mechanics. It provides a characterization
of the generating functions of certain time-dependent canonical
transformations, such that a given Hamiltonian system in such a form
that its solutions are extremely easy to find by reduction to the
equilibrium. Thus,
it is possible in many cases that Hamilton-Jacobi equation provides an
immediate way to integrate the equation of motion of a system, even
when the problem of Hamiltonian system itself has not been or cannot
be solved completely, see Abraham and Marsden \cite{abma78}, Arnold
\cite{ar89} and Marsden and Ratiu \cite{mara99}.
In addition, the Hamilton-Jacobi equation is
also fundamental in the study of the quantum-classical relationship
in quantization, and it also plays an important role
in the study of stochastic dynamical systems, see
Woodhouse \cite{wo92}, Ge and Marsden \cite{gema88},
and L\'{a}zaro-Cam\'{i} and Ortega \cite{laor09}.
For these reasons it is described as a useful tool in the study of
Hamiltonian system theory, and has been extensively developed in
past many years and become one of the most active subjects in the
study of modern applied mathematics and analytical mechanics.\\

Just as we have known that Hamilton-Jacobi theory from the
variational point of view is originally developed by Jacobi in 1866,
which state that the integral of Lagrangian of a mechanical system along the
solution of its Euler-Lagrange equation satisfies the
Hamilton-Jacobi equation. The classical description of this problem
from the generating function and the geometrical point of view is
given by Abraham and Marsden in \cite{abma78} as follows:
Let $Q$ be a smooth manifold and $TQ$
the tangent bundle, $T^* Q$ the cotangent bundle with a canonical
symplectic form $\omega$ and the projection $\pi_Q: T^* Q
\rightarrow Q. $
\begin{theo}
Assume that the triple $(T^*Q,\omega,H)$ is a Hamiltonian system
with Hamiltonian vector field $X_H$, and $W: Q\rightarrow
\mathbb{R}$ is a given generating function. Then the following two assertions
are equivalent:\\
\noindent $(\mathrm{i})$ For every curve $\sigma: \mathbb{R}
\rightarrow Q $ satisfying $\dot{\sigma}(t)= T\pi_Q
(X_H(\mathbf{d}W(\sigma(t))))$, $\forall t\in \mathbb{R}$, then
$\mathbf{d}W \cdot \sigma $ is an integral curve of the Hamiltonian
vector field $X_H$.\\
\noindent $(\mathrm{ii})$ $W$ satisfies the Hamilton-Jacobi equation
$H(q^i,\frac{\partial W}{\partial q^i})=E, $ where $E$ is a
constant.
\end{theo}

From the proof of the above theorem given in
Abraham and Marsden \cite{abma78}, we know that
the assertion $(\mathrm{i})$ with equivalent to
Hamilton-Jacobi equation by the generating function
gives a geometric constraint condition of the canonical symplectic form
on the cotangent bundle $T^*Q$
for Hamiltonian vector field of the system.
Thus, the Hamilton-Jacobi equation reveals the deeply internal relationships of
the generating function, the canonical symplectic form
and the dynamical vector field of a Hamiltonian system.\\

But, from Marsden et al.\cite{mawazh10} we know that,
the set of Hamiltonian systems with symmetries on a cotangent
bundle is not complete under the Marsden-Weinstein reduction.
Since the symplectic reduced system of a
Hamiltonian system with symmetry defined on the cotangent bundle
$T^*Q$ may not be a Hamiltonian system on a cotangent bundle,
then we cannot give the Hamilton-Jacobi theorem for the Marsden-Weinstein
reduced Hamiltonian system just like same as the above Theorem 1.1.
Note that the reduced space of the Marsden-Weinstein
reduced underwater vehicle system is a co-adjoint orbit of a
semidirect product Lie group, which is a submanifold of the dual space
of the corresponding semidirect product Lie algebra, and it may not be a
cotangent bundle. Thus, we cannot yet give the Hamilton-Jacobi equation
for the reduced underwater vehicle system given in Leonard and Marsden \cite{lema97}.
We have to look for a new way.\\

At the same time, we note also that since an RCH system defined on the cotangent bundle
$T^*Q$, in general,  may not be a Hamiltonian system,
and it has yet no generating function,
we cannot yet give the Hamilton-Jacobi theorem for the RCH system
and its regular reduced systems just like same as the above Theorem 1.1.
Thus, it is a natural problem how to describe precisely the Hamilton-Jacobi
theory for an RCH system and its regular reduced systems.
It is worthy of noting that, in Wang \cite{wa13d} the author derive precisely
the geometric constraint conditions
of the (reduced) symplectic forms for the
dynamical vector fields of an RCH system and its regular reduced systems.
These conditions
are called two types of Hamilton-Jacobi equations, which are the
development of the above two types of
Hamilton-Jacobi equations for a Hamiltonian system and its
Marsden-Weinstein reduced Hamiltonian system given in Wang
\cite{wa17}, and also the development of classical Hamilton-Jacobi equation
given in Abraham and Marsden \cite{abma78}.
Now, it is a natural problem
if there is a practical RCH system and how to show the effect on controls
in regular symplectic reductions of the system.
In particular, for the controlled underwater vehicle-rotor system
and its regular reduced systems. under the action of control,
they are not Hamiltonian, and have yet no generating function, however,
we want to derive precisely the geometric constraint conditions of
the (reduced) symplectic forms for the
dynamical vector fields of the regular point reducible controlled
underwater vehicle-rotor system, that is,
the Type I and Type II of Hamilton-Jacobi equation.
Moreover, when the underwater vehicle does not carry any internal rotor,
the reduced controlled underwater vehicle-rotors system is
just the Marsden-Weinstein reduced underwater vehicle system,
we hope that in the two cases, the two types of Hamilton-Jacobi equations
for the corresponding reduced systems are coincident.
These research are our goal in this paper.\\

A brief of outline of this paper is as follows. In the second
section, we first review some relevant basic facts
about underwater vehicle with two internal rotors,
and give the Hamiltonian function of the controlled
underwater vehicle-rotor system, in the cases of coincident
and non-coincident centers of buoyancy and gravity,
respectively, which will be used in subsequent sections.
In the third section, we first give the regular
point reduction of an RCH system with symmetry and momentum map
on the generalization of a semidirect product Lie group, then
derive precisely the geometric constraint conditions of
the reduced symplectic form for the dynamical vector fields
of the regular point reducible RCH system. As an application of the theoretical
result, in the fourth section and the fifth section,
we regard the controlled underwater vehicle-rotor
with a control torque acting on the two internal rotors
as a regular point
reducible RCH system on the generalization of semidirect
product Lie group $\textmd{SE}(3)\times S^1 \times S^1 $ and
$(\textmd{SE}(3)\circledS \mathbb{R}^3)\times S^1 \times S^1 $,
respectively, we give the regular point reduced controlled
underwater vehicle-rotor systems,
in the cases of coincident and non-coincident
centers of buoyancy and gravity. Moreover,
we derive precisely the geometric constraint conditions of
the reduced symplectic forms for the dynamical vector fields
of the regular point reducible controlled
underwater vehicle-rotor systems by calculation in detail,
that is, the two types of Hamilton-Jacobi equations for
these regular point reduced controlled underwater vehicle-rotor
systems.\\

It is worthy of noting that, in the cases of coincident and
non-coincident centers of buoyancy and gravity,
the motions of the controlled underwater vehicle-rotor system are different, and
the configuration spaces, the Hamiltonian functions, the actions of Lie groups,
the reduced symplectic forms and the reduced systems of
the controlled underwater vehicle-rotor system are also different. But,
the two types of Hamilton-Jacobi equations
given by calculation in detail are same,
that is, the internal rules are same by comparing Theorem 4.2 and Theorem 5.2.
Thus, these research work reveal the deeply internal
relationships of the geometrical structures of phase spaces, the dynamical
vector fields and controls of the controlled underwater vehicle-rotor system,
and develop the theory and application of the regular symplectic reduction
and Hamilton-Jacobi theory for the RCH systems
with symmetries, and make us have much deeper understanding and
recognition for the structure of Hamiltonian systems and RCH
systems.

\section{Underwater Vehicle with Two Internal Rotors}

In this section, we first give the the Hamiltonian of underwater vehicle
with two internal rotors, in the cases of
coincident and non-coincident centers of buoyancy and gravity,
respectively. We also review some relevant
basic facts about underwater vehicle
with two internal rotors, which will be used in subsequent sections.
We shall follow the notations and
conventions introduced in Leonard \cite{le97a, le97b}, Leonard and
Marsden \cite{lema97}, Marsden \cite{ma92}, Marsden and Ratiu
\cite{mara99}, and Marsden et al. \cite{mawazh10}. In this paper, we
assume that all manifolds are real, smooth and finite dimensional
and all actions are smooth left actions. For convenience, we also
assume that all controls appearing in this paper are the admissible
controls.

\subsection{With Coincident Centers of Buoyancy and Gravity}

We first describe a underwater vehicle carrying two internal "non-mass" rotors,
which is called a carrier body, moving in a given fluid,
where "non-mass" means that the mass of a rotor
is very very small comparing with the mass of the underwater vehicle.
We consider that the underwater vehicle-rotor system is a
neutrally buoyant, rigid body (often ellipsoidal) submerged in an
infinitely large volume of incompressible, inviscid, irrotational
fluid which is at rest at infinity. The dynamics of the carrier body-fluid
system are described by Kirchhoff's equations.
We first assume that the external forces
and torques acting on the underwater vehicle-rotor system are due to
buoyancy and gravity. In general, it is possible that the
center of buoyancy of the underwater vehicle-rotor system
may not be coincident with its center of gravity.
But, in this subsection we assume that the underwater vehicle
is symmetric and to have uniformly distributed mass, and
its center of buoyancy and its center of gravity are coincident.
Denote by $O$ the center of mass of the system in the carrier body frame and
at $O$ place a set of (orthogonal) body axes with origin located at the center of
buoyancy and axes aligned with the principal axes of the displaced
fluid, see Leonard and Marsden \cite{lema97}.
Next, we put two rotors within the underwater vehicle so that each rotor's rotation
axis is parallel to the first and the second principal axes of the
body, and when the body is oriented so that the body fixed frame is
aligned with the inertial frame, the third principal axis aligns
with the direction of gravity. Assume that the rotors spin
under the influence of a control torque $u$
acting on the rotors. If the carrier body moves in the given fluid,
and translation and rotation are considered, in this case,
the configuration space is $Q=W \times V$, where
$W=\textmd{SE}(3)= \textmd{SO}(3)\circledS \mathbb{R}^3$ and
$V=S^1\times S^1$, with the first factor being the attitude and
position of the underwater vehicle and the second factor being the
angles of rotors. The corresponding phase space is the cotangent
bundle $T^*Q$ and locally, $T^\ast Q=T^\ast \textmd{SE}(3)\times T^\ast V$, where
$T^\ast V = T^\ast (S^1\times S^1)\cong T^\ast \mathbb{R}^2$ locally, with the
canonical symplectic form $\omega_Q$.
By using the local left trivialization, locally,
$T^\ast \textmd{SE}(3)\cong \textmd{SE}(3)\times \mathfrak{se}^\ast(3)$
and $T^*\mathbb{R}^2 \cong \mathbb{R}^2 \times \mathbb{R}^{2*}$,
then we have that locally,
$T^*Q \cong \textmd{SE}(3)\times \mathfrak{se}^\ast(3)
\times \mathbb{R}^2 \times \mathbb{R}^{2*}$.
For convenience, in the following we denote uniformly that, locally,
$Q= \textmd{SE}(3)\times \mathbb{R}^2, $ and
$T^* Q= T^*(\textmd{SE}(3)\times \mathbb{R}^2)\cong \textmd{SE}(3)\times \mathfrak{se}^\ast(3)
\times \mathbb{R}^2 \times \mathbb{R}^{2*}$.\\

Assume that the matrix of inertia moment
of the carrier body-fluid system is denoted by
$I=\textmd{diag}(I_1,I_2,I_3)$ and the mass matrix by
$M=\textmd{diag}(m_1,m_2,m_3)$, in the carrier body fixed frame,
which is a principal body frame. Here these matrices include the
"added" inertias and masses due to the fluid, see Leonard \cite{le97a,
le97b} and Leonard and Marsden \cite{lema97}.
Let $J_k, k=1,2$ be the moments of
inertia of rotors around their rotation axes. Let $J_{ki},\;
k=1,2,\; i=1,2,3,$ be the moments of inertia of the $k$-th rotor
with $k=1,2$ around the $i$-th principal axis with $i=1,2,3,$
respectively, and denote by $\bar{I}_k=I_k+J_{1k}+J_{2k}-J_{kk}, \;
k=1,2$, and $\bar{I}_3=I_3+J_{13}+J_{23}$, see
Marsden \cite{ma92}. Let
$\Omega=(\Omega_1,\Omega_2,\Omega_3)$ and $v=(v_1,v_2,v_3)$ be the
angular and linear velocity vectors of underwater vehicle-rotors system computed
with respect to the axes fixed in the carrier body and
$(\Omega_1,\Omega_2,\Omega_3)\in \mathfrak{so}(3)$,
$v=(v_1,v_2,v_3) \in \mathbb{R}^3$, and
$(\Omega,v)\in \mathfrak{se}(3)= \mathfrak{so}(3)\circledS \mathbb{R}^3$.
Let $\theta_k,\; k=1,2,$ be the relative angles of rotors and
$\dot{\theta}=(\dot{\theta_1},\dot{\theta_2})$ the relative angular
velocity vector of rotor about the principal axes with respect to
the carrier body fixed frame of the underwater vehicle-rotor system.
For convenience, we assume the total mass of the system $m=1$.\\

Now, by the local left trivialization, locally,
$T\textmd{SE}(3)\cong \textmd{SE}(3)\times \mathfrak{se}(3)$
and $T\mathbb{R}^2 \cong \mathbb{R}^2 \times \mathbb{R}^{2}$,
then we have that locally,
$TQ \cong \textmd{SE}(3)\times \mathfrak{se}(3)
\times \mathbb{R}^2 \times \mathbb{R}^{2}$.
We consider the Lagrangian
 $L(A,c,\Omega,v,\theta,\dot{\theta}):
TQ \cong \textmd{SE}(3)\times\mathfrak{se}(3)\times\mathbb{R}^2\times\mathbb{R}^2\to
\mathbb{R}$, which is that the total kinetic energy of the underwater
vehicle-fluid plus the total kinetic energy of rotors of the system, given
by
$$L(A,c,\Omega,v,\theta,\dot{\theta})
=\dfrac{1}{2}[\bar{I}_1\Omega_1^2+\bar{I}_2\Omega_2^2+\bar{I}_3\Omega_3^2
+ m_1v_1^2 + m_2v_2^2+ m_3v_3^2
+J_1(\Omega_1+\dot{\theta}_1)^2+J_2(\Omega_2+\dot{\theta}_2)^2],$$
where $(A,c)\in \textmd{SE}(3)$, $(\Omega,v)\in \mathfrak{se}(3)$
and $\Omega=(\Omega_1,\Omega_2,\Omega_3)\in \mathfrak{so}(3)$,
$v\in\mathbb{R}^3$, $\theta=(\theta_1,\theta_2)\in \mathbb{R}^2$,
$\dot{\theta}=(\dot{\theta}_1,\dot{\theta}_2)\in \mathbb{R}^2$. If
we introduce the conjugate angular momentum and linear momentum,
which is given by
\begin{align*}
& \Pi_k= \dfrac{\partial L}{\partial
\Omega_k}=\bar{I}_k\Omega_k+J_k(\Omega_k+\dot{\theta}_k), \;\;\;\;\;\;
\Pi_3=\dfrac{\partial L}{\partial \Omega_3}=\bar{I}_3\Omega_3, \\
& P_i= \dfrac{\partial L}{\partial v_i}= m_iv_i, \;\;\;\; i=1,2,3, \;\;\;\;\;\;
l_k=\dfrac{\partial L}{\partial
\dot{\theta}_k}=J_k(\Omega_k+\dot{\theta}_k), \;\;\;\; k=1,2.
\end{align*}
By the Legendre transformation \begin{align*} FL:
TQ\cong \textmd{SE}(3)\times\mathfrak{se}(3)\times\mathbb{R}^2\times\mathbb{R}^2
& \to T^* Q\cong \textmd{SE}(3)\times
\mathfrak{se}^\ast(3)\times\mathbb{R}^2\times\mathbb{R}^{2*},\\
(A,c,\Omega,v,\theta,\dot{\theta}) & \to (A,c,\Pi,P,\theta,l),
\end{align*} where
$\Pi=(\Pi_1,\Pi_2,\Pi_3)\in \mathfrak{so}^\ast(3), \; P
=(P_1,P_2,P_3) \in \mathbb{R}^{3*}$,
and $(\Pi, P)\in
\mathfrak{se}^\ast(3)=  \mathfrak{so}^*(3)\circledS \mathbb{R}^{3*}$,
$l=(l_1,l_2)\in \mathbb{R}^{2*}$, we
have the Hamiltonian $H(A,c,\Pi,P,\theta,l): T^*Q \cong\textmd{SE}(3)\times
\mathfrak{se}^\ast (3)\times\mathbb{R}^2\times\mathbb{R}^{2*} \to
\mathbb{R}$ given by
\begin{align} &H(A,c,\Pi,P,\theta,l)=\Omega\cdot
\Pi+ v\cdot P+ \dot{\theta}\cdot l-L(A,c,\Omega,v,\theta,\dot{\theta}) \nonumber \\
&=\frac{1}{2}[\frac{(\Pi_1-l_1)^2}{\bar{I}_1}+\frac{(\Pi_2-l_2)^2}{\bar{I}_2}
+\frac{\Pi_3^2}{\bar{I}_3}+
\frac{P_1^2}{m_1}+\frac{P_2^2}{m_2}+\frac{P_3^2}{m_3}+\frac{l_1^2}{J_1}+\frac{l_2^2}{J_2}].
\label{2.1}
\end{align}

In this case, in order to give the dynamical vector field and the two types of
Hamilton-Jacobi equation of the controlled underwater vehicle-rotor
system, we need to consider the regular point reduction of
the controlled underwater vehicle-rotor system and
give precisely the $R_p$-reduced symplectic form of the cotangent bundle $T^* Q\cong
\textmd{SE}(3)\times \mathfrak{se}^*(3)\times\mathbb{R}^2\times\mathbb{R}^{2*}$.

\subsection{With Non-Coincident Centers of Buoyancy and Gravity }

Since it is possible that the carrier body's center of
buoyancy may not be coincident with its center of gravity, in this
subsection then we consider the underwater vehicle-rotor system with non-coincident centers
of buoyancy and gravity. We fix an orthogonal coordinate frame to
the carrier body with origin located at the center of buoyancy and axes
aligned with the principal axes of the displaced fluid, and these
axes are also the principal axes of the carrier body, since the carrier body
is regarded as an ellipsoidal and it is symmetric about these axes.
We put two rotors within the underwater vehicle so that each rotor's rotation
axis is parallel to the first and the second principal axes of the
body, and when the carrier body is oriented so that
the carrier body-fixed frame is aligned with the inertial frame, the third
principal axis aligns with the direction of gravity.
The vector from the center of buoyancy to the center of gravity with respect to the
carrier body-fixed frame is $h\chi$, where $\chi$ is an unit vector on the
line connecting the two centers which is assumed to be aligned along
the third principal axis, and $h$ is the length of this segment, see
Leonard and Marsden \cite{lema97}. Assume that the total mass of the
underwater vehicle-rotor system
$m=1$, and the magnitude of gravitational acceleration
is denoted $g$, and let $\Gamma$ be the unit vector viewed by an observer
moving with the carrier body, and the rotors
spin under the influence of a control torque $u$ acting on the rotors..
In this case, the configuration space is $Q=W \times V$, where
$W= \textmd{SE}(3)\circledS \mathbb{R}^3= (\textmd{SO}(3)\circledS
\mathbb{R}^3)\circledS \mathbb{R}^3$ is a double semidirect product Lie group
and $V=S^1\times S^1$, with the first factor being the attitude and
position of underwater vehicle as well as the drift of the underwater
vehicle-rotor in the rotational and translational process, and the second
factor being the angles of rotors.
The corresponding phase space is the cotangent bundle
$T^*Q$ and locally, $T^\ast Q=T^\ast (\textmd{SE}(3)\circledS \mathbb{R}^3)\times T^\ast
V$, where $T^\ast V = T^\ast (S^1\times S^1)\cong T^\ast
\mathbb{R}^2 $ locally, with the
canonical symplectic form $\omega_Q$.
By using the local left trivialization, locally,
$T^\ast (\textmd{SE}(3)\circledS \mathbb{R}^3)\cong (\textmd{SE}(3)\circledS \mathbb{R}^3)
\times (\mathfrak{se}^\ast(3)\circledS \mathbb{R}^{3*})$£¬
and $T^*\mathbb{R}^2 \cong \mathbb{R}^2 \times \mathbb{R}^{2*}$,
then we have that locally,
$T^*Q \cong (\textmd{SE}(3)\circledS \mathbb{R}^3)\times (\mathfrak{se}^\ast(3)\circledS \mathbb{R}^{3*})
\times \mathbb{R}^2 \times \mathbb{R}^{2*}$.
For convenience, in the following we denote uniformly that, locally,
$Q= \textmd{SE}(3)\circledS \mathbb{R}^3 \times \mathbb{R}^2, $ and
$T^* Q= T^*(\textmd{SE}(3)\circledS \mathbb{R}^3 \times \mathbb{R}^2)
\cong \textmd{SE}(3)\circledS \mathbb{R}^3\times \mathfrak{se}^\ast(3)\circledS \mathbb{R}^{3*}
\times \mathbb{R}^2 \times \mathbb{R}^{2*}$.\\

Now, by the local left trivialization, locally,
$T(\textmd{SE}(3)\circledS \mathbb{R}^3)\cong
(\textmd{SE}(3)\circledS \mathbb{R}^3)\times (\mathfrak{se}(3)\circledS \mathbb{R}^3)$
and $T\mathbb{R}^2 \cong \mathbb{R}^2 \times \mathbb{R}^{2}$,
then we have that locally,
$TQ \cong (\textmd{SE}(3)\circledS \mathbb{R}^3)\times (\mathfrak{se}(3)\circledS \mathbb{R}^3)
\times \mathbb{R}^2 \times \mathbb{R}^{2}$.
We consider the Lagrangian
$L(A,c,b,\Omega,\Gamma,v,\theta,\dot{\theta}): TQ \cong
(\textmd{SE}(3)\circledS \mathbb{R}^3)\times (\mathfrak{se}(3)\circledS \mathbb{R}^3)
\times\mathbb{R}^2\times\mathbb{R}^2\to \mathbb{R}$,
which is that the total kinetic energy of the underwater vehicle-fluid plus the
total kinetic energy of rotors, and minus potential energy of
the underwater vehicle-rotor system, given by
\begin{align*}
L(A,c,b,\Omega,\Gamma,v,\theta,\dot{\theta}) = &
\dfrac{1}{2}[\bar{I}_1\Omega_1^2+\bar{I}_2\Omega_2^2+\bar{I}_3\Omega_3^2
+ m_1v_1^2 + m_2v_2^2+ m_3v_3^2 \\
& \;\;\;
+J_1(\Omega_1+\dot{\theta}_1)^2+J_2(\Omega_2+\dot{\theta}_2)^2]-
gh\Gamma \cdot \chi,
\end{align*}
where $(A,c)\in \textmd{SE}(3)$,
$(\Omega,\Gamma)\in \mathfrak{se}(3)= \mathfrak{so}(3)\circledS \mathbb{R}^3$ and
$\Omega=(\Omega_1,\Omega_2,\Omega_3)\in \mathfrak{so}(3)$, $b, v,
\in \mathbb{R}^3$, $\theta=(\theta_1,\theta_2)\in
\mathbb{R}^2$, $\dot{\theta}=(\dot{\theta}_1,\dot{\theta}_2)\in
\mathbb{R}^2$, and the variable $\Gamma \in\mathbb{R}^3$ is regarded as a
parameter with respect to potential energy of the system,
$(\Omega, \Gamma)\in \mathfrak{se}(3)$.
If we introduce the conjugate angular momentum and
linear momentum, which is given by
\begin{align*}
& \Pi_k= \dfrac{\partial L}{\partial
\Omega_k}=\bar{I}_k\Omega_k+J_k(\Omega_k+\dot{\theta}_k), \;\;\;\;\;\;
\Pi_3=\dfrac{\partial L}{\partial \Omega_3}=\bar{I}_3\Omega_3, \\
& P_i= \dfrac{\partial L}{\partial v_i}= m_iv_i, \;\;\;\; i=1,2,3, \;\;\;\;\;\;
l_k=\dfrac{\partial L}{\partial \dot{\theta}_k}=J_k(\Omega_k+\dot{\theta}_k), \;\;\;\; k=1,2.
\end{align*}
By the Legendre transformation with the parameter $\Gamma$, that is,
\begin{align*} FL: TQ\cong (\textmd{SE}(3)\circledS \mathbb{R}^3)\times (\mathfrak{se}(3)\circledS \mathbb{R}^3)\times\mathbb{R}^2\times\mathbb{R}^2 & \to T^* Q\cong
(\textmd{SE}(3)\circledS \mathbb{R}^3)\times (\mathfrak{se}^*(3)\circledS \mathbb{R}^{3*})
\times\mathbb{R}^2\times\mathbb{R}^{2*},\\
(A,c,b,\Omega,\Gamma,v,\theta,\dot{\theta}) &
\to (A,c,b,\Pi,\Gamma,P,\theta,l), \end{align*} where
$\Pi=(\Pi_1,\Pi_2,\Pi_3)\in \mathfrak{so}^\ast(3), \; P
=(P_1,P_2,P_3)\in \mathbb{R}^{3*}$, and $(\Pi, \Gamma)\in
\mathfrak{se}^\ast(3)$, and $l=(l_1,l_2)\in \mathbb{R}^{2*}$, we have the
Hamiltonian $H(A,c,b,\Pi,\Gamma,P,\theta,l):T^* Q\cong
(\textmd{SE}(3)\circledS \mathbb{R}^3)\times (\mathfrak{se}^*(3)\circledS \mathbb{R}^{3*})
\times\mathbb{R}^2\times\mathbb{R}^{2*}\to \mathbb{R}$
given by
\begin{align} &H(A,c,b,\Pi,\Gamma,P,\theta,l)=\Omega\cdot
\Pi+ v\cdot P+ \dot{\theta}\cdot l-L(A,c,b,\Omega,\Gamma,v,\theta,\dot{\theta}) \nonumber \\
&=\frac{1}{2}[\frac{(\Pi_1-l_1)^2}{\bar{I}_1}+\frac{(\Pi_2-l_2)^2}{\bar{I}_2}
+\frac{\Pi_3^2}{\bar{I}_3}+
\frac{P_1^2}{m_1}+\frac{P_2^2}{m_2}+\frac{P_3^2}{m_3}+\frac{l_1^2}{J_1}+\frac{l_2^2}{J_2}]+
gh\Gamma \cdot \chi \;.
 \label{2.2}
\end{align}

In this case, in order to give the dynamical vector field and the two types of
Hamilton-Jacobi equation of the controlled underwater vehicle-rotor
system, we need to consider the regular point reduction of
the controlled underwater vehicle-rotor system and
give precisely the $R_p$-reduced symplectic form of the cotangent bundle $T^* Q\cong
(\textmd{SE}(3)\circledS \mathbb{R}^3)\times (\mathfrak{se}^*(3)\circledS \mathbb{R}^{3*})
\times\mathbb{R}^2\times\mathbb{R}^{2*}$.

\section{RCH System on a Generalization of Semidirect Product Lie Group}

In order to describe the regular point reduction and the two types of
Hamilton-Jacobi theorem of the controlled underwater
vehicle-rotor system, in this section we need to first give the regular point
reduction and the two types of Hamilton-Jacobi theorem
of an RCH system with symmetry and momentum map on the
generalization of a semidirect product Lie group $Q=W\times V$,
where $W=G\circledS E$ is a semidirect product Lie group
with Lie algebra $\mathfrak{w}=\mathfrak{g}\circledS E$, $G$ is a
Lie group with Lie algebra $\mathfrak{g}$, $E$ is an $r$-dimensional
vector space and $V$ is a $k$-dimensional vector space. See Marsden
et al. \cite{mamiorpera07, mamora90, marawe84a, marawe84b }.

\subsection{Symmetric Reduction }

The reduction theory for the mechanical system with symmetry is an
important subject and it is widely studied in the theory of
mathematics and mechanics, as well as applications. The main goal of
reduction theory in mechanics is to use conservation laws and the
associated symmetries to reduce the number of dimensions of a
mechanical system required to be described. So, such reduction
theory is regarded as a useful tool for simplifying and studying
concrete mechanical systems. Over forty years ago, the regular
symplectic reduction for the Hamiltonian system with symmetry and
coadjoint equivariant momentum map was set up by famous professors
Jerrold E. Marsden and Alan Weinstein, which is called
Marsden-Weinstein reduction, and great developments have been
obtained around the work in the theoretical study and applications
of mathematics, mechanics and physics; see
Abraham and Marsden \cite{abma78}, Abraham et al.
\cite{abmara88}, Arnold \cite{ar89}, de Le\'{o}n and Rodrigues \cite{lero89},
Libermann and Marle \cite{lima87}, Marsden \cite{ma92}, Marsden et al.
\cite{mamiorpera07, mamora90}, Marsden and Perlmutter \cite{mape00},
Marsden and Ratiu \cite{mara99}, Marsden and
Weinstein \cite{mawe74}, Meyer \cite{me73},
Nijmeijer and Van der Schaft \cite {nivds90} and Ortega and Ratiu \cite{orra04}.
In particular, the authors in Marsden et al. \cite{mawazh10}
set up the regular reduction theory for the RCH systems
with symplectic structures and symmetries on a symplectic fiber
bundle, as an extension of Marsden-Weinstein reduction theory of
Hamiltonian systems under regular controlled Hamiltonian equivalence
conditions, and from the viewpoint of completeness of regular symplectic
reduction and by analyzing carefully the geometrical and topological
structures of the phase space and the reduced phase space of the
corresponding Hamiltonian system.
Some developments around the work are given in
Wang and Zhang \cite{wazh12}, Ratiu and Wang \cite{rawa12},
Wang \cite{wa15a}.\\

In the following we first give the regular point reduction of an RCH system
with symmetry and momentum map on the generalization of a semidirect product Lie group.
Assume that Lie group $G$ acts on the left by linear maps on $E$, and $G$ also acts
on the left on the dual space $E^\ast$ of $E$, and the action by an
element $g$ on $E^\ast$ is the transpose of the action of $g^{-1}$
on $E$. As a set, the underlying manifold of $W$ is $G\times E$ and
the multiplication on $W$ is given by
\begin{equation}
  (g_1,x_1)(g_2,x_2):=(g_1g_2,x_1+\sigma(g_1)x_2),\quad g_1,g_1\in
  G,\quad x_1,x_2\in E, \label{3.1}
\end{equation}
where $\sigma:G\to \operatorname{Aut}(E)$ is a representation of the
Lie group $G$ on $E$, $\operatorname{Aut}(E)$ denotes the Lie group
of linear isomorphisms of $E$ onto itself, whose Lie algebra is
$\operatorname{End}(E)$, the space of all linear maps of $E$ to
itself.\\

The Lie algebra of $W$ is the semidirect product of Lie algebras
$\mathfrak{w}=\mathfrak{g}\circledS E$, $\mathfrak{w}^\ast$ is the
dual of $\mathfrak{w}$, that is, $\mathfrak{w}^\ast=
(\mathfrak{g}\circledS E)^\ast \cong \mathfrak{g}^\ast\circledS E^*$. The underlying vector space of
$\mathfrak{w}$ is $\mathfrak{g}\times E$ and the Lie bracket on
$\mathfrak{w}$ is given by
\begin{equation} [(\xi_1,v_1),(\xi_2,v_2)]=([\xi_1,\xi_2],\sigma
'(\xi_1)v_2-\sigma '(\xi_2)v_1),\quad \forall \xi_1, \; \xi_2 \in
\mathfrak{g},\quad v_1, \; v_2 \in E, \label{3.2}
\end{equation}
where $\sigma ':\mathfrak{g} \to \operatorname{End}(E)$ is the
induced Lie algebra representation given by
\begin{equation}
\sigma'(\xi)v:=\left.\frac{\mathrm{d}}{\mathrm{d}
t}\right|_{t=0}\sigma(\exp t\xi)v,\quad \xi \in \mathfrak{g},\quad v
\in E. \label{3.3}
\end{equation}
Identify the underlying vector space of $\mathfrak{w}^\ast$ with
$\mathfrak{g}^\ast\times E^\ast\cong \mathfrak{g}^\ast\times E, $ by
using the duality pairing on each factor. One can give the formula
for the Lie-Poisson bracket on the dual of semidirect product
$\mathfrak{w}^\ast = (\mathfrak{g}\circledS E)^\ast
\cong \mathfrak{g}^\ast\circledS E^*$ as follows,
that is, for $F,K: \mathfrak{w}^\ast \to \mathbb{R}$, their
(-)-semidirect product Lie-Poisson bracket is given by
\begin{equation}
\{F,K\}_{\mathfrak{w}^*}(\mu,a)= -\langle\mu,[\frac{\delta F}{\delta
\mu},\frac{\delta K}{\delta \mu}]\rangle- \langle a,\frac{\delta
F}{\delta \mu}\cdot\frac{\delta K}{\delta a}- \frac{\delta K}{\delta
\mu}\cdot\frac{\delta F}{\delta a}\rangle, \label{3.4}
\end{equation}
where $(\mu,a)\in \mathfrak{w}^\ast$ and $\dfrac{\delta F}{\delta
\mu}\in \mathfrak{g}$, $\dfrac{\delta F}{\delta a}\in E$ are the
functional derivatives. Moreover, the Hamiltonian vector field of a
smooth function $H:\mathfrak{w}^\ast \to \mathbb{R}$ is given by
\begin{equation}
X_H(\mu,a)=(\operatorname{ad}_{\delta H/\delta \mu}^\ast
\mu-\rho_{\delta H/\delta a}^\ast a, \; \frac{\delta H}{\delta
\mu}\cdot a), \label{3.5}
\end{equation}
where the infinitesimal action of $\mathfrak{g}$ on $E$ can be
denoted by $\xi \cdot v=\rho_v(\xi)$, for any $\xi \in
\mathfrak{g}$, $v \in E$ and the map $\rho_v: \mathfrak{g}\to E$ is
the derivative of the map $g \mapsto gv$ at the identity and
$\rho_v^\ast: E^\ast \to \mathfrak{g}^\ast$ is its dual,
see Marsden et al. \cite{mamiorpera07}.\\

At first, we consider a symplectic action of $W$ on a symplectic
manifold $P$ and assume that this action has an
$\operatorname{Ad}^\ast$-equivariant momentum map $\mathbf{J}_W:P\to
\mathfrak{w}^\ast$. On the one hand, we can regard $E$ as a normal
subgroup of $W$, it also acts on $P$ and has a momentum map
$\mathbf{J}_E:P\to E^\ast$ given by $\mathbf{J}_E=i_E^\ast\cdot
\mathbf{J}_W$, where $i_E:E\to \mathfrak{w};\;x\mapsto (0,x)$ is the
inclusion, and $i_E^\ast: \mathfrak{w}^\ast \to E^\ast$ is its dual.
$\mathbf{J}_E$ is called the second component of $\mathbf{J}_W$. On
the other hand, we can also regard $G$ as a subgroup of $W$, and
have the inclusion $i_{\mathfrak{g}}: \mathfrak{g} \to
\mathfrak{w}$, given by $\xi \mapsto (\xi,0)$. Thus, $G$-action also
has a momentum map $\mathbf{J}_G: P \to \mathfrak{g}^\ast$ given by
$\mathbf{J}_G=i_{\mathfrak{g}}^\ast\cdot \mathbf{J}_W$, where
$i_{\mathfrak{g}}^\ast: \mathfrak{w}^\ast \to \mathfrak{g}^\ast$ is
its dual of $i_{\mathfrak{g}}$, which is called the first component
of $\mathbf{J}_W$. Moreover, from the
$\operatorname{Ad}^\ast$-equivariance of $\mathbf{J}_W$ under
$W$-action, we know that $\mathbf{J}_E$ is also
$\operatorname{Ad}^\ast$-equivariant under $W$-action. Thus, we can
carry out reduction of $P$ by $W$ at a regular value $\tau
=(\mu,a)\in \mathfrak{w}^\ast$ of the momentum map $\mathbf{J}_W$ in
two stages using the following procedure. (i)First reduce $P$ by $E$
at the value $a\in E^\ast$, and get the reduced space
$P_a=\mathbf{J}_E^{-1}(a)/E$. Since the reduction is by the Abelian
group $E$, so the quotient is done using the whole of $E$. (ii)The
isotropy subgroup $W_a\subset W$, consists of elements of $W$ that
leave the point $(0,a)\in \mathfrak{w}^\ast$, where $a\in E^\ast$,
fixed using the action of $W$ on
$E^\ast$. One can prove that the group $W_a$ leaves the set
$\mathbf{J}_E^{-1}(a)\subset P$ invariant, and acts symplectically
on the reduced space $P_a$ and has a naturally induced momentum map
$\mathbf{J}_a:P_a\to \mathfrak{w}_a^\ast$, where $\mathfrak{w}_a$ is
the Lie algebra of the isotropy subgroup $W_a$ and
$\mathfrak{w}_a^\ast$ is its dual. (iii)Reduce the first reduced
space $P_a$ at the point $\mu_a=\mu|_{\mathfrak{w}^\ast_a}\in
\mathfrak{w}_a^\ast$, one can get the second reduced space
$(P_a)_{\mu_a}=\mathbf{J}_a^{-1}(\mu_a)/(W_a)_{\mu_a}$. Thus, we can
give the following proposition on the reduction by stages for
semidirect product Lie group, see Marsden et al.
\cite{mamiorpera07}.
\begin{prop}
The reduced space $(P_a)_{\mu_a}$ is symplectically diffeomorphic to
the reduced space $P_\tau$ obtained by reducing $P$ by $W$ at the
regular point $\tau=(\mu,a)\in\mathfrak{w}^\ast$.
\end{prop}
In particular, we can choose that $P=T^\ast W$, where $W=G\circledS
E$ is a semidirect product Lie group, with the cotangent lift action
of $W$ on $T^\ast W$ induced by left translation of $W$ on itself.
Since the reduction of $T^\ast W$ by the action of $E$ can give a
space which is isomorphic to $T^\ast G$, from the above proposition
on reduction by stages for semidirect product Lie group, we can get the
following proposition.
\begin{prop}
The reduction of $T^\ast G$ by $W_a$ at the regular values
$\mu_a=\mu|_{\mathfrak{w}^\ast_a}$ gives a space which is isomorphic
to the coadjoint orbit $\mathcal{O}_{(\mu,a)}
\subset\mathfrak{w}^\ast$ through the point $(\mu,a)\in
\mathfrak{w}^\ast$, where $\mathfrak{w}^\ast$ is the dual of the Lie
algebra $\mathfrak{w}$ of $W$.
\end{prop}

Next, we consider the action of $W$ on the generalization of a
semidirect product Lie group $Q=W\times V$ and define the left
$W$-action $\Phi: W\times Q \rightarrow Q, \;
\Phi((g_1,x_1),((g_2,x_2), \theta)):=((g_1,x_1)(g_2,x_2),\theta)$,
for any $(g_1,x_1), \; (g_2,x_2) \in W,\; \theta \in V$, that is,
the $W$-action on $Q$ is the left translation on the first factor
$W$, and $W$ acts trivially on the second factor $V$. Because locally,
$T^\ast Q = T^\ast W \times T^\ast V$, by using the left trivializations of $T^\ast W$ and $T^\ast V$,
that is, $T^\ast W \cong W \times \mathfrak{w}^\ast, $ where $\mathfrak{w}^\ast$
is the dual of $\mathfrak{w}= \mathfrak{g}\circledS E$,
and $T^\ast V \cong V\times V^\ast$. Hence we
have that $T^\ast Q \cong W \times \mathfrak{w}^\ast \times V \times
V^\ast$. If the left $W$-action $\Phi: W\times Q \rightarrow Q $ is
free and proper, then the cotangent lift of the action to its
cotangent bundle $T^\ast Q$, given by $\Phi^{T*}: W \times T^*Q
\rightarrow T^*Q, \;
\Phi^{T*}((g_1,x_1),((g_2,x_2),(\mu,a),\theta,\lambda)):=((g_1,x_1)(g_2,x_2),(\mu,a),\theta,\lambda)$,
for any $(g_1,x_1), \; (g_2,x_2) \in W,\; (\mu,a) \in
\mathfrak{w}^\ast, \; \theta \in V, \; \lambda \in V^\ast$, is also
a free and proper action, and the orbit space $(T^\ast Q)/ W $ is a
smooth manifold and $\pi_W: T^*Q \rightarrow (T^*Q )/W $ is a smooth
submersion. Since $W$ acts trivially on $\mathfrak{w}^\ast$, $V$ and
$V^\ast$, it follows that $(T^\ast Q)/ W$ is diffeomorphic to
$\mathfrak{w}^\ast \times V \times V^\ast$.\\

We know that $\mathfrak{w}^*= (\mathfrak{g}\circledS E)^\ast$ is a
Poisson manifold with respect to the (-)-semidirect product
Lie-Poisson bracket $\{\cdot,\cdot\}_{\mathfrak{w}^*}$ defined by
$(3.4)$. For $(\mu,a) \in \mathfrak{w}^\ast$, the co-adjoint orbit
$\mathcal{O}_{(\mu,a)} \subset \mathfrak{w}^\ast$ has the induced
orbit symplectic form $\omega^{-}_{\mathcal{O}_{(\mu,a)}}$ given by
\begin{equation}
\omega_{\mathcal{O}_{(\mu,a)}}^{-}(X_{F_{(\mu,a)}},
X_{K_{(\mu,a)}})=
\{F_{(\mu,a)},K_{(\mu,a)}\}_{\mathfrak{w}^*}|_{\mathcal{O}_{(\mu,a)}},
 \label{3.6}
\end{equation}
for $F_{(\mu,a)},K_{(\mu,a)}: \mathcal{O}_{(\mu,a)} \to \mathbb{R}$,
which is the restriction of the (-)-semidirect product Lie-Poisson
brackets on $\mathfrak{w}^\ast$ to the co-adjoint orbit
$\mathcal{O}_{(\mu,a)}$. From the Symplectic Stratification theorem
we know that the co-adjoint orbits $(\mathcal{O}_{(\mu,a)},
\omega_{\mathcal{O}_{(\mu,a)}}^{-}), \; (\mu,a)\in
\mathfrak{w}^\ast, $ form the symplectic leaves of the Poisson
manifold $(\mathfrak{w}^\ast,\{\cdot,\cdot\}_{\mathfrak{w}^*}). $
Let $\omega_V$ be the canonical symplectic form on $T^\ast V \cong V
\times V^\ast$ given by
$$\omega_V((\theta_1, \lambda_1),(\theta_2,
\lambda_2))=<\lambda_2,\theta_1> -<\lambda_1,\theta_2>,$$ where
$(\theta_i, \lambda_i)\in V\times V^\ast, \; i=1,2$, $<\cdot,\cdot>$
is the natural pairing between $V^\ast$ and $V$. Thus, we can induce
a symplectic form $\tilde{\omega}^{-}_{\mathcal{O}_{(\mu,a)} \times
V \times V^\ast}= \pi_{\mathcal{O}_{(\mu,a)}}^\ast
\omega^{-}_{\mathcal{O}_{(\mu,a)}}+ \pi_V^\ast \omega_V$ on the
smooth manifold $\mathcal{O}_{(\mu,a)} \times V \times V^\ast$,
where the maps $\pi_{\mathcal{O}_{(\mu,a)}}: \mathcal{O}_{(\mu,a)}
\times V \times V^\ast \to \mathcal{O}_{(\mu,a)}$ and $\pi_V:
\mathcal{O}_{(\mu,a)} \times V \times V^\ast \to V\times V^\ast$ are
canonical projections. On the other hand, note that for $F,K: T^\ast
V \cong V \times V^\ast \to \mathbb{R}$, by using the canonical
symplectic form $\omega_V$ on $T^\ast V $, we
can define the Poisson bracket $\{\cdot,\cdot\}_V$ on $T^\ast V$ as
follows
$$\{F,K\}_V(\theta,\lambda)=<\frac{\delta F}{\delta \theta},
\frac{\delta K}{\delta \lambda}>- <\frac{\delta K}{\delta
\theta},\frac{\delta F}{\delta \lambda}>
$$
If $\theta_i, \; i=1,\cdots, k,$ is a base of $V$, and $\lambda_i,
\; i=1,\cdots, k,$ a base of $V^\ast$, then we have that
\begin{equation}
\{F,K\}_V(\theta,\lambda)=\sum_{i=1}^k(\frac{\partial F}{\partial \theta_i}
\frac{\partial K}{\partial \lambda_i}- \frac{\partial K}{\partial
\theta_i}\frac{\partial F}{\partial \lambda_i}).
\label{3.7}
\end{equation}
Thus, by using the (-)-semidirect product Lie-Poisson brackets on
$\mathfrak{w}^\ast$ and the Poisson bracket $\{\cdot,\cdot\}_V$ on
$T^\ast V$, for $F,K: \mathfrak{w}^\ast \times V \times V^\ast \to
\mathbb{R}$, we can define the Poisson bracket on $\mathfrak{w}^\ast
\times V \times V^\ast$ as follows
\begin{align} & \{F,K\}_{-}(\mu,a,\theta,\lambda)= \{F,K\}_{\mathfrak{w}^*}(\mu,a)+
\{F,K\}_V(\theta,\lambda)\nonumber \\ & = -<\mu,[\frac{\delta
F}{\delta \mu},\frac{\delta K}{\delta\mu}]> - \langle a,\frac{\delta
F}{\delta \mu}\cdot\frac{\delta K}{\delta a}- \frac{\delta K}{\delta
\mu}\cdot\frac{\delta F}{\delta a}\rangle +
 \sum_{i=1}^k(\frac{\partial F}{\partial \theta_i} \frac{\partial
K}{\partial \lambda_i}- \frac{\partial K}{\partial
\theta_i}\frac{\partial F}{\partial \lambda_i}).
\label{3.8}
\end{align}
See Krishnaprasad and Marsden \cite{krma87}. In particular, for
$F_{(\mu,a)},K_{(\mu,a)}: \mathcal{O}_{(\mu,a)} \times V \times
V^\ast \to \mathbb{R}$,  the symplectic form
$\tilde{\omega}^{-}_{\mathcal{O}_{(\mu,a)} \times
V \times V^\ast}$ is induced
from the above Poisson bracket on $\mathfrak{w}^\ast
\times V \times V^\ast$, that is,
$\tilde{\omega}_{\mathcal{O}_{(\mu,a)} \times V \times
V^\ast}^{-}(X_{F_{(\mu,a)}}, X_{K_{(\mu,a)}})=
\{F_{(\mu,a)},K_{(\mu,a)}\}_{-}|_{\mathcal{O}_{(\mu,a)} \times V \times V^\ast}, $
such that $(\mathcal{O}_{(\mu,a)} \times
V\times V^\ast,\tilde{\omega}_{\mathcal{O}_{(\mu,a)} \times V \times
V^\ast}^{-})$ is a symplectic leaf of the Poisson manifold
$(\mathfrak{w}^\ast \times V \times V^\ast, \{\cdot,\cdot\}_{-}). $ \\

On the other hand, from $T^\ast Q = T^\ast W \times T^\ast V$ we
know that there is a canonical symplectic form $\omega_Q= \pi^\ast_1
\omega_0 +\pi^\ast_2 \omega_V$ on $T^\ast Q$, where $\omega_0$ is
the canonical symplectic form on $T^\ast W$ and the maps $\pi_1:
Q=W\times V \to W$ and $\pi_2: Q=W\times V \to V$ are canonical
projections. Assume that the cotangent lift of the left $W$-action
$\Phi^{T^*}: W \times T^\ast Q \to T^\ast Q$ is also symplectic, and the action
admits an associated $\operatorname{Ad}^\ast$-equivariant momentum
map $\mathbf{J}_Q: T^\ast Q \to \mathfrak{w}^\ast$ such that
$\mathbf{J}_Q\cdot \pi^\ast_1=\mathbf{J}_W$, where $\mathbf{J}_W:
T^\ast W \rightarrow \mathfrak{w}^\ast$ is a momentum map of left
$W$-action on $T^\ast W$ and we assume that it exists,
and $\pi^\ast_1: T^\ast W \to T^\ast Q$.
If $(\mu,a)\in\mathfrak{w}^\ast $ is a regular value of
$\mathbf{J}_Q$, then $(\mu,a)\in\mathfrak{w}^\ast $ is also a
regular value of $\mathbf{J}_W $ and $\mathbf{J}_Q^{-1}(\mu,a)\cong
\mathbf{J}_W^{-1}(\mu,a)\times V \times V^\ast $. Denote by
$W_{(\mu,a)}=\{(g,x)\in W|\operatorname{Ad}_{(g,x)^{-1}}^\ast
(\mu,a)= (\mu,a) \},$ the isotropy subgroup of co-adjoint $W$-action
at the point $(\mu,a)\in \mathfrak{w}^\ast$. It follows that
$W_{(\mu,a)}$ acts also freely and properly on
$\mathbf{J}_Q^{-1}(\mu,a)$, the regular point reduced space $(T^\ast
Q)_{(\mu,a)}=\mathbf{J}_Q^{-1}(\mu,a)/W_{(\mu,a)}\cong (T^\ast
W)_{(\mu,a)} \times V \times V^\ast$ of $(T^\ast Q,\omega_Q)$ at
$(\mu,a)$, is a symplectic manifold with symplectic form
$\omega_{(\mu,a)}$ uniquely characterized by the relation
$\pi_{(\mu,a)}^\ast \omega_{(\mu,a)}=i_{(\mu,a)}^\ast
\omega_Q=i_{(\mu,a)}^\ast \pi^\ast_1 \omega_0 +i_{(\mu,a)}^\ast
\pi^\ast_2 \omega_V$, where the map
$i_{(\mu,a)}:\mathbf{J}_Q^{-1}(\mu,a)\rightarrow T^\ast Q$ is the
inclusion and $\pi_{(\mu,a)}:\mathbf{J}_Q^{-1}(\mu,a)\rightarrow
(T^\ast Q)_{(\mu,a)}$ is the projection. From the reduction by
stages propositions for semidirect product Lie group, see
Proposition 3.1 and Proposition 3.2, we know that $((T^\ast
W)_{(\mu,a)},\omega_{(\mu,a)})$ is symplectically diffeomorphic to
$(\mathcal{O}_{(\mu,a)},\omega_{\mathcal{O}_{(\mu,a)}}^{-})$, and
hence we have that $((T^\ast Q)_{(\mu,a)},\omega_{(\mu,a)})$ is
symplectically diffeomorphic to $(\mathcal{O}_{(\mu,a)} \times
V\times V^\ast,\tilde{\omega}_{\mathcal{O}_{(\mu,a)} \times V \times
V^\ast}^{-})$, which is a symplectic leaf of the Poisson manifold
$(\mathfrak{w}^\ast \times V \times V^\ast, \{\cdot,\cdot\}_{-}). $\\

In the following we shall give the regular point reduction of the RCH system
$(T^\ast Q,W,\omega_Q,H,F,\mathcal{C})$ with symmetry and momentum map
on the generalization of a semidirect product Lie group.
We identify $TW$ and $W\times \mathfrak{w}$ locally,
by using the left translation, and $TV\cong V\times V$,
then locally, $TQ\cong W\times \mathfrak{w} \times V\times V$.
In consequence,
we consider Lagrangian $L(g,x,\xi,v,\theta,\dot{\theta}):TQ
\cong W\times \mathfrak{w} \times V\times V \to \mathbb{R}$, which is usual
the total kinetic minus potential energy of the system, where $(g,x)
\in W, \; (\xi,v) \in \mathfrak{w}$, and $\theta \in V$,
$(\xi^i,v^i)$ and $\dot{\theta}^j=\frac{\mathrm{d}
\theta^j}{\mathrm{d} t}$, ($i=1,\cdots,n, \; j=1,\cdots,k$, $n=\dim
W$, $k=\dim V$), regarded as the velocity variables of the system. If we introduce
the conjugate momentum $\zeta_i=\frac{\partial L}{\partial \xi^i},
\; p_i=\frac{\partial L}{\partial v^i}, \; l_j=\frac{\partial
L}{\partial \dot{\theta}^j}$, $i=1,\cdots,n, \; j=1,\cdots,k,$ and
by the Legendre transformation
\begin{align*} FL: TQ \cong W\times
\mathfrak{w}\times V \times V & \to T^\ast Q \cong W\times
\mathfrak{w}^\ast \times V \times V^\ast, \\
(g^i, x^i, \xi^i, v^i, \theta^j, \dot{\theta}^j) & \to
(g^i, x^i, \zeta_i, p_i, \theta^j, l_j),
\end{align*} we have the Hamiltonian
$H(g,x,\zeta,p,\theta,l):T^\ast Q \cong W\times \mathfrak{w}^\ast
\times V \times V^\ast \to \mathbb{R}$ given by
\begin{equation}
H(g^i,x^i,\zeta_i,p_i,\theta^j,l_j)=\sum_{i=1}^{n}(\zeta_i\xi^i+
p_iv^i)+\sum_{j=1}^{k}l_j\dot{\theta}^j
-L(g^i,x^i,\xi^i,v^i,\theta^j,\dot{\theta}^j).\label{3.9}
\end{equation}

If Hamiltonian $H(g,x,\zeta,p,\theta,l): T^\ast Q \cong W\times
\mathfrak{w}^\ast \times V \times V^\ast \to \mathbb{R}$ is left
cotangent lifted $W$-action $\Phi^{T^*}$ invariant, for a regular
value $(\mu,a)\in \mathfrak{w}^\ast$ of $\mathbf{J}_Q$,
we have the associated $R_p$-reduced
Hamiltonian $h_{(\mu,a)}(\zeta,p,\theta,l): (T^\ast Q)_{(\mu,a)}
\cong\mathcal{O}_{(\mu,a)} \times V \times V^\ast \to \mathbb{R}$,
defined by $h_{(\mu,a)}\cdot \pi_{(\mu,a)}=H\cdot i_{(\mu,a)}$, and
the $R_p$-reduced Hamiltonian vector field $X_{h_{(\mu,a)}}$ given by
$$X_{h_{(\mu,a)}}(K_{(\mu,a)})=\{K_{(\mu,a)},h_{(\mu,a)}\}_{-}|_{\mathcal{O}_{(\mu,a)}
\times V \times V^\ast}, $$ where
$K_{(\mu,a)}(\zeta,p,\theta,l): \mathcal{O}_{(\mu,a)}
\times V\times V^* \to \mathbb{R}.$ \\

Thus, if we consider that the fiber-preserving map $F:
T^\ast Q \to T^\ast Q$ and the fiber submanifold $\mathcal{C}$ of
$T^\ast Q$ are all left cotangent lifted $W$-action $\Phi^{T^*}$
invariant, then for $u \in \mathcal{C}$, the 6-tuple $(T^\ast Q,W,\omega_Q,H,F,\mathcal{C})$
is a regular point reducible RCH system with a control law $u$,
and its dynamical vector field can be expressed by
\begin{align}
\tilde{X}= X_{(T^\ast Q,W,\omega_Q,H,F,u)}= X_H+ \textnormal{vlift}(F)+ \textnormal{vlift}(u),
\label{3.10}
\end{align}
where $\textnormal{vlift}(F)= \textnormal{vlift}(F)\cdot X_H, \;\;
\textnormal{vlift}(u)= \textnormal{vlift}(u)\cdot X_H $ are
the changes of $X_H$ under the actions of the external force $F$ and the control law $u$.\\

For the point $(\mu,a)\in\mathfrak{w}^\ast$, the regular value of the momentum map
$\mathbf{J}_Q: T^\ast Q \rightarrow \mathfrak{w}^\ast$,
assume that $F(\mathbf{J}_Q^{-1}(\mu,a))\subset \mathbf{J}_Q^{-1}(\mu,a) $,
$f_{(\mu,a)}\cdot \pi_{(\mu,a)}=\pi_{(\mu,a)}\cdot
F\cdot i_{(\mu,a)}$, $u \in (\mathcal{C}\cap \mathbf{J}_Q^{-1}(\mu,a))\neq \emptyset$,
and $u_{(\mu,a)} \in \mathcal{C}_{(\mu,a)}
=\pi_{(\mu,a)}(\mathcal{C}\cap \mathbf{J}_Q^{-1}(\mu,a))\subset
\mathcal{O}_{(\mu,a)} \times V \times V^\ast$, $u_{(\mu,a)}\cdot
\pi_{(\mu,a)}=\pi_{(\mu,a)}\cdot u\cdot i_{(\mu,a)}$, then the
$R_p$-reduced RCH system is the 5-tuple $(\mathcal{O}_{(\mu,a)} \times V\times
V^\ast,\tilde{\omega}_{\mathcal{O}_{(\mu,a)} \times V \times V^\ast
}^{-},h_{(\mu,a)},f_{(\mu,a)},u_{(\mu,a)})$, where
$\mathcal{O}_{(\mu,a)} \subset \mathfrak{w}^\ast$ is the co-adjoint
orbit, $\tilde{\omega}_{\mathcal{O}_{(\mu,a)} \times V \times V^\ast
}^{-}$ is the induced symplectic form on $\mathcal{O}_{(\mu,a)} \times
V\times V^\ast $, such that the Hamiltonian vector field
$$X_{h_{(\mu,a)}}(K_{(\mu,a)})=\tilde{\omega}_{\mathcal{O}_{(\mu,a)} \times V\times
V^{*}}^{-}(X_{K_{(\mu,a)}}, X_{h_{(\mu,a)}})
=\{K_{(\mu,a)},h_{(\mu,a)}\}_{-}|_{\mathcal{O}_{(\mu,a)}
\times V\times V^{*} }, $$
since $(\mathcal{O}_{(\mu,a)} \times
V \times V^*,\tilde{\omega}_{\mathcal{O}_{(\mu,a)}
\times V \times V^*}^{-})$ is a symplectic leaf of
the Poisson manifold
$(\mathfrak{w}^\ast \times V \times V^*, \{\cdot,\cdot\}_{-}). $. Moreover, assume that the
dynamical vector field of the $R_p$-reduced RCH system can be expressed by
\begin{equation}X_{(\mathcal{O}_{(\mu,a)} \times V\times
V^\ast,\tilde{\omega}_{\mathcal{O}_{(\mu,a)} \times V \times V^\ast
}^{-},h_{(\mu,a)},f_{(\mu,a)},u_{(\mu,a)})} =X_{h_{(\mu,a)}}+
\textnormal{vlift}(f_{(\mu,a)})+\textnormal{vlift}(u_{(\mu,a)}), \;
\label{3.11}
\end{equation}
where $X_{h_{(\mu,a)}} \in T(\mathcal{O}_{(\mu,a)} \times V\times
V^*) $ is Hamiltonian vector field of the $R_p$-reduced Hamiltonian
$h_{(\mu,a)}: \mathcal{O}_{(\mu,a)} \times V \times V^\ast \to
\mathbb{R}$, and $\textnormal{vlift}(f_{(\mu,a)})=
\textnormal{vlift}(f_{(\mu,a)})X_{h_{(\mu,a)}} \in
T(\mathcal{O}_{(\mu,a)} \times V\times V^*)$, and
$\textnormal{vlift}(u_{(\mu,a)})=
\textnormal{vlift}(u_{(\mu,a)})X_{h_{(\mu,a)}} \in
T(\mathcal{O}_{(\mu,a)} \times V\times V^*)$
are the changes of $X_{h_{(\mu,a)}}$ under the actions of the
$R_p$-reduced external force $f_{(\mu,a)}$
and the $R_p$-reduced control torque $u_{(\mu,a)}$.
Note that $\mbox{vlift}(u_{(\mu,a)})X_{h_{(\mu,a)}}$ is the vertical
lift of vector field $X_{h_{(\mu,a)}}$ under the action of
$u_{(\mu,a)}$ along fibers, that is,
\begin{align*}
 \textnormal{vlift}(u_{(\mu,a)})X_{h_{(\mu,a)}}(\zeta,p,\theta,l) & =
 \textnormal{vlift}((Tu_{(\mu,a)}
X_{h_{(\mu,a)}})(u_{(\mu,a)}(\zeta,p,\theta,l)),
(\zeta,p,\theta,l))\\ & = (Tu_{(\mu,a)}
X_{h_{(\mu,a)}})^v_{\tilde{\sigma}} (\zeta,p,\theta,l),
\end{align*}
where $\tilde{\sigma}$ is a geodesic in $\mathcal{O}_{(\mu,a)}
\times V\times V^*$ connecting $u_{(\mu,a)}(\zeta,p,\theta,l)$ and
$(\zeta,p,\theta,l)$, and \\ $(Tu_{(\mu,a)}
X_{h_{(\mu,a)}})^v_{\tilde{\sigma}} (\zeta,p,\theta,l) $ is the
parallel displacement of vertical vector $(Tu_{(\mu,a)}
X_{h_{(\mu,a)}})^v (\zeta,p,\theta,l) $ along the geodesic
$\tilde{\sigma}$ from $u_{(\mu,a)}(\zeta,p,\theta,l)$ to
$(\zeta,p,\theta,l)$, and $ \textnormal{vlift}(f_{(\mu,a)})X_{h_{(\mu,a)}}$
is defined in the similar manner.
The dynamical vector fields of the RCH system and the $R_p$-reduced
RCH system satisfy the condition
\begin{equation}X_{(\mathcal{O}_{(\mu,a)} \times V\times
V^\ast,\tilde{\omega}_{\mathcal{O}_{(\mu,a)} \times V \times V^\ast
}^{-},h_{(\mu,a)},f_{(\mu,a)},u_{(\mu,a)})}\cdot
\pi_{(\mu,a)}=T\pi_{(\mu,a)}\cdot X_{(T^\ast
Q,W,\omega_Q,H,F,u)}\cdot i_{(\mu,a)}. \label{3.12}
\end{equation}
See Marsden et al. \cite{mawazh10}, Wang \cite{wa18}
and Wang \cite{wa13d}. Thus, we can get the following theorem.

\begin{theo}
If the 6-tuple $(T^\ast Q,W,\omega_Q,H,F,\mathcal{C})$ is a regular
point reducible RCH system on the generalization of a
semidirect product Lie group $Q=W\times V$, where $W=G\circledS E$
is a semidirect product Lie group with Lie algebra
$\mathfrak{w}=\mathfrak{g}\circledS E$, $G$ is a Lie group with Lie
algebra $\mathfrak{g}$, $E$ is a $r$-dimensional vector space and
$V$ is a $k$-dimensional vector space, and the Hamiltonian $H:T^\ast
Q \to \mathbb{R}$, the fiber-preserving map $F: T^\ast Q \to T^\ast
Q$ and the fiber submanifold $\mathcal{C}$ of $T^\ast Q$ are all
left cotangent lifted $W$-action $\Phi^{T^*}$ invariant, then for the
point $(\mu,a)\in\mathfrak{w}^\ast$, the regular value of the
momentum map $\mathbf{J}_Q: T^\ast Q \rightarrow \mathfrak{w}^\ast$,
the $R_p$-reduced RCH
system is the 5-tuple $(\mathcal{O}_{(\mu,a)} \times V\times
V^\ast,\tilde{\omega}_{\mathcal{O}_{(\mu,a)} \times V \times V^\ast
}^{-},h_{(\mu,a)},f_{(\mu,a)},u_{(\mu,a)})$, where
$\mathcal{O}_{(\mu,a)} \subset \mathfrak{w}^\ast$ is the co-adjoint
orbit, $\tilde{\omega}_{\mathcal{O}_{(\mu,a)} \times V \times V^\ast
}^{-}$ is the induced symplectic form on $\mathcal{O}_{(\mu,a)} \times
V\times V^\ast $, $h_{(\mu,a)}\cdot \pi_{(\mu,a)}=H\cdot
i_{(\mu,a)}$, $F(\mathbf{J}_Q^{-1}(\mu,a))\subset \mathbf{J}_Q^{-1}(\mu,a) $,
$f_{(\mu,a)}\cdot \pi_{(\mu,a)}=\pi_{(\mu,a)}\cdot
F\cdot i_{(\mu,a)}$, $u \in (\mathcal{C}\cap \mathbf{J}_Q^{-1}(\mu,a))\neq \emptyset$,
and $u_{(\mu,a)} \in \mathcal{C}_{(\mu,a)}
=\pi_{(\mu,a)}(\mathcal{C}\cap \mathbf{J}_Q^{-1}(\mu,a))\subset
\mathcal{O}_{(\mu,a)} \times V \times V^\ast$, $u_{(\mu,a)}\cdot
\pi_{(\mu,a)}=\pi_{(\mu,a)}\cdot u\cdot i_{(\mu,a)}$, and the
dynamical vector field of the $R_p$-reduced RCH system satisfies
$(3.11)$ and $(3.12)$.  \hskip 0.3cm $\blacksquare$
\end{theo}

In the following we regard the underwater vehicle-rotor system with the control
torque $u$ acting on the rotors as a regular point reducible RCH system on
the generalization of the semidirect product Lie groups
$(\textmd{SO}(3)\circledS \mathbb{R}^3)\times\mathbb{R}^2 $ and
$(\textmd{SE}(3)\circledS \mathbb{R}^3)\times \mathbb{R}^2$,
respectively, and give explicitly the $R_p$-reduced
controlled underwater vehicle-rotor systems by calculation in
detail.\\

\subsection{Hamilton-Jacobi Theorems }

The Hamilton-Jacobi theory for the regular (controlled) Hamiltonian system
is a very important subject, following the theoretical and applied development of
geometric mechanics, a lot of important problems about this subject
are being explored and studied, see Wang \cite{wa17}, Wang \cite{wa13d}
Wang \cite{wa20a}and de Le\'{o}n and Wang \cite{lewa15}.
In this subsection, we shall give the geometric constraint conditions
of the canonical symplectic form and the $R_p$-reduced symplectic form for the
dynamical vector field of the regular point reducible RCH system
$(T^\ast Q,W,\omega_Q,H,F,\mathcal{C})$ on the
generalization of a semidirect product Lie group, that is, the two types of Hamilton-Jacobi equations for
the regular point reducible RCH system and the $R_p$- reduced RCH system.
We shall follow the notations
and conventions introduced in Marsden \cite{ma92}, Marsden et al \cite{mawazh10},
Wang \cite{wa17} and Wang \cite{wa13d}.\\

In the following we first give
an important notion and a key lemma, which is an important
tool for the proofs of two types of Hamilton-Jacobi theorems for
the regular point reducible RCH system and the $R_p$-reduced RCH system.\\

Denote by $\Omega^i(Q)$ the set of all i-forms on $Q$, $i=1,2.$
For any $\gamma \in \Omega^1(Q),\; q\in Q, $ then $\gamma(q)\in T_q^*Q, $
and we can define a map $\gamma: Q \rightarrow T^*Q, \; q \rightarrow (q, \gamma(q)).$
Hence we say often that the map $\gamma: Q
\rightarrow T^*Q$ is an one-form on $Q$. If the one-form $\gamma$ is closed,
then $\mathbf{d}\gamma(x,y)=0, \; \forall\;
x, y \in TQ$. In the following we give a weaker notion.
\begin{defi}
The one-form $\gamma$ is called to be closed with respect to $T\pi_{Q}:
TT^* Q \rightarrow TQ, $ if for any $v, w \in TT^* Q, $ we have
$\mathbf{d}\gamma(T\pi_{Q}(v),T\pi_{Q}(w))=0. $
\end{defi}

From the above definition we know that, if $\gamma$ is a closed one-form,
then it must be closed with respect to $T\pi_{Q}: TT^* Q \rightarrow
TQ. $ Conversely, if $\gamma$ is closed with respect to
$T\pi_{Q}: TT^* Q \rightarrow TQ, $ then it may not be closed.
Now, we give the following Lemma 3.5, which  is very important tool for our research,
see Wang \cite{wa17} and Wang \cite{wa13d}.
\begin{lemm}
Assume that $\gamma: Q \rightarrow T^*Q$ is an one-form on $Q$, and
$\lambda=\gamma \cdot \pi_{Q}: T^* Q \rightarrow T^* Q .$ Then
we have that the following two assertions hold.\\
\noindent $(\mathrm{i})$ For any $x, y \in TQ, \;
\gamma^*\omega(x,y)= -\mathbf{d}\gamma (x,y),$ and for any $v, w \in
TT^* Q, \; \lambda^*\omega(v,w)=\\ -\mathbf{d}\gamma(T\pi_{Q}(v), \;
T\pi_{Q}(w)),$
since $\omega$ is the canonical symplectic form on $T^*Q$; \\
\noindent $(\mathrm{ii})$ For any $v, w \in TT^* Q, \;
\omega(T\lambda \cdot v,w)= \omega(v, w-T\lambda \cdot
w)-\mathbf{d}\gamma(T\pi_{Q}(v), \; T\pi_{Q}(w)). $
\end{lemm}

For a given RCH system $(T^\ast Q,W,\omega_Q,H,F,\mathcal{C})$
on the generalization of a semidirect product Lie group, by using
the above Lemma 3.5, we can give precisely the geometric constraint
conditions of the canonical symplectic form $\omega_Q$ for the dynamical
vector field $X_{(T^\ast Q,W,\omega_Q,H,F,u)}$ of the RCH system with a control law $u$,
that is, we can obtain the following two types of
Hamilton-Jacobi equation for the RCH system.
\begin{theo}
For the RCH system $(T^\ast Q,W,\omega_Q,H,F,\mathcal{C})$
on the generalization of a semidirect product Lie group,
assume that $\gamma: Q \rightarrow T^*Q$ is an one-form on $Q$, and
$\lambda=\gamma\cdot\pi_{Q}: T^* Q \rightarrow T^* Q $, and for any
symplectic map $\varepsilon: T^* Q \rightarrow T^* Q $,
denote by $\tilde{X}^\gamma = T\pi_{Q}\cdot \tilde{X} \cdot \gamma$,
and $\tilde{X}^\varepsilon = T\pi_{Q}\cdot \tilde{X} \cdot \varepsilon$,
where $\tilde{X}=X_{(T^\ast Q,W,\omega_Q,H,F,u)}$ is the dynamical vector field
of the RCH system $(T^*Q,W,\omega_Q,H,F,\mathcal{C})$ with a control law $u$.
Then the following two assertions hold:\\
\noindent $(\mathbf{i})$
If the one-form $\gamma: Q \rightarrow T^*Q $ is closed with respect to $T\pi_{Q}:
TT^* Q \rightarrow TQ, $ then $\gamma$ is a solution of the
Type I of Hamilton-Jacobi equation
$T\gamma\cdot \tilde{X}^\gamma= X_H\cdot \gamma ,$ where $X_H$
is the Hamiltonian vector field
of the corresponding Hamiltonian system $(T^*Q,W,\omega_Q,H).$ \\
\noindent $(\mathbf{ii})$
The $\varepsilon$ is a solution of the Type II of
Hamilton-Jacobi equation $T\gamma\cdot \tilde{X}^\varepsilon= X_H\cdot
\varepsilon ,$ if and only if it is a solution of the equation
$T\varepsilon\cdot X_{H\cdot\varepsilon}= T\lambda \cdot \tilde{X} \cdot \varepsilon,$
where $X_H$ and $ X_{H\cdot\varepsilon} \in
TT^*Q $ are the Hamiltonian vector fields of the functions $H$ and $H\cdot\varepsilon:
T^*Q\rightarrow \mathbb{R}, $ respectively.
Here the maps involved in the theorem and its proof are shown
in the following Diagram-1.
\begin{center}
\hskip 0cm \xymatrix{ T^* Q \ar[r]^\varepsilon & T^* Q
\ar[d]_{X_{H\cdot \varepsilon}} \ar[dr]^{\tilde{X}^\varepsilon} \ar[r]^{\pi_Q}
& Q \ar[d]^{\tilde{X}^\gamma} \ar[r]^{\gamma}
& T^*Q \ar[d]^{\tilde{X}} \\
& T(T^*Q)  & TQ \ar[l]^{T\gamma}
& T(T^*Q) \ar[l]^{T\pi_Q} }
\end{center}
$$\mbox{Diagram-1}$$
\end{theo}
\noindent{\bf Proof: } Since $\tilde{X}=X_{(T^\ast Q,W,\omega_Q,H,F,u)}=X_H
+\textnormal{vlift}(F)+\textnormal{vlift}(u), $ and
$T\pi_{Q}\cdot \textnormal{vlift}(F)=T\pi_{Q}\cdot \textnormal{vlift}(u)=0, $
then we have that $T\pi_{Q}\cdot \tilde{X}\cdot \gamma=T\pi_{Q}\cdot X_H\cdot \gamma. $
If we take that $v= X_H\cdot \gamma \in TT^* Q, $ and for
any $w \in TT^* Q, \; T\pi_{Q}(w)\neq 0, $ from Lemma 3.5(ii) we have that
\begin{align*}
& \omega_Q(T\gamma \cdot \tilde{X}^\gamma, \; w)= \omega_Q(T\gamma \cdot
T\pi_Q\cdot \tilde{X}\cdot \gamma, \; w)= \omega_Q(T\gamma \cdot
T\pi_Q\cdot X_H\cdot \gamma, \; w)\\ &=\omega_Q(T(\gamma \cdot
\pi_Q)\cdot X_H\cdot \gamma, \; w)= \omega_Q(X_H\cdot \gamma, \;
w-T(\gamma \cdot \pi_Q)\cdot
w)-\mathbf{d}\gamma(T\pi_{Q}(X_H\cdot \gamma), \; T\pi_{Q}(w))\\
& =\omega_Q(X_H\cdot \gamma, \; w) - \omega_Q(X_H\cdot \gamma, \;
T\lambda \cdot w)-\mathbf{d}\gamma(T\pi_{Q}(X_H\cdot \gamma), \; T\pi_{Q}(w)).
\end{align*}
Because the one-form $\gamma: Q \rightarrow T^*Q $ is closed with respect to
$T\pi_Q: TT^* Q \rightarrow TQ, $ then we have that
$$
\mathbf{d}\gamma(T\pi_{Q}(X_H\cdot \gamma), \; T\pi_{Q}(w))=0,
$$
and hence
\begin{equation}
\omega_Q(T\gamma \cdot \tilde{X}^\gamma, \; w)- \omega_Q(X_H\cdot \gamma, \; w)
= -\omega_Q(X_H\cdot \gamma, \; T\lambda \cdot w). \;\; \label{3.13}
\end{equation}
If $\gamma$ satisfies the equation $T\gamma\cdot \tilde{X}^\gamma= X_H\cdot \gamma ,$
from Lemma 3.5(i) we can obtain that
\begin{align*}
-\omega_Q(X_H\cdot \gamma, \; T\lambda \cdot w) &
= -\omega_Q(T\gamma \cdot \tilde{X}^\gamma, \; T\lambda \cdot w)\\
& =-\omega_Q(T\gamma \cdot T\pi_{Q} \cdot \tilde{X}\cdot\gamma, \; T\lambda \cdot w)
=-\omega_Q(T\lambda \cdot \tilde{X}\cdot\gamma, \; T\lambda \cdot w)\\
& = -\lambda^*\omega_Q( \tilde{X}\cdot\gamma, \; w)=
\mathbf{d}\gamma(T\pi_{Q}( \tilde{X}\cdot\gamma ), \; T\pi_{Q}(w))=0,
\end{align*}
since the one-form $\gamma: Q \rightarrow T^*Q $ is closed with respect to
$T\pi_Q: TT^* Q \rightarrow TQ. $ But, because the symplectic form $\omega_Q$ is non-degenerate,
the left side of (3.13) equals zero, only when
$\gamma$ satisfies the equation $T\gamma\cdot \tilde{X}^\gamma= X_H\cdot \gamma .$ Thus,
if the one-form $\gamma: Q \rightarrow T^*Q $ is closed with respect to
$T\pi_Q: TT^* Q \rightarrow TQ, $ then $\gamma$ must be a solution of
the Type I of Hamilton-Jacobi equation
$T\gamma\cdot \tilde{X}^\gamma= X_H\cdot \gamma .$\\

Next, we prove the assertion $(\mathrm{ii})$.
For any symplectic map $\varepsilon: T^* Q \rightarrow T^* Q $,
if we take that $v= X_H\cdot \varepsilon \in TT^* Q, $ and for
any $w \in TT^* Q, \; T\lambda(w)\neq 0, $ from Lemma 3.5 we have that
\begin{align*}
&\omega_Q(T\gamma \cdot \tilde{X}^\varepsilon, \; w)
= \omega_Q(T\gamma \cdot
T\pi_Q\cdot \tilde{X}\cdot \varepsilon, \; w)= \omega_Q(T\gamma \cdot
T\pi_Q\cdot X_H\cdot \varepsilon, \; w)\\ &
= \omega_Q(T(\gamma \cdot
\pi_Q)\cdot X_H\cdot \varepsilon, \; w)= \omega_Q(X_H\cdot \varepsilon, \;
w-T(\gamma \cdot \pi_Q)\cdot
w)-\mathbf{d}\gamma(T\pi_{Q}(X_H\cdot \varepsilon), \; T\pi_{Q}(w))\\
& =\omega_Q(X_H\cdot \varepsilon, \; w) - \omega_Q(X_H\cdot \varepsilon, \;
T\lambda \cdot w)+\lambda^*\omega_Q(X_H\cdot \varepsilon, \; w)\\
& =\omega_Q(X_H\cdot \varepsilon, \; w) - \omega_Q(X_H\cdot \varepsilon, \;
T\lambda \cdot w)+ \omega_Q(T\lambda \cdot X_H\cdot \varepsilon, \; T\lambda \cdot w).
\end{align*}
Because $\varepsilon: T^* Q
\rightarrow T^* Q $ is symplectic, and hence $ X_H\cdot \varepsilon=
T\varepsilon \cdot X_{H\cdot\varepsilon}, $ along $\varepsilon$. Note that
$T\lambda \cdot X_H\cdot \varepsilon=T\gamma \cdot
T\pi_Q\cdot X_H\cdot \varepsilon=T\gamma \cdot
T\pi_Q\cdot \tilde{X}\cdot \varepsilon=T\lambda\cdot \tilde{X}\cdot \varepsilon.$
From the above arguments, we can obtain that
\begin{align*}
&\omega_Q(T\gamma \cdot \tilde{X}^\varepsilon, \; w)- \omega_Q(X_H\cdot \varepsilon, \; w)\\
& =- \omega_Q(X_H\cdot \varepsilon, \;
T\lambda \cdot w)+ \omega_Q(T\lambda \cdot X_H\cdot \varepsilon, \; T\lambda \cdot w)\\
& =-\omega_Q(T\varepsilon \cdot X_{H\cdot\varepsilon}, \; T\lambda \cdot w)
+ \omega_Q(T\lambda \cdot \tilde{X}\cdot \varepsilon, \; T\lambda \cdot w)\\
& = \omega_Q(T\lambda \cdot \tilde{X}\cdot \varepsilon -T\varepsilon \cdot X_{H\cdot\varepsilon}, \; T\lambda \cdot w).
\end{align*}
Because the symplectic form $\omega_Q$ is non-degenerate,
it follows that $T\gamma\cdot \tilde{X}^\varepsilon= X_H\cdot
\varepsilon ,$ is equivalent to $T\varepsilon \cdot X_{H\cdot\varepsilon} = T\lambda\cdot \tilde{X}\cdot \varepsilon $.
Thus, $\varepsilon$ is a solution of the equation
$T\varepsilon\cdot X_{H\cdot\varepsilon}= T\lambda \cdot \tilde{X} \cdot\varepsilon,$
if and only if it is a solution of the Type II of Hamilton-Jacobi equation
$T\gamma\cdot \tilde{X}^\varepsilon= X_H\cdot \varepsilon .$
\hskip 0.3cm $\blacksquare$\\

Moreover, we can give precisely the geometric constraint
conditions of the $R_p$-reduced symplectic form
$\tilde{\omega}_{\mathcal{O}_{(\mu,a)} \times V\times
V*}^{-}$ for the dynamical vector field
of the regular point reducible RCH system
$(T^\ast Q,W,\omega_Q,H,F,\mathcal{C})$ on the
generalization of a semidirect product Lie group,
that is, Type I and Type II of
Hamilton-Jacobi equation for the $R_p$-reduced RCH system
$(\mathcal{O}_{(\mu,a)} \times V\times
V^\ast,\tilde{\omega}_{\mathcal{O}_{(\mu,a)} \times V \times V^\ast
}^{-},\\ h_{(\mu,a)},f_{(\mu,a)},u_{(\mu,a)})$.
For convenience, the maps involved in
the following theorem and its proof are shown in Diagram-2.
\begin{center}
\hskip 0cm \xymatrix{ \mathbf{J}^{-1}(\mu,a) \ar[r]^{i_{(\mu,a)}} & T^* Q
\ar[d]_{X_{H\cdot \varepsilon}} \ar[dr]^{\tilde{X}^\varepsilon} \ar[r]^{\pi_Q}
& Q \ar[d]^{\tilde{X}^\gamma} \ar[r]^{\gamma}
& T^*Q \ar[d]_{\tilde{X}} \ar[dr]_{X_{h_{(\mu,a)} \cdot\bar{\varepsilon}}} \ar[r]^{\pi_{(\mu,a)}}
&\;\;\; \mathcal{O}_{(\mu,a)} \times V \times V^\ast \ar[d]^{X_{h_{(\mu,a)}}} \\
& T(T^*Q)  & TQ \ar[l]^{T\gamma}
& T(T^*Q) \ar[l]^{T\pi_Q} \ar[r]_{T\pi_{(\mu,a)}} & \;\;\; T(\mathcal{O}_{(\mu,a)} \times V \times V^\ast) }
\end{center}
$$\mbox{Diagram-2}$$

\begin{theo}
If the 6-tuple $(T^\ast Q,W,\omega_Q,H,F,\mathcal{C})$ is a regular point reducible RCH system
on the generalization of a semidirect
product Lie group $Q=W\times V$, where $W=G\circledS E$ is a
semidirect product Lie group with Lie algebra
$\mathfrak{w}=\mathfrak{g}\circledS E$, $G$ is a Lie group with Lie
algebra $\mathfrak{g}$, $E$ is a $r$-dimensional vector space and
$V$ is a $k$-dimensional vector space, and for $(\mu,a) \in \mathfrak{w}^\ast $, the regular value of
the momentum map $\mathbf{J}_Q: T^\ast Q \rightarrow
\mathfrak{w}^\ast$,
the $R_p$-reduced system is the 5-tuple
$(\mathcal{O}_{(\mu,a)} \times V\times
V^\ast,\tilde{\omega}_{\mathcal{O}_{(\mu,a)} \times V \times V^\ast
}^{-},h_{(\mu,a)},f_{(\mu,a)},\mathcal{C}_{(\mu,a)})$.
Assume that $\gamma: Q \rightarrow T^* Q$ is an one-form on
$Q$, and $\lambda=\gamma \cdot \pi_{Q}: T^* Q \rightarrow T^* Q $, and
$\varepsilon: T^* Q \rightarrow T^* Q $ is a $W_{(\mu,a)}$-invariant
symplectic map, where $W_{(\mu,a)}$ is
the isotropy subgroup of co-adjoint $W$-action at the point $(\mu,a)$.
Denote $\tilde{X}^\gamma = T\pi_{Q}\cdot \tilde{X} \cdot \gamma$,
and $\tilde{X}^\varepsilon = T\pi_{Q}\cdot \tilde{X} \cdot \varepsilon$,
where $\tilde{X}=X_{(T^\ast Q,W,\omega_Q,H,F,u)}$ is the dynamical
vector field of the regular point reducible RCH system
$(T^*Q,W,\omega_Q,H,F,\mathcal{C})$ with a control law $u$. Moreover, assume
that $\textmd{Im}(\gamma)\subset
\mathbf{J}_Q^{-1}{(\mu,a)}, $ and it is $W_{(\mu,a)}$-invariant, and
$\varepsilon(\mathbf{J}_Q^{-1}{(\mu,a)}\subset \mathbf{J}_Q^{-1}{(\mu,a)}.$
Denote
$\bar{\gamma}=\pi_{(\mu,a)}(\gamma): Q \rightarrow \mathcal{O}_{(\mu,a)} \times
V\times V^\ast, $ and $\bar{\lambda}=\pi_{(\mu,a)}(\lambda):
\mathbf{J}_Q^{-1}{(\mu,a)} \rightarrow \mathcal{O}_{(\mu,a)}\times V\times
V^\ast, $ and
$\bar{\varepsilon}=\pi_{(\mu,a)}(\varepsilon):
\mathbf{J}_Q^{-1}{(\mu,a)}\rightarrow \mathcal{O}_{(\mu,a)}\times V\times
V^\ast. $
Then the following two assertions hold:\\
\noindent $(\mathbf{i})$
If the one-form $\gamma: Q \rightarrow T^*Q $ is closed with respect to
$T\pi_Q: TT^* Q \rightarrow TQ, $
then $\bar{\gamma}$ is a solution of the Type I of Hamilton-Jacobi equation
$T\bar{\gamma}\cdot \tilde{X}^\gamma= X_{h_{(\mu,a)}}\cdot \bar{\gamma}; $\\
\noindent $(\mathbf{ii})$
The $\varepsilon$ and $\bar{\varepsilon} $ satisfy the Type II of Hamilton-Jacobi equation
$T\bar{\gamma}\cdot \tilde{X}^\varepsilon= X_{h_{(\mu,a)}}\cdot \bar{\varepsilon}, $
if and only if they satisfy
the equation $T\bar{\varepsilon}\cdot(X_{h_{(\mu,a)} \cdot \bar{\varepsilon}})
= T\bar{\lambda}\cdot \tilde{X} \cdot\varepsilon. $
Here $X_{h_{(\mu,a)}}$ and $ X_{h_{(\mu,a)} \cdot \bar{\varepsilon}} \in TT^*Q $
are the Hamiltonian vector fields of the $R_p$-reduced Hamiltonian
functions $h_{(\mu,a)}$ and $h_{(\mu,a)} \cdot \bar{\varepsilon}: T^*Q\rightarrow
\mathbb{R}, $ respectively.
\end{theo}

\noindent{\bf Proof: } At first, from the above Theorem 3.6, we know that
$\gamma$ is a solution of the Type I of Hamilton-Jacobi equation
$T\gamma\cdot \tilde{X}^\gamma= X_H\cdot \gamma. $ Next, we note that
$\textmd{Im}(\gamma)\subset \mathbf{J}_Q^{-1}(\mu,a), $ and it
is $W_{(\mu,a)}$-invariant, in this case,
$\pi_{(\mu,a)}^*\tilde{\omega}_{\mathcal{O}_{(\mu,a)} \times V \times V^\ast}^{-}
= i_{(\mu,a)}^*\omega_Q= \omega_Q, $ along $\textmd{Im}(\gamma)$.
Since $\tilde{X}=X_{(T^\ast Q,W,\omega_Q,H,F,u)}=X_H
+\textnormal{vlift}(F)+\textnormal{vlift}(u), $ and
$T\pi_{Q}\cdot \textnormal{vlift}(F)=T\pi_{Q}\cdot \textnormal{vlift}(u)=0, $
then we have that $T\pi_{Q}\cdot \tilde{X}\cdot \gamma=T\pi_{Q}\cdot X_H\cdot \gamma. $
By using the $R_p$-reduced symplectic form
$\tilde{\omega}_{\mathcal{O}_{(\mu,a)} \times V\times
V*}^{-}$, if we take
that $v= X_H\cdot \gamma \in TT^* Q,$ and for any $w \in TT^* Q, \; T\pi_{Q}(w)\neq 0,$
and $T\pi_{(\mu,a)} (w) \neq 0, $ from Lemma 3.5(ii) we have that
\begin{align*}
& \tilde{\omega}_{\mathcal{O}_{(\mu,a)} \times V\times
V*}^{-}(T\bar{\gamma} \cdot \tilde{X}^\gamma, \; T\pi_{(\mu,a)} \cdot w) \\
 & =\tilde{\omega}_{\mathcal{O}_{(\mu,a)} \times V\times
V*}^{-}(T(\pi_{(\mu,a)} \cdot \gamma) \cdot T\pi_Q\cdot \tilde{X}\cdot \gamma, \; T\pi_{(\mu,a)} \cdot w )\\
& = \pi_{(\mu,a)}^*\tilde{\omega}_{\mathcal{O}_{(\mu,a)} \times V\times
V*}^{-}(T\gamma \cdot T\pi_Q\cdot X_H \cdot\gamma, \; w)
= \omega_Q(T(\gamma \cdot \pi_Q)\cdot X_H\cdot \gamma, \; w)\\
& = \omega_Q(X_H\cdot \gamma, \; w-T(\gamma \cdot \pi_Q)\cdot w)
-\mathbf{d}\gamma(T\pi_{Q}(X_H\cdot \gamma), \; T\pi_{Q}(w))\\
& =\omega_Q(X_H\cdot \gamma, \; w) - \omega_Q(X_H\cdot \gamma, \;
T(\gamma \cdot \pi_Q)\cdot w)-\mathbf{d}\gamma(T\pi_{Q}(X_H\cdot \gamma), \; T\pi_{Q}(w))\\
& =\pi_{(\mu,a)}^*\tilde{\omega}_{\mathcal{O}_{(\mu,a)} \times V\times
V*}^{-}(X_H\cdot \gamma, \; w)
- \pi_{(\mu,a)}^*\tilde{\omega}_{\mathcal{O}_{(\mu,a)} \times V\times
V*}^{-}(X_H\cdot \gamma, \; T(\gamma \cdot \pi_Q)\cdot w) \\
& \;\;\;\; -\mathbf{d}\gamma(T\pi_{Q}(X_H\cdot \gamma), \; T\pi_{Q}(w))\\
& = \tilde{\omega}_{\mathcal{O}_{(\mu,a)} \times V\times
V*}^{-}(T\pi_{(\mu,a)}(X_H\cdot \gamma), \;
T\pi_{(\mu,a)} \cdot w) \\
& \;\;\;\; - \tilde{\omega}_{\mathcal{O}_{(\mu,a)} \times V\times
V*}^{-}(T\pi_{(\mu,a)}\cdot(X_H\cdot \gamma), \;
T(\pi_{(\mu,a)} \cdot\gamma \cdot \pi_Q) \cdot w)
 -\mathbf{d}\gamma(T\pi_{Q}(X_H\cdot \gamma), \; T\pi_{Q}(w))\\
& = \tilde{\omega}_{\mathcal{O}_{(\mu,a)} \times V\times
V*}^{-}(T\pi_{(\mu,a)}(X_H)\cdot \pi_{(\mu,a)}(\gamma), \; T\pi_{(\mu,a)} \cdot w)\\
& \;\;\;\; - \tilde{\omega}_{\mathcal{O}_{(\mu,a)} \times V\times
V*}^{-}(T\pi_{(\mu,a)}(X_H)\cdot \pi_{(\mu,a)}(\gamma), \; T\bar{\gamma}\cdot T\pi_{Q}(w)) -\mathbf{d}\gamma(T\pi_{Q}(X_H\cdot \gamma), \; T\pi_{Q}(w))\\
& = \tilde{\omega}_{\mathcal{O}_{(\mu,a)} \times V\times
V*}^{-}(X_{h_{(\mu,a)}} \cdot
\bar{\gamma}, \; T\pi_{(\mu,a)} \cdot w) \\
& \;\;\;\; - \tilde{\omega}_{\mathcal{O}_{(\mu,a)} \times V\times
V*}^{-}(X_{h_{(\mu,a)}} \cdot
\bar{\gamma}, \; T\bar{\gamma} \cdot T\pi_{Q}(w))
-\mathbf{d}\gamma(T\pi_{Q}(X_H\cdot \gamma), \; T\pi_{Q}(w)),
\end{align*}
where we have used that $T\pi_{(\mu,a)}(X_H)= X_{h_{(\mu,a)}}. $
Since the one-form $\gamma: Q \rightarrow T^*Q $ is closed with respect to
$T\pi_Q: TT^* Q \rightarrow TQ, $ then we have that $
\mathbf{d}\gamma(T\pi_{Q}(X_H\cdot \gamma), \; T\pi_{Q}(w))=0, $
and hence
\begin{align}
& \tilde{\omega}_{\mathcal{O}_{(\mu,a)} \times V\times
V*}^{-}(T\bar{\gamma} \cdot X_H^\gamma, \; T\pi_{(\mu,a)} \cdot w)-
\tilde{\omega}_{\mathcal{O}_{(\mu,a)} \times V\times
V*}^{-}(X_{h_{(\mu,a)}} \cdot \bar{\gamma}, \; T\pi_{(\mu,a)} \cdot w) \nonumber \\
& = - \tilde{\omega}_{\mathcal{O}_{(\mu,a)} \times V\times
V*}^{-}(X_{h_{(\mu,a)}} \cdot
\bar{\gamma}, \; T\bar{\gamma} \cdot T\pi_{Q}(w)). \; \label{3.14}
\end{align}
If $\bar{\gamma}$ satisfies the equation
$T\bar{\gamma}\cdot \tilde{X}^\gamma= X_{h_{(\mu,a)}}\cdot \bar{\gamma}, $
from Lemma 3.5(i) we can obtain that
\begin{align*}
& - \tilde{\omega}_{\mathcal{O}_{(\mu,a)} \times V\times
V*}^{-}(X_{h_{(\mu,a)}} \cdot
\bar{\gamma}, \; T\bar{\gamma} \cdot T\pi_{Q}(w)) \\
& = -\tilde{\omega}_{\mathcal{O}_{(\mu,a)} \times V\times
V*}^{-} (T\bar{\gamma} \cdot \tilde{X}^\gamma, \; T\bar{\gamma} \cdot T\pi_{Q}(w))\\
& = -\bar{\gamma}^*\tilde{\omega}_{\mathcal{O}_{(\mu,a)} \times V\times
V*}^{-} (T\pi_{Q} \cdot \tilde{X}\cdot\gamma, \; T\pi_Q(w))\\
& = -\gamma^* \cdot \pi_{(\mu,a)}^*\tilde{\omega}_{\mathcal{O}_{(\mu,a)} \times V\times
V*}^{-} (T\pi_{Q} \cdot X_{H}\cdot\gamma, \; T\pi_{Q}(w))\\
& = -\gamma^*\omega_Q( T\pi_{Q}(X_{H}\cdot\gamma), \; T\pi_{Q}(w))\\
& = \mathbf{d}\gamma(T\pi_{Q}( X_{H}\cdot\gamma ), \; T\pi_{Q}(w))=0.
\end{align*}
Because the $R_p$-reduced symplectic form
$\tilde{\omega}_{\mathcal{O}_{(\mu,a)} \times V\times
V*}^{-}$ is non-degenerate, the left side of (3.14) equals zero, only when
$\bar{\gamma}$ satisfies the equation
$T\bar{\gamma}\cdot \tilde{X}^\gamma= X_{h_{(\mu,a)}}\cdot \bar{\gamma}.$ Thus,
if the one-form $\gamma: Q \rightarrow T^*Q $ is closed with respect to
$T\pi_Q: TT^* Q \rightarrow TQ, $ then $\bar{\gamma}$ must be
a solution of the Type I of Hamilton-Jacobi equation
$T\bar{\gamma}\cdot \tilde{X}^\gamma= X_{h_{(\mu,a)}}\cdot \bar{\gamma}.$\\

Next, we prove the assertion $(\mathrm{ii})$.
For any $W_{(\mu,a)}$-invariant symplectic map $\varepsilon: T^* Q \rightarrow T^* Q $,
using the $R_p$-reduced symplectic form $\tilde{\omega}_{\mathcal{O}_{(\mu,a)} \times V\times
V*}^{-}$, if we take
that $v= X_H\cdot \varepsilon \in TT^* Q,$ and for any $w \in TT^* Q, \; T\bar{\lambda}(w)\neq 0,$
and $T\pi_{(\mu,a)} (w) \neq 0, $ from Lemma 3.5(ii) we have that
\begin{align*}
& \tilde{\omega}_{\mathcal{O}_{(\mu,a)} \times V\times
V*}^{-}(T\bar{\gamma} \cdot \tilde{X}^\varepsilon, \; T\pi_{(\mu,a)} \cdot w) \\
& = \tilde{\omega}_{\mathcal{O}_{(\mu,a)} \times V\times
V*}^{-}(T(\pi_{(\mu,a)} \cdot \gamma) \cdot T\pi_Q \cdot \tilde{X}\cdot \varepsilon, \; T\pi_{(\mu,a)}
\cdot w )\\
& = \pi_{(\mu,a)}^*\tilde{\omega}_{\mathcal{O}_{(\mu,a)} \times V\times
V*}^{-}(T\gamma \cdot T\pi_Q \cdot X_H \cdot\varepsilon, \; w)
= \omega_Q(T(\gamma \cdot \pi_Q)\cdot X_H\cdot \varepsilon, \; w)\\
& = \omega_Q(X_H\cdot \varepsilon, \; w-T(\gamma \cdot \pi_Q)\cdot w)
-\mathbf{d}\gamma(T\pi_{Q}(X_H\cdot \varepsilon), \; T\pi_{Q}(w))\\
& =\omega_Q(X_H\cdot \varepsilon, \; w) - \omega_Q(X_H\cdot \varepsilon, \;
T\lambda\cdot w)-\mathbf{d}\gamma(T\pi_{Q}(\tilde{X}\cdot \varepsilon), \; T\pi_{Q}(w))\\
& =\pi_{(\mu,a)}^*\tilde{\omega}_{\mathcal{O}_{(\mu,a)} \times V\times
V*}^{-}(X_H\cdot \varepsilon, \; w)
- \pi_{(\mu,a)}^*\tilde{\omega}_{\mathcal{O}_{(\mu,a)} \times V\times
V*}^{-}(X_H\cdot \varepsilon, \; T\lambda\cdot w)+ \lambda^*\omega_Q(\tilde{X}\cdot \varepsilon, \; w)\\
& = \tilde{\omega}_{\mathcal{O}_{(\mu,a)} \times V\times
V*}^{-}(T\pi_{(\mu,a)}(X_H\cdot \varepsilon), \;
T\pi_{(\mu,a)} \cdot w) \\
& \;\;\;\; - \tilde{\omega}_{\mathcal{O}_{(\mu,a)} \times V\times
V*}^{-}(T\pi_{(\mu,a)}\cdot(X_H\cdot \varepsilon), \;
T(\pi_{(\mu,a)} \cdot\lambda) \cdot w)+ \lambda^*\cdot \pi_{(\mu,a)}^*\cdot \tilde{\omega}_{\mathcal{O}_{(\mu,a)} \times V\times
V*}^{-}(\tilde{X}\cdot \varepsilon, \; w)\\
& = \tilde{\omega}_{\mathcal{O}_{(\mu,a)} \times V\times
V*}^{-}(T\pi_{(\mu,a)}(X_H)\cdot \pi_{(\mu,a)}(\varepsilon), \; T\pi_{(\mu,a)} \cdot w)\\
& \;\;\;\; - \tilde{\omega}_{\mathcal{O}_{(\mu,a)} \times V\times
V*}^{-}(T\pi_{(\mu,a)}(X_H)\cdot \pi_{(\mu,a)}(\varepsilon), \; T\bar{\lambda}\cdot w)
+ (\pi_{(\mu,a)}\cdot\lambda)^*\cdot \tilde{\omega}_{\mathcal{O}_{(\mu,a)} \times V\times
V*}^{-}(\tilde{X}\cdot \varepsilon, \; w)\\
& = \tilde{\omega}_{\mathcal{O}_{(\mu,a)} \times V\times
V*}^{-}(X_{h_{(\mu,a)}} \cdot
\bar{\varepsilon}, \; T\pi_{(\mu,a)} \cdot w)\\
& \;\;\;\; - \tilde{\omega}_{\mathcal{O}_{(\mu,a)} \times V\times
V*}^{-}(X_{h_{(\mu,a)}} \cdot
\bar{\varepsilon}, \; T\bar{\lambda} \cdot w)
+ \tilde{\omega}_{\mathcal{O}_{(\mu,a)} \times V\times
V*}^{-}(T\bar{\lambda}\cdot \tilde{X}\cdot \varepsilon, \; T\bar{\lambda}\cdot w),
\end{align*}
where we have used that $T\pi_{(\mu,a)}(X_H)= X_{h_{(\mu,a)}}. $
Note that $\varepsilon: T^* Q \rightarrow T^* Q $ is
symplectic, and $\pi_{(\mu,a)}^*\tilde{\omega}_{\mathcal{O}_{(\mu,a)} \times V\times
V*}^{-}= i_{(\mu,a)}^*\omega_Q
= \omega_Q, $ along $\textmd{Im}(\gamma)$, and hence
$\bar{\varepsilon}= \pi_{(\mu,a)}(\varepsilon): \mathbf{J}^{-1}(\mu,a)(\subset T^* Q) \rightarrow (T^*
Q)_{(\mu,a)} $ is also symplectic along
$\textmd{Im}(\gamma)$, and $X_{h_{(\mu,a)}}\cdot \bar{\varepsilon}=
T\bar{\varepsilon} \cdot X_{h_{(\mu,a)} \cdot \bar{\varepsilon}}, $
along $\textmd{Im}(\gamma)\cap\textmd{Im}(\varepsilon)$.
From the above arguments, we can obtain that
\begin{align*}
& \tilde{\omega}_{\mathcal{O}_{(\mu,a)} \times V\times
V*}^{-}(T\bar{\gamma} \cdot \tilde{X}^\varepsilon, \; T\pi_{(\mu,a)} \cdot w)-
\tilde{\omega}_{\mathcal{O}_{(\mu,a)} \times V\times
V*}^{-}(X_{h_{(\mu,a)}} \cdot \bar{\varepsilon}, \; T\pi_{(\mu,a)} \cdot w) \\
& =- \tilde{\omega}_{\mathcal{O}_{(\mu,a)} \times V\times
V*}^{-}(X_{h_{(\mu,a)}} \cdot
\bar{\varepsilon}, \; T\bar{\lambda} \cdot w)
+ \tilde{\omega}_{\mathcal{O}_{(\mu,a)} \times V\times
V*}^{-}(T\bar{\lambda}\cdot \tilde{X}\cdot \varepsilon, \; T\bar{\lambda}\cdot w)\\
& = \tilde{\omega}_{\mathcal{O}_{(\mu,a)} \times V\times
V*}^{-}(T\bar{\lambda}\cdot \tilde{X}\cdot \varepsilon, \; T\bar{\lambda}\cdot w)
-\tilde{\omega}_{\mathcal{O}_{(\mu,a)} \times V\times
V*}^{-}(T\bar{\varepsilon} \cdot X_{h_{(\mu,a)} \cdot \bar{\varepsilon}}, \; T\bar{\lambda} \cdot w)\\
& = \tilde{\omega}_{\mathcal{O}_{(\mu,a)} \times V\times
V*}^{-}(T\bar{\lambda}\cdot \tilde{X}\cdot \varepsilon- T\bar{\varepsilon} \cdot X_{h_{(\mu,a)} \cdot \bar{\varepsilon}}, \; T\bar{\lambda}\cdot w).
\end{align*}
Because the $R_p$-reduced symplectic form
$\tilde{\omega}_{\mathcal{O}_{(\mu,a)} \times V\times
V*}^{-}$ is non-degenerate, it follows that
$T\bar{\gamma}\cdot \tilde{X}^\varepsilon=
X_{h_{(\mu,a)}}\cdot \bar{\varepsilon}, $ is equivalent to
$T\bar{\lambda}\cdot \tilde{X}\cdot \varepsilon
= T\bar{\varepsilon} \cdot X_{h_{(\mu,a)} \cdot \bar{\varepsilon}}. $
Thus, we know that the $\varepsilon$ and $\bar{\varepsilon}$ satisfy
the equation $T\bar{\varepsilon}\cdot(X_{h_{(\mu,a)} \cdot \bar{\varepsilon}})
= T\bar{\lambda}\cdot \tilde{X} \cdot \varepsilon, $
if and only if they satisfy the Type II of Hamilton-Jacobi equation
$T\bar{\gamma}\cdot \tilde{X}^\varepsilon= X_{h_{(\mu,a)}}\cdot \bar{\varepsilon}. $
\hskip 0.3cm $\blacksquare$\\

In particular, when $Q=W, $ we can obtain the two types of Hamilton-Jacobi
theorems for the regular point reducible RCH system and the $R_p$-reduced RCH system on
semidirect product Lie group $W$. In this case, note that the
symplecitic structure on the coadjoint orbit $\mathcal{O}_{(\mu,a)}$
is induced by the (-)-semidirect product Lie-Poisson brackets on
Lie algebra $\mathfrak{w}^\ast$, then the two types of Hamilton-Jacobi equations for
the regular point reducible RCH system $(T^\ast W,W,\omega_W,H,F,\mathcal{C})$
and the $R_p$-reduced RCH system
$(\mathcal{O}_{(\mu,a)}, \omega_{\mathcal{O}_{(\mu,a)}}^{-},
h_{(\mu,a)}, f_{(\mu,a)}, u_{(\mu,a)})$ are also called the two types of Lie-Poisson
Hamilton-Jacobi equations. See Wang \cite{wa13d}, Marsden and Ratiu
\cite{mara99}, and Ge and Marsden \cite{gema88}.\\

It is worthy of noting that, for a given regular point reducible RCH system
on the generalization of a semidirect product Lie group
$(T^\ast Q,W,\omega_Q,H,F,\mathcal{C})$ with the $R_p$-reduced RCH system
\\ $(\mathcal{O}_{(\mu,a)} \times V\times
V^\ast, \tilde{\omega}_{\mathcal{O}_{(\mu,a)} \times V \times V^\ast
}^{-},h_{(\mu,a)},f_{(\mu,a)},\mathcal{C}_{(\mu,a)})$, we know that the Hamiltonian vector fields
$X_{H}$ and $X_{h_{(\mu,a)}}$ for the corresponding
Hamiltonian system $(T^*Q,W,\omega_Q, H)$
and its $R_p$-reduced system $(\mathcal{O}_{(\mu,a)} \times V\times
V^{*},\tilde{\omega}_{\mathcal{O}_{(\mu,a)} \times V
\times V^{*}}^{-}, h_{(\mu,a)} )$, are $\pi_{(\mu,a)}$-related, that is,
$X_{h_{(\mu,a)}}\cdot \pi_{(\mu,a)}=T\pi_{(\mu,a)}\cdot X_{H}\cdot i_{(\mu,a)}.$
By using the $R_p$-reduced symplectic form
$\tilde{\omega}_{\mathcal{O}_{(\mu,a)} \times V \times V^\ast}^{-}$,
then we can
prove the following Theorem 3.8, which states the relationship
between the solutions of Type II of Hamilton-Jacobi equations and the
regular point reduction.

\begin{theo}
For a given regular point reducible RCH system
on the generalization of a semidirect product Lie group
$(T^\ast Q,W,\omega_Q,H,F,\mathcal{C})$ with the $R_p$-reduced RCH system
$(\mathcal{O}_{(\mu,a)} \times V\times
V^\ast, \\ \tilde{\omega}_{\mathcal{O}_{(\mu,a)} \times V \times V^\ast
}^{-},h_{(\mu,a)},f_{(\mu,a)}, u_{(\mu,a)}) $,
assume that $\gamma: Q \rightarrow T^*Q$ is an one-form on $Q$,
and $\varepsilon: T^* Q \rightarrow T^* Q $ is a
$W_{(\mu,a)}$-invariant symplectic map,
$\bar{\varepsilon}=\pi_{(\mu,a)}(\varepsilon): \mathbf{J}_Q^{-1}(\mu, a )\rightarrow
\mathcal{O}_{(\mu,a)}\times V\times V^{*}. $
Under the hypotheses and notations of Theorem 3.7, then we have that
$\varepsilon$ is a solution of the Type II of Hamilton-Jacobi equation
$T\gamma\cdot \tilde{X}^\varepsilon= X_H\cdot \varepsilon, $ for the
regular point reducible RCH system $(T^\ast Q,W,\omega_Q,H,F,u), $ if and only if
$\varepsilon$ and $\bar{\varepsilon} $ satisfy the Type II of Hamilton-Jacobi equation
$T\bar{\gamma}\cdot \tilde{X}^\varepsilon= X_{h_{(\mu,a)}}\cdot \bar{\varepsilon}, $ for the
$R_p$-reduced RCH system
$(\mathcal{O}_{(\mu,a)} \times V \times
V^{*},\tilde{\omega}_{\mathcal{O}_{(\mu,a)} \times V
\times V^{*}}^{-}, h_{(\mu,a)}, f_{(\mu,a)}, u_{(\mu,a)}) $.
\end{theo}

\noindent{\bf Proof: }
Note that $\textmd{Im}(\gamma)\subset \mathbf{J}_Q^{-1}(\mu,a), $ and it
is $W_{(\mu,a)}$-invariant, in this case,
$\pi_{(\mu,a)}^*\tilde{\omega}_{\mathcal{O}_{(\mu,a)} \times V \times V^\ast}^{-}
= i_{(\mu,a)}^*\omega_Q= \omega_Q, $ along $\textmd{Im}(\gamma)$.
Since the Hamiltonian vector fields
$X_{H}$ and $X_{h_{(\mu,a)}}$ are $\pi_{(\mu,a)}$-related, that is,
$X_{h_{(\mu,a)}}\cdot \pi_{(\mu,a)}= T\pi_{(\mu,a)}\cdot X_{H}\cdot i_{(\mu,a)}, $ and
by using the $R_p$-reduced symplectic form
$\tilde{\omega}_{\mathcal{O}_{(\mu,a)} \times V \times V^\ast}^{-}$,
for any $w \in TT^* Q,$ and $T\pi_{(\mu,a)} \cdot w
\neq 0, $ we have that
\begin{align*}
& \tilde{\omega}_{\mathcal{O}_{(\mu,a)} \times V \times V^\ast}^{-}
(T\bar{\gamma} \cdot \tilde{X}^\varepsilon
- X_{h_{(\mu,a)}} \cdot \bar{\varepsilon}, \; T\pi_{(\mu,a)} \cdot w) \\
& = \tilde{\omega}_{\mathcal{O}_{(\mu,a)} \times V \times V^\ast}^{-}
(T\bar{\gamma} \cdot \tilde{X}^\varepsilon, \; T\pi_{(\mu,a)} \cdot w)-
\tilde{\omega}_{\mathcal{O}_{(\mu,a)} \times V \times V^\ast}^{-}
(X_{h_{(\mu,a)}} \cdot \bar{\varepsilon}, \; T\pi_{(\mu,a)} \cdot w) \\
& = \tilde{\omega}_{\mathcal{O}_{(\mu,a)} \times V \times V^\ast}^{-}
(T\pi_{(\mu,a)} \cdot T\gamma \cdot \tilde{X}^\varepsilon, \; T\pi_{(\mu,a)} \cdot w)-
\tilde{\omega}_{\mathcal{O}_{(\mu,a)} \times V \times V^\ast}^{-}
(X_{h_{(\mu,a)}} \cdot \pi_{(\mu,a)} \cdot \varepsilon, \; T\pi_{(\mu,a)} \cdot w) \\
& = \pi_{(\mu,a)}^*\tilde{\omega}_{\mathcal{O}_{(\mu,a)} \times V \times V^\ast}^{-}
(T\gamma \cdot \tilde{X}^\varepsilon, \; w)
-\tilde{\omega}_{\mathcal{O}_{(\mu,a)} \times V \times V^\ast}^{-}
(T\pi_{(\mu,a)}\cdot X_{H}\cdot \varepsilon, \; T\pi_{(\mu,a)} \cdot w) \\
& = \pi_{(\mu,a)}^*\tilde{\omega}_{\mathcal{O}_{(\mu,a)} \times V \times V^\ast}^{-}
(T\gamma \cdot \tilde{X}^\varepsilon, \; w)
-\pi_{(\mu,a)}^*\tilde{\omega}_{\mathcal{O}_{(\mu,a)} \times V \times V^\ast}^{-}
(X_{H}\cdot \varepsilon, \; w)\\
& = \omega_Q(T\gamma \cdot \tilde{X}^\varepsilon- X_{H}\cdot \varepsilon, \; w).
\end{align*}
Because both the symplectic form $\omega_Q$ and the $R_p$-reduced symplectic form
$\tilde{\omega}_{\mathcal{O}_{(\mu,a)} \times V \times V^\ast}^{-}$
are non-degenerate, it follows that the equation
$T\bar{\gamma}\cdot \tilde{X}^\varepsilon= X_{h_{(\mu,a)}}\cdot \bar{\varepsilon}, $
is equivalent to the equation $T\gamma\cdot \tilde{X}^\varepsilon= X_H\cdot \varepsilon$. Thus,
$\varepsilon$ is a solution of the Type II of Hamilton-Jacobi equation
$T\gamma\cdot \tilde{X}^\varepsilon= X_H\cdot \varepsilon, $ for the
regular point reducible RCH system
on the generalization of a semidirect product Lie group
$(T^\ast Q,W,\omega_Q,H,F,\mathcal{C}), $ if and only if
$\varepsilon$ and $\bar{\varepsilon} $ satisfy the Type II of Hamilton-Jacobi equation
$T\bar{\gamma}\cdot \tilde{X}^\varepsilon= X_{h_{(\mu,a)}}\cdot \bar{\varepsilon}, $ for the
$R_p$-reduced RCH system
$(\mathcal{O}_{(\mu,a)} \times V \times V^\ast,
\tilde{\omega}_{\mathcal{O}_{(\mu,a)} \times V \times V^\ast}^{-},
h_{(\mu,a)}, u_{(\mu,a)}) $.  \hskip 0.3cm
$\blacksquare$ \\

As an application of the theoretical result, in the following we
consider the controlled underwater vehicle-rotor system as a
regular point reducible RCH system on the generalizations
of the semidirect product Lie groups $(\textmd{SO}(3)\circledS
\mathbb{R}^3)\times \mathbb{R}^2$ and $(\textmd{SE}(3)\circledS
\mathbb{R}^3)\times \mathbb{R}^2$, respectively, and derive precisely
the geometric constraint conditions
of the $R_p$-reduced symplectic forms for the
dynamical vector fields of the regular point reducible controlled
underwater vehicle-rotor systems by calculation in detail,
that is, the two types of Hamilton-Jacobi equations of their $R_p$-reduced systems,
which reveal the deeply internal relationships of the geometrical structures
of phase spaces, the dynamical vector fields and controls
of the controlled underwater vehicle-rotor system.

\section{Underwater Vehicle-Rotor System with Coincident Centers
of Buoyancy and Gravity}

In this section, we consider the controlled underwater vehicle-rotor system
with coincident centers of buoyancy and gravity as an RCH
system on the generalization of semidirect product Lie group
$Q=W\times V$, where $W=\textmd{SO}(3)\circledS \mathbb{R}^3$
and $V= \mathbb{R}^2$. In this case, we give precisely the regular point reduction and
the two types of Hamilton-Jacobi equation for
the $R_p$-reduced controlled underwater vehicle-rotor system.

\subsection{Symmetric Reduction }

In the following we assume that semidirect product Lie group $W=
\textmd{SO}(3)\circledS \mathbb{R}^3=\textmd{SE}(3)$ acts freely and
properly on $Q= \textmd{SE}(3)\times \mathbb{R}^2 $ by the left translation on
$\textmd{SE}(3)$, then the action of $\textmd{SE}(3)$ on the phase
space $T^\ast Q$ is by cotangent lift of the left translation on
$\textmd{SE}(3)$ at the identity, that is, $\Phi^{T*}:
\textmd{SE}(3)\times T^\ast \textmd{SE}(3)\times T^\ast \mathbb{R}^2 \cong
\textmd{SE}(3)\times \textmd{SE}(3)\times \mathfrak{se}^\ast(3)
\times \mathbb{R}^2 \times \mathbb{R}^{2*}\to \textmd{SE}(3)\times
\mathfrak{se}^\ast(3)\times \mathbb{R}^2 \times \mathbb{R}^{2*},$ given
by $$\Phi^{T*}((B,s)((A,c),(\Pi,w),\theta,l))=((BA,s+Bc),(\Pi,w),\theta,l), $$
for any $A,B\in \textmd{SO}(3), \; \Pi \in \mathfrak{so}^\ast(3), \;
s,c \in \mathbb{R}^3, \; (A,c),(B,s) \in \textmd{SE}(3), \; w \in \mathbb{R}^{3*},
\; (\Pi, w) \in \mathfrak{se}^\ast(3), \; \theta \in \mathbb{R}^2, \; l \in \mathbb{R}^{2*} $,
and $\textmd{SO}(3)$ acts on $\mathbb{R}^3$ in the standard way. Assume that
the action is free, proper and symplectic, and admits an associated
$\operatorname{Ad}^\ast$-equivariant momentum map $\mathbf{J}_Q:
T^\ast Q \cong \textmd{SE}(3)\times \mathfrak{se}^\ast(3) \times
\mathbb{R}^2 \times \mathbb{R}^{2*} \to \mathfrak{se}^\ast(3)$ for the
cotangent lift $\textmd{SE}(3)$-action. If $(\Pi,w) \in \mathfrak{se}^\ast(3)$
is a regular value of $\mathbf{J}_Q$, then the $R_p$-reduced
space $(T^\ast Q)_{(\Pi,w)}=
\mathbf{J}^{-1}_Q(\Pi,w)/\textmd{SE}(3)_{(\Pi,w)}$ is symplectically
diffeomorphic to the orbit space $\mathcal{O}_{(\Pi,w)} \times
\mathbb{R}^2 \times \mathbb{R}^{2*} \subset \mathfrak{se}^\ast(3)
\times \mathbb{R}^2 \times \mathbb{R}^{2*} $, where
$\textmd{SE}(3)_{(\Pi,w)}$ is the isotropy subgroup of co-adjoint
$\textmd{SE}(3)$-action at the point
$(\Pi,w)\in\mathfrak{se}^\ast(3)$.\\

From Marsden et al. \cite{mamiorpera07} we know that
$\mathfrak{se}^\ast(3)= \mathfrak{so}^\ast(3)\circledS
\mathbb{R}^{3*}$ is a Poisson manifold with respect to its semidirect
product Lie-Poisson bracket, that is, heavy top Lie-Poisson bracket
defined by
\begin{equation}   \{F,K\}_{\mathfrak{se}^\ast(3)}(\Pi,P)= -\Pi\cdot (\nabla_\Pi
F\times \nabla_\Pi K)-P\cdot(\nabla_\Pi F\times \nabla_P
K-\nabla_\Pi K\times \nabla_P F), \;\; \label{4.1}
\end{equation}
where $F,K: \mathfrak{se}^\ast(3)\to \mathbb{R}, $ and $ (\Pi,P) \in
\mathfrak{se}^\ast(3)$. For a point $(\Pi_0,P_0)=(\mu,a) \in \mathfrak{se}^\ast(3)$, the
co-adjoint orbit $\mathcal{O}_{(\mu,a)} \subset
\mathfrak{se}^\ast(3)$ has the induced orbit symplectic form
$\omega^{-}_{\mathcal{O}_{(\mu,a)}}$, which coincides with the
restriction of the Lie-Poisson bracket on $\mathfrak{se}^\ast(3)$ to
the co-adjoint orbit $\mathcal{O}_{(\mu,a)}$, and the co-adjoint
orbits $(\mathcal{O}_{(\mu,a)}, \omega_{\mathcal{O}_{(\mu,a)}}^{-}) $,
$ (\mu,a)\in \mathfrak{se}^\ast(3), $ form the symplectic leaves of
the Poisson manifold
$(\mathfrak{se}^\ast(3),\{\cdot,\cdot\}_{\mathfrak{se}^\ast(3)}). $
Let $\omega_{\mathbb{R}^2}$ be the canonical symplectic form on
$T^\ast \mathbb{R}^2 \cong \mathbb{R}^2 \times \mathbb{R}^{2*}$, and it
induces a canonical Poisson bracket $\{\cdot,\cdot\}_{\mathbb{R}^2}$
on $T^\ast \mathbb{R}^2$ given by $(3.7)$. Thus, we can induce a
symplectic form $\tilde{\omega}^{-}_{\mathcal{O}_{(\mu,a)} \times
\mathbb{R}^2 \times \mathbb{R}^{2*}}= \pi_{\mathcal{O}_{(\mu,a)}}^\ast
\omega^{-}_{\mathcal{O}_{(\mu,a)}}+ \pi_{\mathbb{R}^2}^\ast
\omega_{\mathbb{R}^2}$ on the smooth manifold $\mathcal{O}_{(\mu,a)}
\times \mathbb{R}^2 \times \mathbb{R}^{2*} $, where the maps
$\pi_{\mathcal{O}_{(\mu,a)}}: \mathcal{O}_{(\mu,a)} \times
\mathbb{R}^2 \times \mathbb{R}^{2*} \to \mathcal{O}_{(\mu,a)}$ and
$\pi_{\mathbb{R}^2}: \mathcal{O}_{(\mu,a)} \times \mathbb{R}^2
\times \mathbb{R}^{2*} \to \mathbb{R}^2 \times \mathbb{R}^{2*}$ are
canonical projections, and induce a Poisson bracket
$\{\cdot,\cdot\}_{-}= \pi_{\mathfrak{se}^\ast(3)}^\ast
\{\cdot,\cdot\}_{\mathfrak{se}^\ast(3)}+ \pi_{\mathbb{R}^2}^\ast
\{\cdot,\cdot\}_{\mathbb{R}^2}$ on the smooth manifold
$\mathfrak{se}^\ast(3)\times \mathbb{R}^2 \times \mathbb{R}^{2*} $,
where the maps $\pi_{\mathfrak{se}^\ast(3)}: \mathfrak{se}^\ast(3)
\times \mathbb{R}^2 \times \mathbb{R}^{2*} \to \mathfrak{se}^\ast(3)$
and $\pi_{\mathbb{R}^2}: \mathfrak{se}^\ast(3) \times \mathbb{R}^2
\times \mathbb{R}^{2*} \to \mathbb{R}^2 \times \mathbb{R}^{2*}$ are
canonical projections, and such that $(\mathcal{O}_{(\mu,a)} \times
\mathbb{R}^2 \times
\mathbb{R}^{2*},\tilde{\omega}_{\mathcal{O}_{(\mu,a)} \times
\mathbb{R}^2 \times \mathbb{R}^{2*}}^{-})$ is a symplectic leaf of the
Poisson manifold
$(\mathfrak{se}^\ast(3) \times \mathbb{R}^2 \times \mathbb{R}^{2*}, \{\cdot,\cdot\}_{-}). $\\

From the above expression $(2.1)$ of the Hamiltonian, we know that
$H(A,c,\Pi,P,\theta,l)$ is invariant under the cotangent lift of the left
$\textmd{SE}(3)$-action $\Phi^{T*}:\textmd{SE}(3)\times T^\ast Q\to
T^\ast Q$. Moreover, from the heavy top Lie-Poisson
bracket on $\mathfrak{se}^\ast(3)$ and the Poisson bracket on
$T^\ast \mathbb{R}^2$, we can get the Poisson bracket
on $\mathfrak{se}^\ast(3)\times \mathbb{R}^2\times \mathbb{R}^{2*} $,
that is, for $F,K: \mathfrak{se}^\ast(3)\times \mathbb{R}^2\times
\mathbb{R}^{2*} \to \mathbb{R}, $ we have that
\begin{align} & \{F,K\}_{-}(\Pi, P,\theta,l) \nonumber \\
& = -\Pi\cdot(\nabla_\Pi F\times\nabla_\Pi K)-P\cdot(\nabla_\Pi
F\times \nabla_P K-\nabla_\Pi K\times \nabla_P F)+
\{F,K\}_{\mathbb{R}^2}(\theta,l). \label{4.2}
\end{align}
such that the Hamiltonian vector fields of underwater vehicle-rotor system are given by
\begin{align*}
& X_{H}(\Pi) = \{\Pi,\; H \}_{-}= -\Pi\cdot(\nabla_\Pi\Pi\times\nabla_\Pi
H) -P\cdot(\nabla_\Pi\Pi\times\nabla_P
H-\nabla_\Pi H \times\nabla_P \Pi) + \{\Pi,\; H \}_{\mathbb{R}^2}\\
& =(\Pi_1,\Pi_2,\Pi_3)\times (\frac{\Pi_1-l_1}{ \bar{I}_1}, \frac{\Pi_2-l_2}{
\bar{I}_2}, \frac{\Pi_3}{\bar{I}_3})
+(P_1,P_2,P_3)\times (\frac{P_1}{m_1},\frac{P_2}{m_2},\frac{P_3}{m_3})
+ \sum_{k=1}^2(\frac{\partial \Pi}{\partial \theta_k}
\frac{\partial H}{\partial l_k}- \frac{\partial
H}{\partial \theta_k}\frac{\partial \Pi}{\partial l_k})\\
& = ( \frac{(\bar{I}_2-\bar{I}_3)\Pi_2\Pi_3+
\bar{I}_3\Pi_3l_2}{\bar{I}_2\bar{I}_3}+
 \frac{(m_2-m_3)P_2P_3}{m_2m_3}, \;\;
\frac{(\bar{I}_3-\bar{I}_1)\Pi_3\Pi_1-
\bar{I}_3\Pi_3l_1}{\bar{I}_3\bar{I}_1}+ \frac{(m_3-m_1)P_3P_1}{m_3m_1}, \\
& \;\;\;\;\;\;
\frac{(\bar{I}_1-\bar{I}_2)\Pi_1\Pi_2-
\bar{I}_1\Pi_1l_2 + \bar{I}_2\Pi_2l_1}{\bar{I}_1\bar{I}_2}
 +  \frac{(m_1-m_2)P_1P_2}{m_1m_2}),
\end{align*}
since $\nabla_{\Pi_i}\Pi_i=1,\; \nabla_{\Pi_i}\Pi_j=0, \; i\neq j, \;  \nabla_{\Pi_i}P_j=\nabla_{P_i}\Pi_j=0$
and $ \nabla_{P_i} H= \frac{P_i}{m_i}, \; \nabla_{\Pi_k} H= (\Pi_k-l_k)/\bar{I}_k ,\; \nabla_{\Pi_3}
H= \Pi_3/\bar{I}_3 , \; \frac{\partial
\Pi}{\partial \theta_k}= \frac{\partial H}{\partial
\theta_k}=0, \; i, j=1,2,3, \; k=1,2. $

\begin{align*}
 X_{H}(P)&
= \{P,\; H\}_{-} =-\Pi\cdot(\nabla_\Pi P\times\nabla_\Pi
H)-P\cdot(\nabla_\Pi P\times\nabla_P
H-\nabla_\Pi H \times\nabla_P P)+  \{P,\; H\}_{\mathbb{R}^2} \\
& =\nabla_P P\cdot (P \times\nabla_\Pi
H)+ \sum_{k=1}^2(\frac{\partial P}{\partial \theta_k}
\frac{\partial H}{\partial l_k}- \frac{\partial
H}{\partial \theta_k}\frac{\partial P}{\partial l_k})\\
& = (P_1,P_2,P_3)\times (\frac{(\Pi_1- l_1)}{ \bar{I}_1},\;\;
\frac{(\Pi_2- l_2)}{ \bar{I}_2}, \;\; \frac{\Pi_3}{\bar{I}_3}) \\
&= ( \frac{\bar{I}_2P_2\Pi_3-\bar{I}_3P_3\Pi_2+
\bar{I}_3\Pi_3l_2}{\bar{I}_2\bar{I}_3}, \;\;
\frac{\bar{I}_3P_3\Pi_1-\bar{I}_1P_1\Pi_3-
\bar{I}_3\Pi_3l_1}{\bar{I}_3\bar{I}_1}, \\
& \;\;\;\;\;\;
\frac{\bar{I}_1P_1\Pi_2-\bar{I}_2P_2\Pi_1-
\bar{I}_1\Pi_1l_2 + \bar{I}_2\Pi_2l_1}{\bar{I}_1\bar{I}_2} ),
\end{align*}
since $\nabla_{P_i}P_i =1, \; \nabla_{P_i}P_j=0, \; i\neq j, \; \nabla_{\Pi_i}P_j =0, $ and
$\nabla_{\Pi_k} H= (\Pi_k-l_k)/\bar{I}_k ,\; \nabla_{\Pi_3}
H= \Pi_3/\bar{I}_3 , \; \frac{\partial
P_j}{\partial \theta_k}= \frac{\partial H}{\partial
\theta_k}=0, \; i, j= 1,2,3, \; k=1,2.$

\begin{align*}
 X_{H}(\theta)& = \{\theta,\; H \}_{-}= -\Pi\cdot(\nabla_\Pi\theta\times\nabla_\Pi
H) -P\cdot(\nabla_\Pi\theta\times\nabla_P
H-\nabla_\Pi H \times\nabla_P \theta) + \{\theta,\; H \}_{\mathbb{R}^2}\\
& =\sum_{k=1}^2(\frac{\partial \theta}{\partial \theta_k}
\frac{\partial H}{\partial l_k}- \frac{\partial
H}{\partial \theta_k}\frac{\partial \theta}{\partial l_k})= ( -\frac{(\Pi_1- l_1)}{ \bar{I}_1}
+\frac{l_1}{J_1},\;\; -\frac{(\Pi_2-
l_2)}{ \bar{I}_2} +\frac{l_2}{J_2} ),
\end{align*}
since $\nabla_{\Pi_i}\theta=\nabla_{P_i}\theta=0, $ $\frac{\partial \theta_k}{\partial
\theta_k}= 1, \; \frac{\partial \theta_n}{\partial \theta_k}=
0, \; n\neq k, \;\; \frac{\partial H}{\partial \theta_k}=0, $
and $\frac{\partial H}{\partial l_k}= -(\Pi_k-l_k)/\bar{I}_k
+\frac{l_k}{J_k}, \; i= 1,2,3, \; k, n=1,2. $

\begin{align*}
 X_{H}(l)& = \{l,\; H \}_{-}= -\Pi\cdot(\nabla_\Pi l \times\nabla_\Pi
H) -P\cdot(\nabla_\Pi l \times\nabla_P
H-\nabla_\Pi H \times\nabla_P l) + \{l,\; H \}_{\mathbb{R}^2}\\
& = \sum_{k=1}^2(\frac{\partial l}{\partial \theta_k}
\frac{\partial H}{\partial l_k}- \frac{\partial
H}{\partial \theta_k}\frac{\partial l}{\partial l_k})=(0,0),
\end{align*}
since $\nabla_{\Pi_i} l=\nabla_{P_i} l=0,$ and $\frac{\partial l}{\partial \theta_k}=
\frac{\partial H}{\partial \theta_k}=0, \; i=1,2,3, \; k=1,2. $\\

Moreover, if we consider the underwater vehicle-rotor system with a control torque $u: T^\ast
Q \to \mathcal{C}$ acting on the rotors, where the control subset $\mathcal{C}\subset T^* Q $ is a fiber submanifold,
and assume that $u\in \mathcal{C}$ is invariant under the cotangent lift of the left
$\textmd{SE}(3)$-action, and
the dynamical vector field of the regular point reducible
controlled underwater vehicle-rotor system $(T^\ast Q,\textmd{SE}(3),\omega_Q,H,u)$
can be expressed by
\begin{align}
\tilde{X}= X_{(T^\ast Q,\textmd{SE}(3),\omega_Q,H,u)}= X_H+ \textnormal{vlift}(u),
\label{4.3} \end{align}
where $\textnormal{vlift}(u)= \textnormal{vlift}(u)\cdot X_H $
is the change of $X_H$ under the action of the control torque $u$.\\

Since the Hamiltonian
$H(A,c,\Pi,P,\theta,l)$ is invariant under the cotangent lift $\Phi^{T^*}$ of the left
$\textmd{SE}(3)$-action, for the point $(\Pi_0,P_0)=(\mu,a)\in
\mathfrak{se}^\ast(3)$, the regular value of $\mathbf{J}_Q$, we
have the $R_p$-reduced Hamiltonian
$h_{(\mu,a)}(\Pi,P,\theta,l):\mathcal{O}_{(\mu,a)}\times\mathbb{R}^2
\times \mathbb{R}^{2*} (\subset \mathfrak{se}^\ast (3)\times
\mathbb{R}^2\times \mathbb{R}^{2*}) \to \mathbb{R}$ given by
$$h_{(\mu,a)}(\Pi,P,\theta,l)\cdot \pi_{(\mu,a)}
=H(A,c,\Pi,P,\theta,l)|_{\mathcal{O}_{(\mu,a)}\times
\mathbb{R}^2\times \mathbb{R}^{2*}}, $$ where $\pi_{(\mu,a)}:
\mathbf{J}_Q^{-1}(\mu,a) \rightarrow \mathcal{O}_{(\mu,a)}
\times\mathbb{R}^2\times\mathbb{R}^{2*}$.
Moreover, for the $R_p$-reduced
Hamiltonian $h_{(\mu,a)}(\Pi,P, \theta, l): \mathcal{O}_{(\mu,a)}\times
\mathbb{R}^2\times \mathbb{R}^{2*} \to \mathbb{R}$, we have the
Hamiltonian vector field
$$X_{h_{(\mu,a)}}(K_{(\mu,a)})=\{K_{(\mu,a)},h_{(\mu,a)}\}_{-}|_{\mathcal{O}_{(\mu,a)}
\times \mathbb{R}^2\times \mathbb{R}^{2*}}, $$ where
$K_{(\mu,a)}(\Pi,P, \theta, l): \mathcal{O}_{(\mu,a)}\times
\mathbb{R}^2\times \mathbb{R}^{2*} \to \mathbb{R}.$
Assume that $u\in \mathcal{C} \cap
\mathbf{J}^{-1}_Q(\mu,a)$ and the  $R_p$-reduced control torque $u_{(\mu,a)}:
\mathcal{O}_{(\mu,a)} \times\mathbb{R}^2\times\mathbb{R}^{2*} \to \mathcal{C}_{(\mu,a)}
(\subset \mathcal{O}_{(\mu,a)} \times\mathbb{R}^2\times\mathbb{R}^{2*}) $ is
given by $u_{(\mu,a)}(\Pi,P,\theta,l)\cdot \pi_{(\mu,a)}=
u(A,c,\Pi,P, \theta,l )|_{\mathcal{O}_{(\mu,a)}
\times\mathbb{R}^2\times\mathbb{R}^{2*} }, $ where $ \mathcal{C}_{(\mu,a)}
= \pi_{(\mu,a)}(\mathcal{C}\cap \mathbf{J}_Q^{-1}(\mu,a)). $
The $R_p$-reduced controlled underwater vehicle-rotor
system is the 4-tuple $(\mathcal{O}_{(\mu,a)} \times \mathbb{R}^2 \times
\mathbb{R}^{2*},\tilde{\omega}_{\mathcal{O}_{(\mu,a)} \times \mathbb{R}^2
\times \mathbb{R}^{2*}}^{-},h_{(\mu,a)},u_{(\mu,a)}), $ where
$\tilde{\omega}_{\mathcal{O}_{(\mu,a)} \times \mathbb{R}^2 \times
\mathbb{R}^{2*}}^{-}$ is the induced symplectic form
from the above Poisson bracker on  $\mathcal{O}_{(\mu,a)}
\times \mathbb{R}^2\times \mathbb{R}^{2*} ,$ such that
the Hamiltonian vector field
$$X_{h_{(\mu,a)}}(K_{(\mu,a)})=\tilde{\omega}_{\mathcal{O}_{(\mu,a)} \times \mathbb{R}^2 \times
\mathbb{R}^{2*}}^{-}(X_{K_{(\mu,a)}}, X_{h_{(\mu,a)}})
=\{K_{(\mu,a)},h_{(\mu,a)}\}_{-}|_{\mathcal{O}_{(\mu,a)}
\times\mathbb{R}^2\times\mathbb{R}^{2*} }.$$
Moreover, assume that the dynamical vector field of the $R_p$-reduced controlled
underwater vehicle-rotor system $(\mathcal{O}_{(\mu,a)} \times \mathbb{R}^2 \times
\mathbb{R}^{2*},\tilde{\omega}_{\mathcal{O}_{(\mu,a)} \times
\mathbb{R}^2 \times
\mathbb{R}^{2*}}^{-},h_{(\mu,a)},u_{(\mu,a)})$ can be expressed
by
\begin{align} X_{(\mathcal{O}_{(\mu,a)} \times \mathbb{R}^2 \times
\mathbb{R}^{2*},\tilde{\omega}_{\mathcal{O}_{(\mu,a)} \times \mathbb{R}^2 \times
\mathbb{R}^{2*}}^{-},h_{(\mu,a)},u_{(\mu,a)})} = X_{h_{(\mu,a)}} +
\mbox{vlift}(u_{(\mu,a)}), \label{4.4}
\end{align}
where $\mbox{vlift}(u_{(\mu,a)})=
\mbox{vlift}(u_{(\mu,a)})X_{h_{(\mu,a)}} \in T(\mathcal{O}_{(\mu,a)}
\times \mathbb{R}^2 \times
\mathbb{R}^{2*}), $ is the change of $X_{h_{(\mu,a)}}$
under the action of the $R_p$-reduced control torque $u_{(\mu,a)}$,
The dynamical vector fields of the controlled
underwater vehicle-rotor system and the $R_p$-reduced controlled
underwater vehicle-rotor system satisfy the condition
\begin{equation}X_{(\mathcal{O}_{(\mu,a)} \times \mathbb{R}^2 \times
\mathbb{R}^{2*},\tilde{\omega}_{\mathcal{O}_{(\mu,a)} \times \mathbb{R}^2 \times
\mathbb{R}^{2*}}^{-}, h_{(\mu,a)}, u_{(\mu,a)})}\cdot \pi_{(\mu,a)}
=T\pi_{(\mu,a)}\cdot X_{(T^\ast Q,\textmd{SE}(3),\omega_Q,H,u)}\cdot i_{(\mu,a)}.
\label{4.5}\end{equation}
See Marsden et al \cite{mawazh10} and Wang \cite{wa18}.\\

To sum up the above discussion, we have the following Theorem 4.1.
\begin{theo}
In the case of coincident centers of buoyancy and gravity,
the underwater vehicle-rotor system with the control torque $u$ acting on the
rotors, that is, the 5-tuple $(T^\ast Q,\textmd{SE}(3),\omega_Q,H,u), $ where $Q=
\textmd{SE}(3)\times \mathbb{R}^2, $ is a regular point reducible
RCH system. For a point $(\mu,a) \in \mathfrak{se}^\ast(3)$, the regular
value of the momentum map $\mathbf{J}_Q: T^* Q \cong \textmd{SE}(3)\times
\mathfrak{se}^\ast(3) \times \mathbb{R}^2 \times \mathbb{R}^{2*} \to
\mathfrak{se}^\ast(3)$, the $R_p$-reduced controlled underwater vehicle-rotor
system is the 4-tuple $(\mathcal{O}_{(\mu,a)} \times \mathbb{R}^2 \times
\mathbb{R}^{2*},\tilde{\omega}_{\mathcal{O}_{(\mu,a)} \times \mathbb{R}^2
\times \mathbb{R}^{2*}}^{-}, \\ h_{(\mu,a)},u_{(\mu,a)}), $ where $\mathcal{O}_{(\mu,a)}
\subset \mathfrak{se}^\ast(3)$ is the co-adjoint orbit,
$\tilde{\omega}_{\mathcal{O}_{(\mu,a)} \times \mathbb{R}^2 \times
\mathbb{R}^{2*}}^{-}$ is the induced symplectic form on $\mathcal{O}_{(\mu,a)}
\times \mathbb{R}^2\times \mathbb{R}^{2*} $,
$h_{(\mu,a)}(\Pi,P,\theta,l)\cdot \pi_{(\mu,a)}=H(A,c,\Pi,P,\theta,l)|_{\mathcal{O}_{(\mu,a)}
\times\mathbb{R}^2\times\mathbb{R}^{2*}}$,
$u_{(\mu,a)}(\Pi,P,\theta,l)\cdot \pi_{(\mu,a)}= u(A,c,\Pi,P,\theta,l)|_{\mathcal{O}_{(\mu,a)}
\times\mathbb{R}^2\times\mathbb{R}^{2*}}$,
and the dynamical vector field of the $R_p$-reduced controlled
underwater vehicle-rotor system satisfies $(4.4)$ and $(4.5)$.
\hskip 0.3cm $\blacksquare$
\end{theo}

\subsection{Hamilton-Jacobi Theorems }

In the following we shall give precisely the geometric constraint conditions
of the $R_p$-reduced symplectic form
$\tilde{\omega}_{\mathcal{O}_{(\mu,a)} \times \mathbb{R}^2 \times
\mathbb{R}^{2*}}^{-}$ for the
dynamical vector field of the regular point reducible controlled
underwater vehicle-rotors system with
coincident centers of buoyancy and gravity,
that is, Type I and Type II of
Hamilton-Jacobi equation for the $R_p$-reduced controlled underwater vehicle-rotor system
$(\mathcal{O}_{(\mu,a)}\times \mathbb{R}^2\times \mathbb{R}^{2*},
\omega^{-}_{\mathcal{O}_{(\mu,a)}\times \mathbb{R}^2\times \mathbb{R}^{2*}},
h_{(\mu,a)}, u_{(\mu,a)}) .$\\

Assume that $\gamma: \textmd{SE}(3)\times \mathbb{R}^2 \rightarrow
T^* (\textmd{SE}(3) \times \mathbb{R}^2) $ is an one-form on
$\textmd{SE}(3) \times \mathbb{R}^2 $, and
$\gamma(A,c,\theta)=(\gamma_1, \cdots, \gamma_{16})$,
and $\gamma$ is closed with respect to $T\pi_{\textmd{SE}(3)\times \mathbb{R}^2}:
TT^* (\textmd{SE}(3) \times \mathbb{R}^2) \rightarrow
T(\textmd{SE}(3) \times \mathbb{R}^2). $
For $(\mu,a) \in \mathfrak{se}^\ast(3),$ the regular value of
$\mathbf{J}_Q: T^* Q \cong \textmd{SE}(3)\times
\mathfrak{se}^\ast(3) \times \mathbb{R}^2 \times \mathbb{R}^{2*} \to
\mathfrak{se}^\ast(3)$,
assume that $\textmd{Im}(\gamma)\subset \mathbf{J}_Q^{-1}(\mu,a), $ and it is
$\textmd{SE}(3)_{(\mu,a)}$-invariant, and $\bar{\gamma}=\pi_{(\mu,a)}(\gamma):
\textmd{SE}(3)\times \mathbb{R}^2 \rightarrow \mathcal{O}_{(\mu,a)} \times
\mathbb{R}^2 \times \mathbb{R}^{2*}$. Denote by
$\bar{\gamma}= (\bar{\gamma}_1,\bar{\gamma}_2,\bar{\gamma}_3,
\bar{\gamma}_4,\bar{\gamma}_5, \bar{\gamma}_6,\bar{\gamma}_7,
\bar{\gamma}_8,\bar{\gamma}_9, \bar{\gamma}_{10}) \in \mathcal{O}_{(\mu,a)} \times
\mathbb{R}^2 \times \mathbb{R}^{2*}(\subset \mathfrak{se}^\ast(3)
\times \mathbb{R}^2 \times \mathbb{R}^{2*}), $ where
$\pi_{(\mu,a)}: \mathbf{J}_Q^{-1}(\mu,a) \rightarrow \mathcal{O}_{(\mu,a)} \times
\mathbb{R}^2 \times \mathbb{R}^{2*}. $ We choose that
$(\Pi,P,\theta,l)\in
\mathcal{O}_{(\mu,a)}\times \mathbb{R}^2 \times \mathbb{R}^{2*}, $ and
$\Pi=(\Pi_1,\Pi_2,\Pi_3)=(\bar{\gamma}_1,\bar{\gamma}_2,\bar{\gamma}_3)$,
$P=(P_1,P_2,P_3)=(\bar{\gamma}_4,\bar{\gamma}_5,\bar{\gamma}_6)$,
$\theta=(\bar{\gamma}_7, \bar{\gamma}_8)$ and $l=(\bar{\gamma}_9,
\bar{\gamma}_{10}) $. Then $h_{(\mu,a)} \cdot \bar{\gamma}:
\textmd{SE}(3) \times \mathbb{R}^2 \rightarrow \mathbb{R} $ is given
by
\begin{align*} & h_{(\mu,a)}(\Pi,P,\theta,l) \cdot \bar{\gamma}=
H(A,c,\Pi,P,\theta,l) |_{\mathcal{O}_{(\mu,a)} \times \mathbb{R}^2 \times \mathbb{R}^2}
\cdot \bar{\gamma}\\ &
=\frac{1}{2}[\frac{(\bar{\gamma}_1-
\bar{\gamma}_9)^2}{\bar{I}_1}+\frac{(\bar{\gamma}_2-
\bar{\gamma}_{10})^2}{\bar{I}_2}
+\frac{\bar{\gamma}_3^2}{\bar{I}_3}+
\frac{\bar{\gamma}_4^2}{m_1}+\frac{\bar{\gamma}_5^2}{m_2}+\frac{\bar{\gamma}_6^2}{m_3}
+\frac{\bar{\gamma}_9^2}{J_1}+\frac{\bar{\gamma}_{10}^2}{J_2}],
\end{align*} and the vector field
\begin{align*}
& X_{h_{(\mu,a)}}(\Pi) \cdot \bar{\gamma}=\{\Pi,\;
h_{(\mu,a)}\}_{-}|_{\mathcal{O}_{(\mu,a)}
\times \mathbb{R}^2\times \mathbb{R}^2}\cdot \bar{\gamma}\\
& = -\Pi\cdot(\nabla_\Pi\Pi\times\nabla_\Pi (h_{(\mu,a)})) \cdot
\bar{\gamma}-P\cdot(\nabla_\Pi\Pi\times\nabla_P
(h_{(\mu,a)})-\nabla_\Pi
(h_{(\mu,a)}) \times\nabla_P \Pi)\cdot \bar{\gamma}\\
& \;\;\;\; + \{\Pi,\;
h_{(\mu,a)}\}_{\mathbb{R}^2}|_{\mathcal{O}_{(\mu,a)} \times
\mathbb{R}^2\times \mathbb{R}^2}\cdot \bar{\gamma}\\
& = -\nabla_\Pi\Pi \cdot(\nabla_\Pi (h_{(\mu,a)}) \times \Pi)\cdot
\bar{\gamma}- \nabla_\Pi\Pi \cdot (\nabla_P
(h_{(\mu,a)}) \times P) \cdot \bar{\gamma}\\
& \;\;\;\; + \sum_{k=1}^2(\frac{\partial \Pi}{\partial \theta_k}
\frac{\partial (h_{(\mu,a)})}{\partial l_k}- \frac{\partial
(h_{(\mu,a)})}{\partial
\theta_k}\frac{\partial \Pi}{\partial l_k})\cdot \bar{\gamma}\\
&=(\Pi_1,\Pi_2,\Pi_3)\times (\frac{(\Pi_1- l_1)}{ \bar{I}_1},
\frac{(\Pi_2- l_2)}{ \bar{I}_2}, \frac{\Pi_3}{\bar{I}_3})\cdot
\bar{\gamma}
+(P_1,P_2,P_3)\times (\frac{P_1}{m_1},\frac{P_2}{m_2},\frac{P_3}{m_3})\cdot \bar{\gamma}\\
&= ( \frac{(\bar{I}_2-\bar{I}_3)\bar{\gamma}_2\bar{\gamma}_3+
\bar{I}_3\bar{\gamma}_3\bar{\gamma}_{10}}{\bar{I}_2\bar{I}_3}
+ \frac{(m_2-m_3)\bar{\gamma}_5\bar{\gamma}_6}{m_2m_3},\\
& \;\;\;\;\;\; \frac{(\bar{I}_3-\bar{I}_1)\bar{\gamma}_3
\bar{\gamma}_1-
\bar{I}_3\bar{\gamma}_3\bar{\gamma}_9}{\bar{I}_3\bar{I}_1}
+ \frac{(m_3-m_1)\bar{\gamma}_6\bar{\gamma}_4}{m_3m_1},\\
& \;\;\;\;\;\;
\frac{(\bar{I}_1-\bar{I}_2)\bar{\gamma}_1\bar{\gamma}_2-\bar{I}_1\bar{\gamma}_1
\bar{\gamma}_{10}+
\bar{I}_2\bar{\gamma}_2\bar{\gamma}_9}{\bar{I}_1\bar{I}_2}
+  \frac{(m_1-m_2)\bar{\gamma}_4\bar{\gamma}_5}{m_1m_2}),
\end{align*}
since $\nabla_{\Pi_i}\Pi_i=1,\; \nabla_{\Pi_i}\Pi_j=0, \; i\neq j, \;  \nabla_{\Pi_i}P_j=\nabla_{P_i}\Pi_j=0$
and $ \nabla_{P_i} (h_{(\mu,a)})= \frac{P_i}{m_i}, \; \nabla_{\Pi_k} (h_{(\mu,a)}
)= (\Pi_k-l_k)/\bar{I}_k ,\; \nabla_{\Pi_3}
(h_{(\mu,a)})= \Pi_3/\bar{I}_3 , \; \frac{\partial
\Pi}{\partial \theta_k}= \frac{\partial (h_{(\mu,a)})}{\partial
\theta_k}=0, \; i, j=1,2,3, \; k=1,2. $

\begin{align*}
& X_{h_{(\mu,a)}}(P) \cdot \bar{\gamma}=\{P,\;
h_{(\mu,a)}\}_{-}|_{\mathcal{O}_{(\mu,a)} \times \mathbb{R}^2\times
\mathbb{R}^2}\cdot \bar{\gamma}\\
& =-\Pi\cdot(\nabla_\Pi P \times\nabla_\Pi (h_{(\mu,a)})) \cdot
\bar{\gamma}-P\cdot(\nabla_\Pi P\times\nabla_P
(h_{(\mu,a)})-\nabla_\Pi
(h_{(\mu,a)}) \times\nabla_P P)\cdot \bar{\gamma}\\
& \;\;\;\; + \{P,\;
h_{(\mu,a)}\}_{\mathbb{R}^2}|_{\mathcal{O}_{(\mu,a)} \times
\mathbb{R}^2\times \mathbb{R}^2}\cdot \bar{\gamma}\\
& =\nabla_P P \cdot(P \times\nabla_\Pi
(h_{(\mu,a)}))\cdot \bar{\gamma}+ \sum_{k=1}^2(\frac{\partial
P}{\partial \theta_k} \frac{\partial (h_{(\mu,a)})}{\partial
l_k}- \frac{\partial (h_{(\mu,a)})}{\partial
\theta_k}\frac{\partial P}{\partial l_k})\cdot \bar{\gamma}\\
&=(P_1,P_2,P_3)\times (\frac{(\Pi_1- l_1)}{
\bar{I}_1}, \frac{(\Pi_2- l_2)}{ \bar{I}_2},
\frac{\Pi_3}{\bar{I}_3})\cdot \bar{\gamma}\\
&= ( \frac{\bar{I}_2\bar{\gamma}_5\bar{\gamma}_3-\bar{I}_3\bar{\gamma}_6\bar{\gamma}_2+
\bar{I}_3\bar{\gamma}_3\bar{\gamma}_{10}}{\bar{I}_2\bar{I}_3}, \;\;
\frac{\bar{I}_3\bar{\gamma}_6\bar{\gamma}_1-\bar{I}_1\bar{\gamma}_4\bar{\gamma}_3-
\bar{I}_3\bar{\gamma}_3\bar{\gamma}_9}{\bar{I}_3\bar{I}_1}, \\
& \;\;\;\;\;\;
\frac{\bar{I}_1\bar{\gamma}_4\bar{\gamma}_2-\bar{I}_2\bar{\gamma}_5\bar{\gamma}_1-
\bar{I}_1\bar{\gamma}_1\bar{\gamma}_{10} + \bar{I}_2\bar{\gamma}_2\bar{\gamma}_9}{\bar{I}_1\bar{I}_2} ),
\end{align*}
since $\nabla_{P_i}P_i =1, \; \nabla_{P_i}P_j=0, \; i\neq j, \; \nabla_{\Pi_i}P_j =0, $ and
$\nabla_{\Pi_k} (h_{(\mu,a)})= (\Pi_k-l_k)/\bar{I}_k ,\; \nabla_{\Pi_3}
(h_{(\mu,a)})= \Pi_3/\bar{I}_3 , \; \frac{\partial
P_j}{\partial \theta_k}= \frac{\partial (h_{(\mu,a)})}{\partial
\theta_k}=0, \; i, j= 1,2,3, \; k=1,2.$

\begin{align*}
& X_{h_{(\mu,a)}}(\theta) \cdot \bar{\gamma}=\{\theta,\;
h_{(\mu,a)}\}_{-}|_{\mathcal{O}_{(\mu,a)} \times \mathbb{R}^2\times
\mathbb{R}^2}\cdot \bar{\gamma}\\
& =-\Pi\cdot(\nabla_\Pi \theta \times\nabla_\Pi (h_{(\mu,a)}))\cdot
\bar{\gamma}-P \cdot(\nabla_\Pi \theta \times\nabla_P
(h_{(\mu,a)})-\nabla_\Pi
(h_{(\mu,a)}) \times\nabla_P \theta)\cdot \bar{\gamma}\\
& \;\;\;\; + \{\theta,\;
h_{(\mu,a)}\}_{\mathbb{R}^2}|_{\mathcal{O}_{(\mu,a)} \times
\mathbb{R}^2\times \mathbb{R}^2}\cdot \bar{\gamma}\\
& = \sum_{k=1}^2(\frac{\partial \theta}{\partial \theta_k}
\frac{\partial (h_{(\mu,a)})}{\partial l_k}- \frac{\partial
(h_{(\mu,a)})}{\partial
\theta_k}\frac{\partial \theta}{\partial l_k})\cdot \bar{\gamma}\\
& = (-\frac{(\bar{\gamma}_1- \bar{\gamma}_9)}{ \bar{I}_1}+
\frac{\bar{\gamma}_9}{J_1}, \; -\frac{(\bar{\gamma}_2-
\bar{\gamma}_{10})}{ \bar{I}_2}+ \frac{\bar{\gamma}_{10}}{J_2}),
\end{align*}
since $\nabla_{\Pi_i}\theta=\nabla_{P_i}\theta=0, $ $\frac{\partial \theta_k}{\partial
\theta_k}= 1, \; \frac{\partial \theta_n}{\partial \theta_k}=
0, \; n\neq k, \;\; \frac{\partial (h_{(\mu,a)})}{\partial \theta_k}=0, $
and $\frac{\partial (h_{(\mu,a)})}{\partial l_k}= -(\Pi_k-l_k)/\bar{I}_k
+\frac{l_k}{J_k}, \; i= 1,2,3, \; k,n=1,2. $

\begin{align*}
& X_{h_{(\mu,a)}}(l) \cdot \bar{\gamma}=\{l,\;
h_{(\mu,a)}\}_{-}|_{\mathcal{O}_{(\mu,a)} \times \mathbb{R}^2\times
\mathbb{R}^2}\cdot \bar{\gamma}\\
& =-\Pi\cdot(\nabla_\Pi l \times\nabla_\Pi (h_{(\mu,a)})) \cdot
\bar{\gamma}-P \cdot(\nabla_\Pi l \times\nabla_P
(h_{(\mu,a)})-\nabla_\Pi
(h_{(\mu,a)}) \times\nabla_P l)\cdot \bar{\gamma}\\
& \;\;\;\; + \{l,\;
h_{(\mu,a)}\}_{\mathbb{R}^2}|_{\mathcal{O}_{(\mu,a)} \times
\mathbb{R}^2\times \mathbb{R}^2}\cdot \bar{\gamma}\\
& = \sum_{k=1}^2(\frac{\partial l}{\partial \theta_k} \frac{\partial
(h_{(\mu,a)})}{\partial l_k}- \frac{\partial (h_{(\mu,a)})}{\partial
\theta_k}\frac{\partial l}{\partial l_k})\cdot \bar{\gamma}=(0,0),
\end{align*}
since $\nabla_{\Pi_i} l=\nabla_{P_i} l=0,$ and $\frac{\partial l}{\partial \theta_k}=
\frac{\partial (h_{(\mu,a)})}{\partial \theta_k}=0, \; i=1,2,3, \; k=1,2. $\\

On the other hand, from the expressions of the dynamical vector field $\tilde{X}$
and Hamiltonian vector field $X_H$, we have that
\begin{align*}
\tilde{X}(\Pi, P, \theta, l)^\gamma & =T\pi_{\textmd{SE}(3)\times \mathbb{R}^2 }\cdot \tilde{X}\cdot\gamma(\Pi, P, \theta, l)\\
& =T\pi_{\textmd{SE}(3)\times \mathbb{R}^2 }\cdot (X_H+ \textnormal{vlift}(u))\cdot\gamma (\Pi, P, \theta, l)\\
&=T\pi_{\textmd{SE}(3)\times \mathbb{R}^2 }\cdot X_H \cdot\gamma (\Pi, P, \theta, l) = X_H\cdot\gamma(\Pi, P, \theta, l),
\end{align*}
that is,
\begin{align*}
\tilde{X}(\Pi)^\gamma & = X_H(\Pi)\cdot\gamma \\
& = ( \frac{(\bar{I}_2-\bar{I}_3)\gamma_8\gamma_9+
\bar{I}_3\gamma_9\gamma_{16}}{\bar{I}_2\bar{I}_3}
+ \frac{(m_2-m_3)\gamma_{11}\gamma_{12}}{m_2m_3},\\
& \;\;\;\;\;\;
\frac{(\bar{I}_3-\bar{I}_1)\gamma_9\gamma_7-
\bar{I}_3\gamma_9\gamma_{15}}{\bar{I}_3\bar{I}_1}
+ \frac{(m_3-m_1)\gamma_{12}\gamma_{10}}{m_3m_1}, \\
& \;\;\;\;\;\;
\frac{(\bar{I}_1-\bar{I}_2)\gamma_7\gamma_8-
\bar{I}_1\gamma_7\gamma_{16} + \bar{I}_2\gamma_8\gamma_{15}}{\bar{I}_1\bar{I}_2}
+ \frac{(m_1-m_2)\gamma_{10}\gamma_{11}}{m_1m_2},
\end{align*}
\begin{align*}
\tilde{X}(P)^\gamma & = X_H(P)\cdot\gamma \\
&= ( \frac{\bar{I}_2\gamma_{11}\gamma_9-\bar{I}_3\gamma_{12}\gamma_8+
\bar{I}_3\gamma_9\gamma_{16}}{\bar{I}_2\bar{I}_3}, \;\;
\frac{\bar{I}_3\gamma_{12}\gamma_7-\bar{I}_1\gamma_{10}\gamma_9-
\bar{I}_3\gamma_9\gamma_{15}}{\bar{I}_3\bar{I}_1}, \\
& \;\;\;\;\;\;
\frac{\bar{I}_1\gamma_{10}\gamma_8-\bar{I}_2\gamma_{11}\gamma_7-
\bar{I}_1\gamma_7\gamma_{16} + \bar{I}_2\gamma_8\gamma_{15}}{\bar{I}_1\bar{I}_2} ),
\end{align*}
\begin{align*}
\tilde{X}(\theta)^\gamma & = X_H(\theta)\cdot\gamma
 = ( -\frac{(\gamma_7- \gamma_{15})}{ \bar{I}_1}
+\frac{\gamma_{15}}{J_1},\;\; -\frac{(\gamma_8- \gamma_{16})}{ \bar{I}_2}
+\frac{\gamma_{16}}{J_2} ),
\end{align*}
\begin{align*}
\tilde{X}(l)^\gamma & = X_H(l)\cdot\gamma=(0,0).
\end{align*}
Since $\gamma$ is closed with respect to
$T\pi_{(\textmd{SE}(3) \times \mathbb{R}^2)}: TT^* (\textmd{SE}(3) \times \mathbb{R}^2)
\rightarrow T(\textmd{SE}(3) \times \mathbb{R}^2), $
then $\pi_{(\textmd{SE}(3) \times \mathbb{R}^2)}^*(\mathbf{d}\gamma)=0.$ We choose that
$(\gamma_7,\gamma_8,\gamma_9)=\Pi=(\Pi_1,\Pi_2,\Pi_3)=
(\bar{\gamma}_1,\bar{\gamma}_2,\bar{\gamma}_3), \;
(\gamma_{10},\gamma_{11},\gamma_{12})=P=(P_1,P_2,P_3)=
(\bar{\gamma}_4,\bar{\gamma}_5,\bar{\gamma}_6), $
and $(\gamma_{13},\gamma_{14})= \theta= (\theta_1,\theta_2)
=(\bar{\gamma}_7,\bar{\gamma}_8),
\; (\gamma_{15},\gamma_{16})= l=(l_1,l_2)
=(\bar{\gamma}_9,\bar{\gamma}_{10}). $
Hence
\begin{align*}
& T\bar{\gamma}\cdot \tilde{X}(\Pi)^\gamma
= X_{h_{(\mu,a)}}(\Pi) \cdot \bar{\gamma}, \;\;\;\;\;\;
T\bar{\gamma}\cdot \tilde{X}(P)^\gamma
= X_{h_{(\mu,a)}}(P) \cdot \bar{\gamma}, \\
& T\bar{\gamma}\cdot \tilde{X}(\theta)^\gamma
= X_{h_{(\mu,a)}}(\theta) \cdot \bar{\gamma}, \;\;\;\;\;\;
T\bar{\gamma}\cdot \tilde{X}(l)^\gamma
= X_{h_{(\mu,a)}}(l) \cdot \bar{\gamma}.
\end{align*}
Thus, the Type I of Hamilton-Jacobi equation for the
$R_p$-reduced controlled underwater vehicle-rotor system
$(\mathcal{O}_{(\mu,a)}\times \mathbb{R}^2\times \mathbb{R}^{2*},
\omega^{-}_{\mathcal{O}_{(\mu,a)}\times \mathbb{R}^2\times \mathbb{R}^{2*}},
h_{(\mu,a)}, u_{(\mu,a)})$ holds.\\

Next, for $(\mu,a) \in \mathfrak{se}^\ast(3),$
the regular value of $\mathbf{J}_Q$, and a $\textmd{SE}(3)_{(\mu,a)}$-invariant symplectic map
$\varepsilon: T^* (\textmd{SE}(3)\times \mathbb{R}^2) \rightarrow T^*(\textmd{SE}(3)\times \mathbb{R}^2),$
assume that $\varepsilon(A,c,\Pi,P, \theta, l) =(\varepsilon_1,\cdots,
\varepsilon_{16}),$ and
$\varepsilon(\mathbf{J}^{-1}((\mu,a)))\subset
\mathbf{J}^{-1}((\mu,a)). $ Denote by
$\bar{\varepsilon}=\pi_{(\mu,a)}(\varepsilon):
\mathbf{J}^{-1}((\mu,a)) \rightarrow \mathcal{O}_{(\mu,a)}\times \mathbb{R}^2\times \mathbb{R}^{2*}, $
and $\bar{\varepsilon}=(\bar{\varepsilon}_1,\cdots,\bar{\varepsilon}_{10})
\in \mathcal{O}_{(\mu,a)}\times \mathbb{R}^2\times \mathbb{R}^{2*}
(\subset \mathfrak{se}^\ast(3)\times \mathbb{R}^2\times \mathbb{R}^{2*}), $
and $\lambda= \gamma \cdot \pi_{\textmd{SE}(3)\times \mathbb{R}^2}: T^* (\textmd{SE}(3)\times \mathbb{R}^2)
\rightarrow T^* (\textmd{SE}(3)\times \mathbb{R}^2),$ and $\lambda(A,c, \Pi,P, \theta, l)
=(\lambda_1,\cdots, \lambda_{16}),$ and
$\bar{\lambda}=\pi_{(\mu,a)}(\lambda): T^* (\textmd{SE}(3)\times \mathbb{R}^2)
\rightarrow \mathcal{O}_{(\mu,a)}\times \mathbb{R}^2\times \mathbb{R}^{2*}, $ and
$\bar{\lambda}=
(\bar{\lambda}_1,\cdots,\bar{\lambda}_{10}) \in \mathcal{O}_{(\mu,a)}\times \mathbb{R}^2\times \mathbb{R}^{2*}
(\subset \mathfrak{se}^\ast(3)\times \mathbb{R}^2\times \mathbb{R}^{2*}). $ We choose that
$(\Pi,P,\theta,l)\in
\mathcal{O}_{(\mu,a)}\times \mathbb{R}^2 \times \mathbb{R}^{2*}, $ and
$\Pi=(\Pi_1,\Pi_2,\Pi_3)
=(\bar{\varepsilon}_1,\bar{\varepsilon}_2,\bar{\varepsilon}_3),
\; P= (P_1,P_2,P_3)=
(\bar{\varepsilon}_4,\bar{\varepsilon}_5,\bar{\varepsilon}_6), $
$\theta= (\theta_1,\theta_2)=(\bar{\varepsilon}_7, \bar{\varepsilon}_8)$ and
$ l=(l_1,l_2)=(\bar{\varepsilon}_9, \bar{\varepsilon}_{10}) $.
Then $h_{(\mu,a)} \cdot \bar{\varepsilon}:
T^*(\textmd{SE}(3)\times \mathbb{R}^2) \rightarrow \mathbb{R} $ is given by
\begin{align*}
& h_{(\mu,a)}(\Pi,P,\theta,l) \cdot \bar{\varepsilon}=
H(A,c,\Pi,P, \theta, l)|_{\mathcal{O}_{(\mu,a)}\times \mathbb{R}^2\times \mathbb{R}^{2*}} \cdot \bar{\varepsilon}\\
& =\frac{1}{2}[\frac{(\bar{\varepsilon}_1-
\bar{\varepsilon}_9)^2}{\bar{I}_1}+\frac{(\bar{\varepsilon}_2-
\bar{\varepsilon}_{10})^2}{\bar{I}_2}
+\frac{\bar{\varepsilon}_3^2}{\bar{I}_3}+
\frac{\bar{\varepsilon}_4^2}{m_1}+\frac{\bar{\varepsilon}_5^2}{m_2}+\frac{\bar{\varepsilon}_6^2}{m_3}
+\frac{\bar{\varepsilon}_9^2}{J_1}+\frac{\bar{\varepsilon}_{10}^2}{J_2}],
\end{align*} and the vector field
\begin{align*}
& X_{h_{(\mu,a)}}(\Pi) \cdot \bar{\varepsilon}
= \{\Pi,h_{(\mu,a)}\}_{-}|_{\mathcal{O}_{(\mu,a)}\times \mathbb{R}^2\times \mathbb{R}^{2*}}\cdot \bar{\varepsilon}\\
&=(\Pi_1,\Pi_2,\Pi_3)\times (\frac{(\Pi_1- l_1)}{ \bar{I}_1},
\frac{(\Pi_2- l_2)}{ \bar{I}_2}, \frac{\Pi_3}{\bar{I}_3})\cdot
\bar{\varepsilon}
+(P_1,P_2,P_3)\times (\frac{P_1}{m_1},\frac{P_2}{m_2},\frac{P_3}{m_3})\cdot \bar{\varepsilon}\\
&= ( \frac{(\bar{I}_2-\bar{I}_3)\bar{\varepsilon}_2\bar{\varepsilon}_3+
\bar{I}_3\bar{\varepsilon}_3\bar{\varepsilon}_{10}}{\bar{I}_2\bar{I}_3}
+\frac{(m_2-m_3)\bar{\varepsilon}_5\bar{\varepsilon}_6}{m_2m_3}, \\
& \;\;\;\;\;\; \frac{(\bar{I}_3-\bar{I}_1)\bar{\varepsilon}_3
\bar{\varepsilon}_1-
\bar{I}_3\bar{\varepsilon}_3\bar{\varepsilon}_9}{\bar{I}_3\bar{I}_1}
+\frac{(m_3-m_1)\bar{\varepsilon}_6\bar{\varepsilon}_4}{m_3m_1},\\
& \;\;\;\;\;\;
\frac{(\bar{I}_1-\bar{I}_2)\bar{\varepsilon}_1\bar{\varepsilon}_2-\bar{I}_1\bar{\varepsilon}_1
\bar{\varepsilon}_{10}+
\bar{I}_2\bar{\varepsilon}_2\bar{\varepsilon}_9}{\bar{I}_1\bar{I}_2}
+\frac{(m_1-m_2)\bar{\varepsilon}_4\bar{\varepsilon}_5}{m_1m_2},
\end{align*}
\begin{align*}
& X_{h_{(\mu,a)}}(P) \cdot
\bar{\varepsilon}
= \{P,h_{(\mu,a)}\}_{-}|_{\mathcal{O}_{(\mu,a)}\times \mathbb{R}^2\times \mathbb{R}^{2*}}\cdot \bar{\varepsilon}\\
&=(P_1,P_2,P_3)\times (\frac{(\Pi_1- l_1)}{
\bar{I}_1}, \frac{(\Pi_2- l_2)}{ \bar{I}_2},
\frac{\Pi_3}{\bar{I}_3})\cdot \bar{\varepsilon}\\
&= ( \frac{\bar{I}_2\bar{\varepsilon}_5\bar{\varepsilon}_3-\bar{I}_3\bar{\varepsilon}_6\bar{\varepsilon}_2+
\bar{I}_3\bar{\varepsilon}_3\bar{\varepsilon}_{10}}{\bar{I}_2\bar{I}_3}, \;\;
\frac{\bar{I}_3\bar{\varepsilon}_6\bar{\varepsilon}_1-\bar{I}_1\bar{\varepsilon}_4\bar{\varepsilon}_3-
\bar{I}_3\bar{\varepsilon}_3\bar{\varepsilon}_9}{\bar{I}_3\bar{I}_1}, \\
& \;\;\;\;\;\;
\frac{\bar{I}_1\bar{\varepsilon}_4\bar{\varepsilon}_2-\bar{I}_2\bar{\varepsilon}_5\bar{\varepsilon}_1-
\bar{I}_1\bar{\varepsilon}_1\bar{\varepsilon}_{10} + \bar{I}_2\bar{\varepsilon}_2\bar{\varepsilon}_9}{\bar{I}_1\bar{I}_2} ),
\end{align*}
\begin{align*}
& X_{h_{(\mu,a)}}(\theta) \cdot \bar{\varepsilon}=\{\theta,\;
h_{(\mu,a)}\}_{-}|_{\mathcal{O}_{(\mu,a)} \times \mathbb{R}^2\times
\mathbb{R}^2}\cdot \bar{\varepsilon}\\
& = \sum_{k=1}^2(\frac{\partial \theta}{\partial \theta_k}
\frac{\partial (h_{(\mu,a)})}{\partial l_k}- \frac{\partial
(h_{(\mu,a)})}{\partial
\theta_k}\frac{\partial \theta}{\partial l_k})\cdot \bar{\varepsilon}
= (-\frac{(\bar{\varepsilon}_1- \bar{\varepsilon}_9)}{ \bar{I}_1}+
\frac{\bar{\varepsilon}_9}{J_1}, \; -\frac{(\bar{\varepsilon}_2-
\bar{\varepsilon}_{10})}{ \bar{I}_2}+ \frac{\bar{\varepsilon}_{10}}{J_2}),
\end{align*}
\begin{align*}
& X_{h_{(\mu,a)}}(l) \cdot \bar{\varepsilon}=\{l,\;
h_{(\mu,a)}\}_{-}|_{\mathcal{O}_{(\mu,a)} \times \mathbb{R}^2\times
\mathbb{R}^2}\cdot \bar{\varepsilon}
= \sum_{k=1}^2(\frac{\partial l}{\partial \theta_k} \frac{\partial
(h_{(\mu,a)})}{\partial l_k}- \frac{\partial (h_{(\mu,a)})}{\partial
\theta_k}\frac{\partial l}{\partial l_k})\cdot \bar{\varepsilon}=(0,0).
\end{align*}
On the other hand, from the expressions of the dynamical vector field $\tilde{X}$
and Hamiltonian vector field $X_H$, we have that
\begin{align*}
\tilde{X}(\Pi, P, \theta, l)^\varepsilon & =T\pi_{\textmd{SE}(3)\times \mathbb{R}^2 }\cdot \tilde{X}\cdot\varepsilon(\Pi, P, \theta, l)\\
& =T\pi_{\textmd{SE}(3)\times \mathbb{R}^2 }\cdot (X_H+ \textnormal{vlift}(u))\cdot\varepsilon (\Pi, P, \theta, l)\\
& =T\pi_{\textmd{SE}(3)\times \mathbb{R}^2 }\cdot X_H \cdot\varepsilon (\Pi, P, \theta, l)= X_H\cdot\varepsilon(\Pi, P, \theta, l),
\end{align*}
that is,
\begin{align*}
\tilde{X}(\Pi)^\varepsilon & = X_H(\Pi)\cdot\varepsilon \\
& = ( \frac{(\bar{I}_2-\bar{I}_3)\varepsilon_8\varepsilon_9+
\bar{I}_3\varepsilon_9\varepsilon_{16}}{\bar{I}_2\bar{I}_3}
+ \frac{(m_2-m_3)\varepsilon_{11}\varepsilon_{12}}{m_2m_3},\\
& \;\;\;\;\;\;
\frac{(\bar{I}_3-\bar{I}_1)\varepsilon_9\varepsilon_7-
\bar{I}_3\varepsilon_9\varepsilon_{15}}{\bar{I}_3\bar{I}_1}
+ \frac{(m_3-m_1)\varepsilon_{12}\varepsilon_{10}}{m_3m_1}, \\
& \;\;\;\;\;\;
\frac{(\bar{I}_1-\bar{I}_2)\varepsilon_7\varepsilon_8-
\bar{I}_1\varepsilon_7\varepsilon_{16} + \bar{I}_2\varepsilon_8\varepsilon_{15}}{\bar{I}_1\bar{I}_2}
 +  \frac{(m_1-m_2)\varepsilon_{10}\varepsilon_{11}}{m_1m_2},
\end{align*}
\begin{align*}
\tilde{X}(P)^\varepsilon & = X_H(P)\cdot\varepsilon \\
&= ( \frac{\bar{I}_2\varepsilon_{11}\varepsilon_9-\bar{I}_3\varepsilon_{12}\varepsilon_8+
\bar{I}_3\varepsilon_9\varepsilon_{16}}{\bar{I}_2\bar{I}_3}, \;\;
\frac{\bar{I}_3\varepsilon_{12}\varepsilon_7-\bar{I}_1\varepsilon_{10}\varepsilon_9-
\bar{I}_3\varepsilon_9\varepsilon_{15}}{\bar{I}_3\bar{I}_1}, \\
& \;\;\;\;\;\;
\frac{\bar{I}_1\varepsilon_{10}\varepsilon_8-\bar{I}_2\varepsilon_{11}\varepsilon_7-
\bar{I}_1\varepsilon_7\varepsilon_{16} + \bar{I}_2\varepsilon_8\varepsilon_{15}}{\bar{I}_1\bar{I}_2} ),
\end{align*}
\begin{align*}
\tilde{X}(\theta)^\varepsilon & = X_H(\theta)\cdot\varepsilon
= ( -\frac{(\varepsilon_7- \varepsilon_{15})}{ \bar{I}_1}
+\frac{\varepsilon_{15}}{J_1},\;\; -\frac{(\varepsilon_8- \varepsilon_{16})}{ \bar{I}_2}
+\frac{\varepsilon_{16}}{J_2} ),
\end{align*}
\begin{align*}
\tilde{X}(l)^\varepsilon & = X_H(l)\cdot\varepsilon =(0,0),
\end{align*}
then we have that
\begin{align*}
T\bar{\gamma}\cdot \tilde{X}(\Pi)^\varepsilon
&= ( \frac{(\bar{I}_2-\bar{I}_3)\bar{\gamma}_2\bar{\gamma}_3+
\bar{I}_3\bar{\gamma}_3\bar{\gamma}_{10}}{\bar{I}_2\bar{I}_3}
+  \frac{(m_2-m_3)\bar{\gamma}_5\bar{\gamma}_6}{m_2m_3},\\
& \;\;\;\;\;\; \frac{(\bar{I}_3-\bar{I}_1)\bar{\gamma}_3
\bar{\gamma}_1-
\bar{I}_3\bar{\gamma}_3\bar{\gamma}_9}{\bar{I}_3\bar{I}_1}
+  \frac{(m_3-m_1)\bar{\gamma}_6\bar{\gamma}_4}{m_3m_1},\\
& \;\;\;\;\;\;
\frac{(\bar{I}_1-\bar{I}_2)\bar{\gamma}_1\bar{\gamma}_2-\bar{I}_1\bar{\gamma}_1
\bar{\gamma}_{10}+
\bar{I}_2\bar{\gamma}_2\bar{\gamma}_9}{\bar{I}_1\bar{I}_2}
+  \frac{(m_1-m_2)\bar{\gamma}_4\bar{\gamma}_5}{m_1m_2},
\end{align*}
\begin{align*}
T\bar{\gamma}\cdot \tilde{X}(P)^\varepsilon
&= ( \frac{\bar{I}_2\bar{\gamma}_5\bar{\gamma}_3-\bar{I}_3\bar{\gamma}_6\bar{\gamma}_2+
\bar{I}_3\bar{\gamma}_3\bar{\gamma}_{10}}{\bar{I}_2\bar{I}_3}, \;\;
\frac{\bar{I}_3\bar{\gamma}_6\bar{\gamma}_1-\bar{I}_1\bar{\gamma}_4\bar{\gamma}_3-
\bar{I}_3\bar{\gamma}_3\bar{\gamma}_9}{\bar{I}_3\bar{I}_1}, \\
& \;\;\;\;\;\;
\frac{\bar{I}_1\bar{\gamma}_4\bar{\gamma}_2-\bar{I}_2\bar{\gamma}_5\bar{\gamma}_1-
\bar{I}_1\bar{\gamma}_1\bar{\gamma}_{10}
+ \bar{I}_2\bar{\gamma}_2\bar{\gamma}_9}{\bar{I}_1\bar{I}_2} ),
\end{align*}
\begin{align*}
T\bar{\gamma}\cdot \tilde{X}(\theta)^\varepsilon
& = (-\frac{(\bar{\gamma}_1- \bar{\gamma}_9)}{ \bar{I}_1}+
\frac{\bar{\gamma}_9}{J_1}, \; -\frac{(\bar{\gamma}_2-
\bar{\gamma}_{10})}{ \bar{I}_2}+ \frac{\bar{\gamma}_{10}}{J_2}),
\end{align*}
\begin{align*}
T\bar{\gamma}\cdot \tilde{X}(l)^\varepsilon=(0,0).
\end{align*}
Note that $$T\bar{\lambda}\cdot \tilde{X} \cdot \varepsilon=T\pi_{(\mu,a)}\cdot T\lambda \cdot (X_H+ \textnormal{vlift}(u))\cdot\varepsilon
=T\pi_{(\mu,a)}\cdot T\gamma \cdot T\pi_{\textmd{SE}(3)\times \mathbb{R}^2}\cdot (X_H+ \textnormal{vlift}(u))\cdot\varepsilon
=T\bar{\lambda}\cdot X_H \cdot \varepsilon,$$ that is,
\begin{align*}
T\bar{\lambda}\cdot \tilde{X}(\Pi) \cdot \varepsilon & =T\bar{\lambda}\cdot X_H(\Pi) \cdot \varepsilon\\
&= ( \frac{(\bar{I}_2-\bar{I}_3)\bar{\lambda}_2\bar{\lambda}_3+
\bar{I}_3\bar{\lambda}_3\bar{\lambda}_{10}}{\bar{I}_2\bar{I}_3}
+  \frac{(m_2-m_3)\bar{\lambda}_5\bar{\lambda}_6}{m_2m_3},\\
& \;\;\;\;\;\; \frac{(\bar{I}_3-\bar{I}_1)\bar{\lambda}_3
\bar{\lambda}_1-
\bar{I}_3\bar{\lambda}_3\bar{\lambda}_9}{\bar{I}_3\bar{I}_1}
+  \frac{(m_3-m_1)\bar{\lambda}_6\bar{\lambda}_4}{m_3m_1},\\
& \;\;\;\;\;\;
\frac{(\bar{I}_1-\bar{I}_2)\bar{\lambda}_1\bar{\lambda}_2-\bar{I}_1\bar{\lambda}_1
\bar{\lambda}_{10}+
\bar{I}_2\bar{\lambda}_2\bar{\lambda}_9}{\bar{I}_1\bar{I}_2}
+  \frac{(m_1-m_2)\bar{\lambda}_4\bar{\lambda}_5}{m_1m_2},
\end{align*}
\begin{align*}
T\bar{\lambda}\cdot \tilde{X}(P) \cdot \varepsilon & =T\bar{\lambda}\cdot X_H(P) \cdot \varepsilon\\
&= ( \frac{\bar{I}_2\bar{\lambda}_5\bar{\lambda}_3-\bar{I}_3\bar{\lambda}_6\bar{\lambda}_2+
\bar{I}_3\bar{\lambda}_3\bar{\lambda}_{10}}{\bar{I}_2\bar{I}_3}, \;\;
\frac{\bar{I}_3\bar{\lambda}_6\bar{\lambda}_1-\bar{I}_1\bar{\lambda}_4\bar{\lambda}_3-
\bar{I}_3\bar{\lambda}_3\bar{\lambda}_9}{\bar{I}_3\bar{I}_1}, \\
& \;\;\;\;\;\;
\frac{\bar{I}_1\bar{\lambda}_4\bar{\lambda}_2-\bar{I}_2\bar{\lambda}_5\bar{\lambda}_1-
\bar{I}_1\bar{\lambda}_1\bar{\lambda}_{10}
+ \bar{I}_2\bar{\lambda}_2\bar{\lambda}_9}{\bar{I}_1\bar{I}_2} ),
\end{align*}
\begin{align*}
T\bar{\lambda}\cdot \tilde{X}(\theta) \cdot \varepsilon=T\bar{\lambda}\cdot X_H(\theta) \cdot \varepsilon
= (-\frac{(\bar{\lambda}_1- \bar{\lambda}_9)}{ \bar{I}_1}+
\frac{\bar{\lambda}_9}{J_1}, \; -\frac{(\bar{\lambda}_2-
\bar{\lambda}_{10})}{ \bar{I}_2}+ \frac{\bar{\lambda}_{10}}{J_2}),
\end{align*}
\begin{align*}
T\bar{\lambda}\cdot \tilde{X}(l) \cdot \varepsilon=T\bar{\lambda}\cdot X_H(l) \cdot \varepsilon=(0,0).
\end{align*}
Thus, when we choose that
$(\Pi,P,\theta,l)\in
\mathcal{O}_{(\mu,a)}\times \mathbb{R}^2 \times \mathbb{R}^{2*}, $ and
 $(\varepsilon_7,\varepsilon_8,\varepsilon_9)=\Pi=(\Pi_1,\Pi_2,\Pi_3)
 =(\bar{\gamma}_1,\bar{\gamma}_2,\bar{\gamma}_3)=
(\bar{\varepsilon}_1,\bar{\varepsilon}_2,\bar{\varepsilon}_3)=
(\bar{\lambda}_1,\bar{\lambda}_2,\bar{\lambda}_3), $ and
$(\varepsilon_{10},\varepsilon_{11},\varepsilon_{12})=P=(P_1,P_2,P_3)
=(\bar{\gamma}_4,\bar{\gamma}_5,\bar{\gamma}_6)=
(\bar{\varepsilon}_4,\bar{\varepsilon}_5,\bar{\varepsilon}_6)=
(\bar{\lambda}_4,\bar{\lambda}_5,\bar{\lambda}_6), $
and $(\varepsilon_{13},\varepsilon_{14})=\theta
= (\theta_1,\theta_2)=(\bar{\gamma}_7,\bar{\gamma}_8)
=(\bar{\varepsilon}_7,\bar{\varepsilon}_8)=(\bar{\lambda}_7,\bar{\lambda}_8),
\; (\varepsilon_{15},\varepsilon_{16})=l=(l_1,l_2)
=(\bar{\gamma}_9,\bar{\gamma}_{10})
=(\bar{\varepsilon}_9,\bar{\varepsilon}_{10})
=(\bar{\lambda}_9,\bar{\lambda}_{10}). $ we must have that
\begin{align*}
& T\bar{\gamma}\cdot \tilde{X}(\Pi)^\varepsilon =X_{h_{(\mu,a)}}(\Pi) \cdot \bar{\varepsilon}
=T\bar{\lambda}\cdot \tilde{X}(\Pi) \cdot \varepsilon, \\
& T\bar{\gamma}\cdot \tilde{X}(P)^\varepsilon=X_{h_{(\mu,a)}}(P) \cdot \bar{\varepsilon}
=T\bar{\lambda}\cdot \tilde{X}(P) \cdot \varepsilon,\\
& T\bar{\gamma}\cdot \tilde{X}(\theta)^\varepsilon =X_{h_{(\mu,a)}}(\theta) \cdot \bar{\varepsilon}
=T\bar{\lambda}\cdot \tilde{X}(\theta) \cdot \varepsilon, \\
& T\bar{\gamma}\cdot \tilde{X}(l)^\varepsilon =X_{h_{(\mu,a)}}(l) \cdot \bar{\varepsilon}
=T\bar{\lambda}\cdot \tilde{X}(l) \cdot \varepsilon.
\end{align*}
Since the map $\varepsilon: T^* (\textmd{SE}(3)\times \mathbb{R}^2)
\rightarrow T^* (\textmd{SE}(3)\times \mathbb{R}^2)$ is symplectic, then
$T\bar{\varepsilon}\cdot X_{h_{(\mu,a)} \cdot \bar{\varepsilon}}
=X_{h_{(\mu,a)}} \cdot \bar{\varepsilon}. $
Thus, in this case, we must have that
$\varepsilon$ and $\bar{\varepsilon} $ are the solution of the Type II of
Hamilton-Jacobi equation
$T\bar{\gamma}\cdot \tilde{X}^\varepsilon= X_{h_{(\mu,a)}}\cdot \bar{\varepsilon}, $
for the $R_p$-reduced controlled underwater vehicle-rotor system
$(\mathcal{O}_{(\mu,a)}\times \mathbb{R}^2\times \mathbb{R}^{2*},
\omega^{-}_{\mathcal{O}_{(\mu,a)}\times \mathbb{R}^2\times \mathbb{R}^{2*}},
h_{(\mu,a)}, u_{(\mu,a)})$, if and only if they satisfy
the equation $T\bar{\varepsilon}\cdot(X_{h_{(\mu,a)} \cdot \bar{\varepsilon}})
= T\bar{\lambda}\cdot \tilde{X}\cdot\varepsilon. $\\

To sum up the above discussion, we have the following Theorem 4.2.
For convenience, the maps involved in
the following theorem are shown in Diagram-3.

\begin{center}
\hskip 0cm \xymatrix{ \mathbf{J}_Q^{-1}(\mu,a) \ar[r]^{i_{(\mu,a)}} & T^* Q
\ar[d]_{X_{H\cdot \varepsilon}} \ar[dr]^{\tilde{X}^\varepsilon} \ar[r]^{\pi_Q}
& Q \ar[d]^{\tilde{X}^\gamma} \ar[r]^{\gamma}
& T^*Q \ar[d]_{\tilde{X}} \ar[dr]_{X_{h_{(\mu,a)} \cdot\bar{\varepsilon}}} \ar[r]^{\pi_{(\mu,a)}}
& \;\;\; \mathcal{O}_{(\mu,a)}\times \mathbb{R}^2\times \mathbb{R}^{2*} \ar[d]^{X_{h_{(\mu,a)}}} \\
& T(T^*Q)  & TQ \ar[l]^{T\gamma}
& T(T^*Q) \ar[l]^{T\pi_Q} \ar[r]_{T\pi_{(\mu,a)}}
& \;\;\; T(\mathcal{O}_{(\mu,a)}\times \mathbb{R}^2\times \mathbb{R}^{2*})}
\end{center}
$$\mbox{Diagram-3}$$

\begin{theo}
In the case of coincident centers of buoyancy and gravity,
if the 5-tuple $(T^\ast Q,\textmd{SE}(3), \\ \omega_Q,H,u), $ where $Q=
\textmd{SE}(3)\times \mathbb{R}^2, $ is a regular point reducible
underwater vehicle-rotor system with the control torque $u$ acting on the rotors,
then for a point $(\mu,a) \in \mathfrak{se}^\ast(3)$, the regular
value of the momentum map $\mathbf{J}_Q:  T^* Q \cong\textmd{SE}(3)\times
\mathfrak{se}^\ast(3) \times \mathbb{R}^2 \times \mathbb{R}^{2*} \to
\mathfrak{se}^\ast(3)$, the $R_p$-reduced controlled underwater vehicle-rotor system is the 4-tuple
$(\mathcal{O}_{(\mu,a)} \times \mathbb{R}^2 \times
\mathbb{R}^{2*},\tilde{\omega}_{\mathcal{O}_{(\mu,a)} \times \mathbb{R}^2
\times \mathbb{R}^{2*}}^{-},h_{(\mu,a)},u_{(\mu,a)}). $
Assume that
$\gamma: Q \rightarrow T^* Q $ is an one-form on
$Q=\textmd{SE}(3)\times \mathbb{R}^2$,
and $\lambda=\gamma \cdot \pi_Q: T^* Q \rightarrow T^* Q, $
and $\varepsilon: T^* Q \rightarrow T^* Q $
is a $\textmd{SE}(3)_{(\mu,a)}$-invariant symplectic map,
where $(\textmd{SE}(3))_{(\mu,a)}$ is
the isotropy subgroup of co-adjoint $\textmd{SE}(3)$-action at the point $(\mu,a)$..
Denote $\tilde{X}^\gamma = T\pi_Q \cdot \tilde{X}\cdot \gamma$, and
$\tilde{X}^\varepsilon = T\pi_Q\cdot \tilde{X}\cdot \varepsilon$,
where $\tilde{X}=X_{(T^\ast Q,\textmd{SE}(3),\omega_Q,H,u)}$ is the dynamical vector field of
the controlled underwater vehicle-rotor system $(T^\ast Q,\textmd{SE}(3),\omega_Q,H,u)$.
Moreover, assume that $\textmd{Im}(\gamma)\subset \mathbf{J}_Q^{-1}(\mu,a), $ and it is
$\textmd{SE}(3)_{(\mu,a)}$-invariant,
and $\varepsilon(\mathbf{J}_Q^{-1}(\mu,a))\subset \mathbf{J}_Q^{-1}(\mu,a). $
Denote $\bar{\gamma}=\pi_{(\mu,a)}(\gamma):
Q \rightarrow \mathcal{O}_{(\mu,a)}\times \mathbb{R}^2\times \mathbb{R}^{2*}, $ and
$\bar{\lambda}=\pi_{(\mu,a)}(\lambda): T^* Q \rightarrow
\mathcal{O}_{(\mu,a)}\times \mathbb{R}^2\times \mathbb{R}^{2*}, $ and
$\bar{\varepsilon}=\pi_{(\mu,a)}(\varepsilon): \mathbf{J}_Q^{-1}(\mu,a)\rightarrow
\mathcal{O}_{(\mu,a)}\times \mathbb{R}^2\times \mathbb{R}^{2*}. $
Then the following two assertions hold:\\
\noindent $(\mathbf{i})$
If the one-form $\gamma: Q \rightarrow T^* Q $ is closed with respect to
$T\pi_Q: TT^* Q \rightarrow TQ, $
then $\bar{\gamma}$ is a solution of the Type I of Hamilton-Jacobi equation
$T\bar{\gamma}\cdot \tilde{X}^\gamma= X_{h_{(\mu,a)}}\cdot \bar{\gamma}; $\\
\noindent $(\mathbf{ii})$
The $\varepsilon$ and $\bar{\varepsilon} $ satisfy the Type II of Hamilton-Jacobi equation
$T\bar{\gamma}\cdot \tilde{X}^\varepsilon= X_{h_{(\mu,a)}}\cdot \bar{\varepsilon}, $
if and only if they satisfy
the equation $T\bar{\varepsilon}\cdot(X_{h_{(\mu,a)} \cdot \bar{\varepsilon}})
= T\bar{\lambda}\cdot \tilde{X}\cdot\varepsilon. $
\hskip 0.3cm $\blacksquare$
\end{theo}

It is worthy of noting that, for the regular point reducible controlled
underwater vehicle-rotor system $(T^\ast Q,\textmd{SE}(3),\omega_Q, H,u)$
with the $R_p$-reduced controlled underwater vehicle-rotor system
$(\mathcal{O}_{(\mu,a)} \times \mathbb{R}^2 \times
\mathbb{R}^{2*},\tilde{\omega}_{\mathcal{O}_{(\mu,a)} \times \mathbb{R}^2
\times \mathbb{R}^{2*}}^{-},
h_{(\mu,a)}, u_{(\mu,a)}) $, we know that the Hamiltonian vector fields
$X_{H}$ and $X_{h_{(\mu,a)}}$ for the corresponding
Hamiltonian system $(T^*Q,\textmd{SE}(3),\omega_Q, H)$
and its $R_p$-reduced system $(\mathcal{O}_{(\mu,a)} \times \mathbb{R}^2 \times
\mathbb{R}^{2*},\tilde{\omega}_{\mathcal{O}_{(\mu,a)} \times \mathbb{R}^2
\times \mathbb{R}^{2*}}^{-}, h_{(\mu,a)} )$, are $\pi_{(\mu,a)}$-related, that is,
$X_{h_{(\mu,a)}}\cdot \pi_{(\mu,a)}=T\pi_{(\mu,a)}\cdot X_{H}\cdot i_{(\mu,a)}.$
From the above Theorem 3.8 we can obtain the following Theorem 4.3,
which states the relationship between the solutions of Type II of
Hamilton-Jacobi equations and the regular point reduction.

\begin{theo}
In the case of coincident centers of buoyancy and gravity, for the regular point reducible controlled
underwater vehicle-rotor system $(T^\ast Q,\textmd{SE}(3),\omega_Q,H,u)$
with the $R_p$-reduced controlled underwater vehicle-rotor system
$(\mathcal{O}_{(\mu,a)} \times \mathbb{R}^2 \times
\mathbb{R}^{2*},\tilde{\omega}_{\mathcal{O}_{(\mu,a)} \times \mathbb{R}^2
\times \mathbb{R}^{2*}}^{-}, h_{(\mu,a)}, u_{(\mu,a)}) $,
assume that $\gamma: Q \rightarrow T^* Q $ is an one-form on
$Q=\textmd{SE}(3)\times \mathbb{R}^2 $, and $\varepsilon:
T^* Q \rightarrow T^* Q $ is a
$\textmd{SE}(3)_{(\mu,a)}$-invariant symplectic map,
$\bar{\varepsilon}=\pi_{(\mu,a)}(\varepsilon): \mathbf{J}_Q^{-1}(\mu, a )\rightarrow
\mathcal{O}_{(\mu,a)}\times \mathbb{R}^2\times \mathbb{R}^{2*}. $
Under the hypotheses and notations of Theorem 4.2, then we have that
$\varepsilon$ is a solution of the Type II of Hamilton-Jacobi equation
$T\gamma\cdot \tilde{X}^\varepsilon= X_H\cdot \varepsilon, $ for the
regular point reducible controlled underwater vehicle-rotor system
$(T^\ast Q,\textmd{SE}(3),\omega_Q,H,u), $ if and only if
$\varepsilon$ and $\bar{\varepsilon} $ satisfy the Type II of Hamilton-Jacobi equation
$T\bar{\gamma}\cdot \tilde{X}^\varepsilon= X_{h_{(\mu,a)}}\cdot \bar{\varepsilon}, $ for the
$R_p$-reduced controlled underwater vehicle-rotor system
$(\mathcal{O}_{(\mu,a)} \times \mathbb{R}^2 \times
\mathbb{R}^{2*},\tilde{\omega}_{\mathcal{O}_{(\mu,a)} \times \mathbb{R}^2
\times \mathbb{R}^{2*}}^{-}, h_{(\mu,a)}, u_{(\mu,a)}) $.
 \hskip 0.3cm $\blacksquare$
\end{theo}

\begin{rema}
When the underwater vehicle does not carry any internal rotor, in
this case the configuration space is $Q=W=\textmd{SE}(3), $ the
motion of underwater vehicle is just the rotation and translation
motion of a rigid body, the
above $R_p$-reduced controlled underwater vehicle-rotors system is
just the Marsden-Weinstein reduced underwater vehicle system, that is, 3-tuple
$(\mathcal{O}_{(\mu,a)},
\omega_{\mathcal{O}_{(\mu,a)}}^{-}, h_{(\mu,a)} )$,
where $\mathcal{O}_{(\mu,a)}\subset \mathfrak{se}^\ast(3)$
is the co-adjoint orbit of the semidirect product
Lie group $\textmd{SE}(3)=\textmd{SO}(3)\circledS \mathbb{R}^3$,
$\omega_{\mathcal{O}_{(\mu,a)}}^{-}$
is the orbit symplectic form on $\mathcal{O}_{(\mu,a)}$,
which is induced by the (-)-semidirect product Lie-Poisson brackets on
$\mathfrak{so}^\ast(3)\circledS \mathbb{R}^3$, that is,
the heavy top Lie-Poisson bracket on $\mathfrak{se}^\ast(3)$,
see Marsden et al. \cite{mamiorpera07, marawe84a, marawe84b},
$h_{(\mu,a)}(\Pi,P)\cdot \pi_{(\mu,a)}=H(A,c,\Pi,P)|_{\mathcal{O}_{(\mu,a)}}$.
From the above Theorem 4.2 we can obtain the two
types of Hamilton-Jacobi equation for the Marsden-Weinstein
reduced underwater vehicle system $(\mathcal{O}_{(\mu,a)},
\omega_{\mathcal{O}_{(\mu,a)}}^{-}, h_{(\mu,a)} )$,
see Wang \cite{wa17}. Moreover, from the above Theorem 4.3 we can also
state the relationship between the solutions of Type II of Hamilton-Jacobi equations
and the Marsden-Weinstein reduction.
\end{rema}

\section{Underwater Vehicle-Rotor System with Non-coincident Centers
of Buoyancy and Gravity}

In this section, we consider the controlled underwater vehicle-rotor system
with non-coincident centers of buoyancy and gravity as an RCH
system on the generalization of semidirect product Lie group
$Q=W\times V$, where $W=\textmd{SE}(3)\circledS \mathbb{R}^3
\cong (\textmd{SO}(3)\circledS \mathbb{R}^3) \circledS \mathbb{R}^3$
and $V= \mathbb{R}^2$. In this case, we give the regular point reduction
and the two types of Hamilton-Jacobi equation for the controlled underwater
vehicle-rotor system with non-coincident centers of buoyancy and gravity.

\subsection{Symmetric Reduction }

In the following we first give the regular point reduction of the controlled underwater
vehicle-rotors system with non-coincident centers of buoyancy and
gravity, and give precisely the $R_p$-reduced symplectic form
and the $R_p$-reduced controlled underwater
vehicle-rotors system. \\

Because the drift in the direction of gravity
breaks the symmetry and the underwater vehicle-rotors system is no
longer $\textmd{SE}(3)$ invariant. In this case, its physical phase
space is $T^\ast \textmd{SE}(3)\times T^* V$ and the symmetry group
is $S^1\circledS \mathbb{R}^3\cong \textmd{SE}(2)\times \mathbb{R}^3$,
regarded as rotations about the third principal axis, that is, the
axis of gravity, and translation. By the semidirect product
reduction theorem, see Marsden et al. \cite{mamiorpera07} or Leonard
and Marsden \cite{lema97}, we know that the reduction of $T^\ast
\textmd{SE}(3)$ by $S^1\circledS \mathbb{R}^3$ gives a space which
is symplectically diffeomorphic to the reduced space obtained by the
reduction of $T^* W=T^\ast (\textmd{SE}(3)\circledS \mathbb{R}^3)$
by left action of $W=\textmd{SE}(3)\circledS \mathbb{R}^3$, that is
the co-adjoint orbit $\mathcal{O}_{(\mu,a_1,a_2)} \subset
\mathfrak{se}^\ast(3)\circledS \mathbb{R}^{3*}\cong T^\ast W/W $,
where $(\mu,a_1,a_2) \in \mathfrak{se}^\ast(3)\circledS \mathbb{R}^{3*}$. In fact,
in this case, we can identify the phase space $T^\ast
\textmd{SE}(3)$ with the reduction of the cotangent bundle of the
double semidirect product Lie group $W= \textmd{SE}(3)\circledS
\mathbb{R}^3 =(\textmd{SO}(3)\circledS \mathbb{R}^3)\circledS
\mathbb{R}^3$ by the Euclidean translation subgroup $\mathbb{R}^3$
and identifies the symmetry group $S^1\circledS \mathbb{R}^3$ with
isotropy group $(W_{a_2})_{(\mu,a_1)}\circledS \mathbb{R}^3$, where
$W_{a_2}=\{ A\in \textmd{SO}(3)\mid Aa_2=a_2 \}=S^1$, which is
Abelian and $(W_{a_2})_{(\mu,a_1)}= W_{a_2} =S^1,\; \forall
(\mu,a_1) \in \mathfrak{se}^\ast(3)_{a_2}$, and $a_2$ is a vector
aligned with the direction of gravity and where
$\textmd{SO}(3)$ acts on $\mathbb{R}^3$ in the standard way.
See the Proposition 3.1and Proposition 3.2 on the reduction by stages for
semidirect product Lie group. \\

Assume that semidirect product Lie group $W=\textmd{SE}(3)\circledS
\mathbb{R}^3 $ acts freely and properly on $Q=
(\textmd{SE}(3)\circledS \mathbb{R}^3 )\times \mathbb{R}^2 $ by the left
translation on $\textmd{SE}(3)\circledS \mathbb{R}^3 $, then the
action of $W$ on the phase space $T^\ast Q$ is by cotangent lift of
left translation on $Q$ at the identity, that is, $\Phi^{T*}: W \times
T^\ast (\textmd{SE}(3)\circledS \mathbb{R}^3)\times T^\ast V \cong
\textmd{SE}(3)\circledS \mathbb{R}^3 \times \textmd{SE}(3)\circledS
\mathbb{R}^3 \times \mathfrak{se}^\ast(3)\circledS \mathbb{R}^{3*}
\times \mathbb{R}^2 \times \mathbb{R}^{2*} \to T^* Q \cong \textmd{SE}(3)\circledS
\mathbb{R}^3 \times \mathfrak{se}^\ast(3)\circledS
\mathbb{R}^{3*}\times \mathbb{R}^2 \times \mathbb{R}^{2*},$ given by
$$\Phi^{T*}((B,s_1,s_2)((A,c_1,c_2),(\Pi,w_1,w_2),\theta,l))=((BA,s_1+Bc_1,s_2+Bc_2),(\Pi,w_1,w_2),\theta,l), $$
for any $A,B\in \textmd{SO}(3), \; \Pi \in \mathfrak{so}^\ast(3), \;
s_i, c_i \in \mathbb{R}^3, \; (A,c_1, c_2), (B,s_1,s_2) \in \textmd{SE}(3)\circledS \mathbb{R}^3,
\; w_i \in \mathbb{R}^{3*},  \; i=1,2,
\; (\Pi, w_1, w_2) \in \mathfrak{se}^\ast(3)\circledS \mathbb{R}^{3*},
\; \theta \in \mathbb{R}^2, \; l \in \mathbb{R}^{2*} $, and
$\textmd{SO}(3)$ acts on $\mathbb{R}^3$ in the standard way.
Assume that the action is free, proper and symplectic, and admits an
associated $\operatorname{Ad}^\ast$-equivariant momentum map
$\mathbf{J}_Q: T^\ast Q \cong \textmd{SE}(3)\circledS
\mathbb{R}^3\times \mathfrak{se}^\ast(3)\circledS \mathbb{R}^{3*}
\times \mathbb{R}^2 \times \mathbb{R}^{2*} \to
\mathfrak{se}^\ast(3)\circledS \mathbb{R}^{3*}$ for the cotangent lift
$\textmd{SE}(3)\circledS \mathbb{R}^3$-action. If $(\Pi,w_1,w_2) \in
\mathfrak{se}^\ast(3)\circledS \mathbb{R}^{3*}$ is a regular value of
$\mathbf{J}_Q$, then the $R_p$-reduced space $(T^\ast
Q)_{(\Pi,w_1,w_2)}=
\mathbf{J}^{-1}_Q(\Pi,w_1,w_2)/(\textmd{SE}(3)\circledS \mathbb{R}^3
)_{(\Pi,w_1,w_2)}$ is symplectically diffeomorphic to the
orbit space $\mathcal{O}_{(\Pi,w_1,w_2)} \times \mathbb{R}^2 \times
\mathbb{R}^{2*} \subset \mathfrak{se}^\ast(3)\circledS \mathbb{R}^{3*}
\times \mathbb{R}^2 \times \mathbb{R}^{2*} $, where
$(\textmd{SE}(3)\circledS \mathbb{R}^3)_{(\Pi,w_1,w_2)}$ is the
isotropy subgroup of co-adjoint $\textmd{SE}(3)\circledS
\mathbb{R}^3$-action at the point $(\Pi,w_1,w_2)\in
\mathfrak{se}^\ast(3)\circledS \mathbb{R}^{3*}$.\\

From Marsden et al. \cite{mamiorpera07} we know that $\mathfrak{w}^\ast=\mathfrak{se}^\ast(3)\circledS
\mathbb{R}^{3*}$ is a Poisson manifold with respect to its semidirect
product Lie-Poisson bracket defined by
\begin{align}   \{F,K\}_{\mathfrak{w}^\ast}(\Pi,\Gamma, P)= & -\Pi\cdot (\nabla_\Pi
F\times \nabla_\Pi K) -\Gamma\cdot(\nabla_\Pi F\times \nabla_\Gamma
K-\nabla_\Pi K\times \nabla_\Gamma F)\nonumber \\ &
- P \cdot(\nabla_\Pi F\times \nabla_P K-\nabla_\Pi K\times
\nabla_P F), \;\; \label{5.1}
\end{align}
where $F,K: \mathfrak{se}^\ast(3)\circledS \mathbb{R}^{3*}\to
\mathbb{R}, $ and $ (\Pi, \Gamma, P) \in
\mathfrak{w}^\ast=\mathfrak{se}^\ast(3)\circledS \mathbb{R}^{3*}$. For
a point $ (\Pi_0, \Gamma_0, P_0) =(\mu,a_1,a_2) \in \mathfrak{se}^\ast(3)\circledS \mathbb{R}^{3*}$, the
co-adjoint orbit $\mathcal{O}_{(\mu,a_1,a_2)} \subset
\mathfrak{se}^\ast(3)\circledS \mathbb{R}^{3*}$ has the induced orbit
symplectic form $\omega^{-}_{\mathcal{O}_{(\mu,a_1,a_2)}}$, which
coincides with the restriction of the semidirect product Lie-Poisson
bracket on $\mathfrak{se}^\ast(3)\circledS \mathbb{R}^{3*}$ to the
co-adjoint orbit $\mathcal{O}_{(\mu,a_1,a_2)}$, and the co-adjoint
orbits $(\mathcal{O}_{(\mu,a_1,a_2)},
\omega_{\mathcal{O}_{(\mu,a_1,a_2)}}^{-}), \\ (\mu,a_1,a_2)\in
\mathfrak{se}^\ast(3)\circledS \mathbb{R}^{3*},$ form the symplectic
leaves of the Poisson manifold $(\mathfrak{se}^\ast(3)\circledS
\mathbb{R}^{3*},\{\cdot,\cdot\}_{\mathfrak{w}^\ast}). $ Let
$\omega_{\mathbb{R}^2}$ be the canonical symplectic form on $T^\ast
\mathbb{R}^2 \cong \mathbb{R}^{2} \times \mathbb{R}^{2*}$, and it induces
a canonical Poisson bracket $\{\cdot,\cdot\}_{\mathbb{R}^2}$ on
$T^\ast \mathbb{R}^2$ given by $(3.7)$. Thus, we can induce a
symplectic form $\tilde{\omega}^{-}_{\mathcal{O}_{(\mu,a_1,a_2)} \times \mathbb{R}^2
\times \mathbb{R}^{2*}}=
\pi_{\mathcal{O}_{(\mu,a_1,a_2)}}^\ast
\omega^{-}_{\mathcal{O}_{(\mu,a_1,a_2)}}+ \pi_{\mathbb{R}^2}^\ast
\omega_{\mathbb{R}^2}$ on the smooth manifold
$\mathcal{O}_{(\mu,a_1,a_2)} \times \mathbb{R}^2
\times \mathbb{R}^{2*}$, where the maps
$\pi_{\mathcal{O}_{(\mu,a_1,a_2)}}: \mathcal{O}_{(\mu,a_1,a_2)} \times \mathbb{R}^2
\times \mathbb{R}^{2*} \to
\mathcal{O}_{(\mu,a_1,a_2)}$ and $\pi_{\mathbb{R}^2}:
\mathcal{O}_{(\mu,a_1,a_2)} \times \mathbb{R}^2
\times \mathbb{R}^{2*} \to \mathbb{R}^2 \times \mathbb{R}^{2*}$ are
canonical projections, and induce a Poisson bracket
$\{\cdot,\cdot\}_{-}= \pi_{\mathfrak{w}^\ast}^\ast
\{\cdot,\cdot\}_{\mathfrak{w}^\ast}+ \pi_{\mathbb{R}^2}^\ast
\{\cdot,\cdot\}_{\mathbb{R}^2}$ on the smooth manifold
$\mathfrak{se}^\ast(3)\circledS \mathbb{R}^{3*}\times \mathbb{R}^2
\times \mathbb{R}^{2*}$, where the maps $\pi_{\mathfrak{w}^\ast}:
\mathfrak{se}^\ast(3)\circledS \mathbb{R}^{3*} \times \mathbb{R}^2
\times \mathbb{R}^{2*} \to \mathfrak{w}^*=\mathfrak{se}^\ast(3)\circledS \mathbb{R}^{3*}$
and $\pi_{\mathbb{R}^2}: \mathfrak{se}^\ast(3)\circledS \mathbb{R}^{3*}
\times \mathbb{R}^2 \times \mathbb{R}^{2*} \to \mathbb{R}^2 \times
\mathbb{R}^{2*}$ are canonical projections, and such that
$(\mathcal{O}_{(\mu,a_1,a_2)} \times \mathbb{R}^2
\times \mathbb{R}^{2*}, \tilde{\omega}_{\mathcal{O}_{(\mu,a_1,a_2)} \times \mathbb{R}^2
\times \mathbb{R}^{2*}} ^{-})$ is
a symplectic leaf of the Poisson manifold
$(\mathfrak{se}^\ast(3)\circledS
\mathbb{R}^{3*} \times \mathbb{R}^2 \times \mathbb{R}^{2*}, \{\cdot,\cdot\}_{-}). $\\

From the above expression $(2.2)$ of the Hamiltonian, we know that
$H(A,c,b,\Pi,\Gamma, P, \theta,l)$ is invariant under
is invariant under the cotangent lift of the left
$\textmd{SE}(3)\circledS \mathbb{R}^3$-action
$\Phi^{T^*}:\textmd{SE}(3)\circledS \mathbb{R}^3\times T^\ast Q\to
T^\ast Q$. Moreover, from the  semidirect product Poisson
bracket on $\mathfrak{w}^*=\mathfrak{se}^\ast(3)\circledS \mathbb{R}^{3*}$ and the Poisson bracket on
$T^\ast \mathbb{R}^2$, we can get the Poisson bracket
on $\mathfrak{se}^\ast(3)\circledS \mathbb{R}^{3*}\times \mathbb{R}^2\times \mathbb{R}^{2*} $,
that is, for $F,K: \mathfrak{se}^\ast(3)\circledS \mathbb{R}^{3*}\times \mathbb{R}^2\times
\mathbb{R}^{2*} \to \mathbb{R}, $ we have that
\begin{align} \{F,K\}_{-}(\Pi,\Gamma,P,\theta,l) &
= -\Pi\cdot(\nabla_\Pi F\times\nabla_\Pi K)-\Gamma\cdot(\nabla_\Pi
F\times \nabla_\Gamma K-\nabla_\Pi K\times \nabla_\Gamma F) \nonumber \\
& \;\;\; -P \cdot(\nabla_\Pi
F\times \nabla_P K-\nabla_\Pi K\times \nabla_P F)+
\{F,K\}_{\mathbb{R}^2}(\theta,l).
\label{5.2}
\end{align}
Hence, the Hamiltonian vector fields of underwater vehicle-rotor system are given by
\begin{align*}
X_{H}(\Pi) = \{\Pi,\; H \}_{-}& = -\Pi\cdot(\nabla_\Pi\Pi\times\nabla_\Pi
H) -\Gamma\cdot(\nabla_\Pi\Pi\times\nabla_\Gamma
H-\nabla_\Pi H \times\nabla_\Gamma\Pi) \\
& \;\;\;\;\;\; -P\cdot(\nabla_\Pi\Pi\times\nabla_P
H-\nabla_\Pi H \times\nabla_P \Pi)+ \{\Pi,\; H \}_{\mathbb{R}^2}\\
& =(\Pi_1,\Pi_2,\Pi_3)\times (\frac{\Pi_1-l_1}{ \bar{I}_1}, \frac{\Pi_2-l_2}{
\bar{I}_2}, \frac{\Pi_3}{\bar{I}_3})+mgh(\Gamma_1,\Gamma_2,\Gamma_3)\times (\chi_1,\chi_2,\chi_3)\\
& \;\;\;\;\;\; + (P_1,P_2,P_3)\times (\frac{P_1}{m_1},\frac{P_2}{m_2},\frac{P_3}{m_3})
+ \sum_{k=1}^2(\frac{\partial \Pi}{\partial \theta_k}
\frac{\partial H}{\partial l_k}- \frac{\partial
H}{\partial \theta_k}\frac{\partial \Pi}{\partial l_k})\\
& = ( \frac{(\bar{I}_2-\bar{I}_3)\Pi_2\Pi_3+
\bar{I}_3\Pi_3l_2}{\bar{I}_2\bar{I}_3}+ \frac{(m_2-m_3)P_2P_3}{m_2m_3}+
gh(\Gamma_2\chi_3-\Gamma_3\chi_2), \\
& \;\;\;\;\;\;
\frac{(\bar{I}_3-\bar{I}_1)\Pi_3\Pi_1-
\bar{I}_3\Pi_3l_1}{\bar{I}_3\bar{I}_1}+  \frac{(m_3-m_1)P_3P_1}{m_3m_1}
+gh(\Gamma_3\chi_1-\Gamma_1\chi_3), \\
& \;\;\;\;\;\;
\frac{(\bar{I}_1-\bar{I}_2)\Pi_1\Pi_2-
\bar{I}_1\Pi_1l_2 + \bar{I}_2\Pi_2l_1}{\bar{I}_1\bar{I}_2}+ \frac{(m_1-m_2)P_1P_2}{m_1m_2}
 + gh(\Gamma_1\chi_2-\Gamma_2\chi_1) ),
\end{align*}
since $\nabla_{\Pi_i}\Pi_i=1,\; \nabla_{\Pi_i}\Pi_j=0, \; i\neq j, \;  \nabla_{\Pi_i}\Gamma_j=\nabla_{\Gamma_i}\Pi_j=0, \; \nabla_{\Pi_i}P_j=\nabla_{P_i}\Pi_j=0, $
and $\chi=(\chi_1,\chi_2,\chi_3), \; \nabla_{\Gamma_j} H= mgh\chi_j,
\; \nabla_{P_i} H= P_i /{m_i}, \; \nabla_{\Pi_k} H= (\Pi_k-l_k)/\bar{I}_k ,\; \nabla_{\Pi_3}
H= \Pi_3/\bar{I}_3 , \; \frac{\partial
\Pi}{\partial \theta_k}= \frac{\partial H}{\partial
\theta_k}=0, \; i, j=1,2,3, \; k=1,2. $

\begin{align*}
 X_{H}(\Gamma)
= \{\Gamma,\; H\}_{-} &= -\Pi\cdot(\nabla_\Pi\Gamma\times\nabla_\Pi
H)-\Gamma\cdot(\nabla_\Pi\Gamma\times\nabla_\Gamma
H-\nabla_\Pi H \times\nabla_\Gamma\Gamma)\\
& \;\;\;\;\;\;
-P\cdot(\nabla_\Pi\Gamma\times\nabla_P
H-\nabla_\Pi H \times\nabla_P\Gamma)
+  \{\Gamma,\; H\}_{\mathbb{R}^2} \\
& =\nabla_\Gamma\Gamma\cdot (\Gamma\times\nabla_\Pi
H)+ \sum_{k=1}^2(\frac{\partial \Gamma}{\partial \theta_k}
\frac{\partial H}{\partial l_k}- \frac{\partial
H}{\partial \theta_k}\frac{\partial \Gamma}{\partial l_k})\\
& = (\Gamma_1,\Gamma_2,\Gamma_3)\times (\frac{(\Pi_1- l_1)}{ \bar{I}_1},\;\;
\frac{(\Pi_2- l_2)}{ \bar{I}_2}, \;\; \frac{\Pi_3}{\bar{I}_3}) \\
&= ( \frac{\bar{I}_2\Gamma_2\Pi_3-\bar{I}_3\Gamma_3\Pi_2+
\bar{I}_3\Pi_3l_2}{\bar{I}_2\bar{I}_3}, \;\;
\frac{\bar{I}_3\Gamma_3\Pi_1-\bar{I}_1\Gamma_1\Pi_3-
\bar{I}_3\Pi_3l_1}{\bar{I}_3\bar{I}_1}, \\
& \;\;\;\;\;\;
\frac{\bar{I}_1\Gamma_1\Pi_2-\bar{I}_2\Gamma_2\Pi_1-
\bar{I}_1\Pi_1l_2 + \bar{I}_2\Pi_2l_1}{\bar{I}_1\bar{I}_2} ),
\end{align*}
since $\nabla_{\Gamma_i}\Gamma_i =1, \; \nabla_{\Gamma_i}\Gamma_j=0, \; i\neq j, \;
\nabla_{\Pi_i}\Gamma_j =\nabla_{P_i}\Gamma_j =0, $ and
$\nabla_{\Pi_k} H= (\Pi_k-l_k)/\bar{I}_k ,\; \nabla_{\Pi_3}
H= \Pi_3/\bar{I}_3 , \; \frac{\partial
\Gamma_j}{\partial \theta_k}= \frac{\partial H}{\partial
\theta_k}=0, \; i, j= 1,2,3, \; k=1,2.$

\begin{align*}
 X_{H}(P) = \{P,\; H\}_{-} &=-\Pi\cdot(\nabla_\Pi P\times\nabla_\Pi
H)- \Gamma\cdot(\nabla_\Pi P\times\nabla_\Gamma
H-\nabla_\Pi H \times\nabla_\Gamma P) \\
& \;\;\;\;\;\; -P\cdot(\nabla_\Pi P\times\nabla_P
H-\nabla_\Pi H \times\nabla_P P)+  \{P,\; H\}_{\mathbb{R}^2} \\
& =\nabla_P P\cdot (P\times\nabla_\Pi
H)+ \sum_{k=1}^2(\frac{\partial P}{\partial \theta_k}
\frac{\partial H}{\partial l_k}- \frac{\partial
H}{\partial \theta_k}\frac{\partial P}{\partial l_k})\\
& = (P_1,P_2,P_3)\times (\frac{(\Pi_1- l_1)}{ \bar{I}_1},\;\;
\frac{(\Pi_2- l_2)}{ \bar{I}_2}, \;\; \frac{\Pi_3}{\bar{I}_3}) \\
&= ( \frac{\bar{I}_2P_2\Pi_3-\bar{I}_3P_3\Pi_2+
\bar{I}_3\Pi_3l_2}{\bar{I}_2\bar{I}_3}, \;\;
\frac{\bar{I}_3P_3\Pi_1-\bar{I}_1P_1\Pi_3-
\bar{I}_3\Pi_3l_1}{\bar{I}_3\bar{I}_1}, \\
& \;\;\;\;\;\;
\frac{\bar{I}_1P_1\Pi_2-\bar{I}_2P_2\Pi_1-
\bar{I}_1\Pi_1l_2 + \bar{I}_2\Pi_2l_1}{\bar{I}_1\bar{I}_2} ),
\end{align*}
since $\nabla_{P_i}P_i =1, \; \nabla_{P_i}P_j=0, \; i\neq j, \; \nabla_{\Pi_i}P_j =\nabla_{\Gamma_i}P_j=0, $ and
$\nabla_{\Pi_k} H= (\Pi_k-l_k)/\bar{I}_k ,\; \nabla_{\Pi_3}
H= \Pi_3/\bar{I}_3 , \; \frac{\partial
P_j}{\partial \theta_k}= \frac{\partial H}{\partial
\theta_k}=0, \; i, j= 1,2,3, \; k=1,2.$

\begin{align*}
 X_{H}(\theta) = \{\theta,\; H \}_{-}
 & = -\Pi\cdot(\nabla_\Pi\theta\times\nabla_\Pi
H) -\Gamma\cdot(\nabla_\Pi\theta\times\nabla_\Gamma
H-\nabla_\Pi H \times\nabla_\Gamma\theta) \\
& \;\;\;\;\;\;
-P\cdot(\nabla_\Pi\theta\times\nabla_P
H-\nabla_\Pi H \times\nabla_P \theta)
+ \{\theta,\; H \}_{\mathbb{R}^2}\\
& =\sum_{k=1}^2(\frac{\partial \theta}{\partial \theta_k}
\frac{\partial H}{\partial l_k}- \frac{\partial
H}{\partial \theta_k}\frac{\partial \theta}{\partial l_k})= ( -\frac{(\Pi_1- l_1)}{ \bar{I}_1}
+\frac{l_1}{J_1},\;\; -\frac{(\Pi_2-
l_2)}{ \bar{I}_2} +\frac{l_2}{J_2} ),
\end{align*}
since $\nabla_{\Pi_i}\theta=\nabla_{\Gamma_i}\theta=\nabla_{P_i}\theta=0, $ $\frac{\partial \theta_k}{\partial
\theta_k}= 1, \; \frac{\partial \theta_n}{\partial \theta_k}=
0, \; n\neq k, \;\; \frac{\partial H}{\partial \theta_k}=0, $
and $\frac{\partial H}{\partial l_k}= -(\Pi_k-l_k)/\bar{I}_k
+\frac{l_k}{J_k}, \; i= 1,2,3, \; k,n=1,2. $

\begin{align*}
 X_{H}(l) = \{l,\; H \}_{-}
 & = -\Pi\cdot(\nabla_\Pi l \times\nabla_\Pi
H) -\Gamma\cdot(\nabla_\Pi l \times\nabla_\Gamma
H-\nabla_\Pi H \times\nabla_\Gamma l) \\
& \;\;\;\;\;\;
-P\cdot(\nabla_\Pi l \times\nabla_P
H-\nabla_\Pi H \times\nabla_P l)
+ \{l,\; H \}_{\mathbb{R}^2}\\
& = \sum_{k=1}^2(\frac{\partial l}{\partial \theta_k}
\frac{\partial H}{\partial l_k}- \frac{\partial
H}{\partial \theta_k}\frac{\partial l}{\partial l_k})=(0,0),
\end{align*}
since $\nabla_{\Pi_i} l=\nabla_{\Gamma_i} l=\nabla_{P_i} l=0,$ and $\frac{\partial l}{\partial \theta_k}=
\frac{\partial H}{\partial \theta_k}=0, \; i=1,2,3, \; k=1,2. $\\

Moreover, if we consider the underwater vehicle-rotor system with a control torque
$u: T^\ast Q \to \mathcal{C}$ acting on the rotors, where the control subset
$\mathcal{C}\subset T^* Q $ is a fiber submanifold,
and assume that $u\in \mathcal{C}$ is invariant under the cotangent lift
of the left $\textmd{SE}(3)\circledS \mathbb{R}^3$-action, and
the dynamical vector field of the regular point reducible
controlled underwater vehicle-rotor system $(T^\ast Q,\textmd{SE}(3)\circledS \mathbb{R}^3,\omega_Q,H,u)$
can be expressed by
\begin{align}
\tilde{X}= X_{(T^\ast Q,\textmd{SE}(3)\circledS \mathbb{R}^3,\omega_Q,H,u)}= X_H+ \textnormal{vlift}(u),
\label{5.3}
\end{align}
where $\textnormal{vlift}(u)= \textnormal{vlift}(u)\cdot X_H $
is the change of $X_H$ under the action of the control torque $u$.\\

Since the Hamiltonian
$H(A,c,b,\Pi,\Gamma, P,\theta,l)$ is invariant under the cotangent lift $\Phi^{T^*}$ of the left
$\textmd{SE}(3)\circledS \mathbb{R}^3$-action, for the point
$(\Pi_0,\Gamma_0, P_0)=(\mu,a_1,a_2)\in \mathfrak{se}^\ast(3)\circledS \mathbb{R}^{3*}$,
the regular value of $\mathbf{J}_Q$, we
have the $R_p$-reduced Hamiltonian
$h_{(\mu,a_1,a_2)}(\Pi,\Gamma,P,\theta,l):\mathcal{O}_{(\mu,a_1,a_2)}\times\mathbb{R}^2
\times \mathbb{R}^{2*} (\subset \mathfrak{se}^\ast (3)\circledS \mathbb{R}^{3*}\times
\mathbb{R}^2\times \mathbb{R}^{2*}) \to \mathbb{R}$ given by
$$h_{(\mu,a_1,a_2)}(\Pi,\Gamma,P,\theta,l)\cdot \pi_{(\mu,a_1,a_2)}
=H(A,c,b,\Pi,\Gamma,P,\theta,l)|_{\mathcal{O}_{(\mu,a_1,a_2)}\times
\mathbb{R}^2\times \mathbb{R}^{2*}}, $$ where $\pi_{(\mu,a_1,a_2)}:
\mathbf{J}_Q^{-1}(\mu,a_1,a_2) \rightarrow \mathcal{O}_{(\mu,a_1,a_2)}
\times\mathbb{R}^2\times\mathbb{R}^{2*}.$
Moreover, for the $R_p$-reduced
Hamiltonian $h_{(\mu,a_1,a_2)}(\Pi,\Gamma, P,\theta, l): \mathcal{O}_{(\mu,a_1,a_2)}\times
\mathbb{R}^2\times \mathbb{R}^{2*} \to \mathbb{R}$, we have the
Hamiltonian vector field
$$X_{h_{(\mu,a_1,a_2)}}(K_{(\mu,a_1,a_2)})
=\{K_{(\mu,a_1,a_2)},h_{(\mu,a_1,a_2)}\}_{-}|_{\mathcal{O}_{(\mu,a_1,a_2)}
\times \mathbb{R}^2\times \mathbb{R}^{2*}}, $$ where
$K_{(\mu,a_1,a_2)}(\Pi,\Gamma, P,\theta, l): \mathcal{O}_{(\mu,a_1,a_2)}\times
\mathbb{R}^2\times \mathbb{R}^{2*} \to \mathbb{R}.$
Assume that $u\in \mathcal{C} \cap
\mathbf{J}^{-1}_Q(\mu,a_1,a_2)$ and the  $R_p$-reduced control torque $u_{(\mu,a_1,a_2)}:
\mathcal{O}_{(\mu,a_1,a_2)} \times\mathbb{R}^2\times\mathbb{R}^{2*} \to \mathcal{C}_{(\mu,a_1,a_2)}
(\subset \mathcal{O}_{(\mu,a_1,a_2)} \times\mathbb{R}^2\times\mathbb{R}^{2*}) $ is
given by $u_{(\mu,a_1,a_2)}(\Pi,\Gamma,P,\theta,l)\cdot \pi_{(\mu,a_1,a_2)}=
u(A,c,b,\Pi,\Gamma, P,\theta,l )|_{\mathcal{O}_{(\mu,a_1,a_2)}
\times\mathbb{R}^2\times\mathbb{R}^{2*} }, $ where $\mathcal{C}_{(\mu,a_1,a_2)}
= \pi_{(\mu,a_1,a_2)}(\mathcal{C}\cap \mathbf{J}_Q^{-1}(\mu,a_1,a_2)). $
The $R_p$-reduced controlled underwater vehicle-rotor
system is the 4-tuple $(\mathcal{O}_{(\mu,a_1,a_2)} \times \mathbb{R}^2 \times
\mathbb{R}^{2*},\tilde{\omega}_{\mathcal{O}_{(\mu,a_1,a_2)} \times \mathbb{R}^2
\times \mathbb{R}^{2*}}^{-},h_{(\mu,a_1,a_2)},u_{(\mu,a_1,a_2)}), $ where
$\tilde{\omega}_{\mathcal{O}_{(\mu,a_1,a_2)} \times \mathbb{R}^2 \times
\mathbb{R}^{2*}}^{-}$ is the induced symplectic form
from the above Poisson bracker on  $\mathcal{O}_{(\mu,a_1,a_2)}
\times \mathbb{R}^2\times \mathbb{R}^{2*} ,$ such that
the Hamiltonian vector field
\begin{align*}
X_{h_{(\mu,a_1,a_2)}}(K_{(\mu,a_1,a_2)}) &
=\tilde{\omega}_{\mathcal{O}_{(\mu,a_1,a_2)} \times \mathbb{R}^2 \times
\mathbb{R}^{2*}}^{-}(X_{K_{(\mu,a_1,a_2)}}, X_{h_{(\mu,a_1,a_2)}})\\
& =\{K_{(\mu,a_1,a_2)},h_{(\mu,a_1,a_2)}\}_{-}|_{\mathcal{O}_{(\mu,a_1,a_2)}
\times\mathbb{R}^2\times\mathbb{R}^{2*} }.
\end{align*}
Moreover, assume that the dynamical vector field of the $R_p$-reduced controlled
underwater vehicle-rotor system $(\mathcal{O}_{(\mu,a_1,a_2)} \times \mathbb{R}^2 \times
\mathbb{R}^{2*},\tilde{\omega}_{\mathcal{O}_{(\mu,a_1,a_2)} \times
\mathbb{R}^2 \times
\mathbb{R}^{2*}}^{-},h_{(\mu,a_1,a_2)},u_{(\mu,a_1,a_2)})$ can be expressed
by
\begin{align} X_{(\mathcal{O}_{(\mu,a_1,a_2)} \times \mathbb{R}^2 \times
\mathbb{R}^{2*},\tilde{\omega}_{\mathcal{O}_{(\mu,a_1,a_2)} \times \mathbb{R}^2 \times
\mathbb{R}^{2*}}^{-},h_{(\mu,a_1,a_2)},u_{(\mu,a_1,a_2)})} = X_{h_{(\mu,a_1,a_2)}} +
\mbox{vlift}(u_{(\mu,a_1,a_2)}),
\label{5.4}
\end{align}
where $\mbox{vlift}(u_{(\mu,a_1,a_2)})=
\mbox{vlift}(u_{(\mu,a_1,a_2)})X_{h_{(\mu,a_1,a_2)}} \in T(\mathcal{O}_{(\mu,a_1,a_2)}
\times \mathbb{R}^2 \times
\mathbb{R}^{2*}), $ is the change of $X_{h_{(\mu,a_1,a_2)}}$
under the action of the $R_p$-reduced control torque $u_{(\mu,a_1,a_2)}$.
The dynamical vector fields of the controlled
underwater vehicle-rotor system and the $R_p$-reduced controlled
underwater vehicle-rotor system satisfy the condition
\begin{align}
& X_{(\mathcal{O}_{(\mu,a_1,a_2)} \times \mathbb{R}^2 \times
\mathbb{R}^{2*},\tilde{\omega}_{\mathcal{O}_{(\mu,a_1,a_2)} \times \mathbb{R}^2 \times
\mathbb{R}^{2*}}^{-}, h_{(\mu,a_1,a_2)}, u_{(\mu,a_1,a_2)})}\cdot \pi_{(\mu,a_1,a_2)} \nonumber \\
& =T\pi_{(\mu,a_1,a_2)}\cdot X_{(T^\ast Q,\textmd{SE}(3)\circledS
\mathbb{R}^3,\omega_Q,H,u)}\cdot i_{(\mu,a_1,a_2)}.
\label{5.5}\end{align}
See Marsden et al \cite{mawazh10} and Wang \cite{wa18}.\\

To sum up the above discussion, we have the following Theorem 5.1.
\begin{theo}
In the case of non-coincident centers of buoyancy and gravity, the
underwater vehicle-rotor system with the control torque $u$ acting
on the rotors, that is, the 5-tuple \\ $(T^\ast
Q,\textmd{SE}(3)\circledS \mathbb{R}^3,\omega_Q,H,u ), $ where $Q=
\textmd{SE}(3)\circledS \mathbb{R}^3 \times \mathbb{R}^2, $ is a
regular point reducible RCH system. For a point
$(\mu,a_1,a_2) \in \mathfrak{se}^\ast(3)\circledS \mathbb{R}^{3*}$, the
regular value of the momentum map $\mathbf{J}_Q:
T^* Q \cong \textmd{SE}(3)\circledS \mathbb{R}^3\times
\mathfrak{se}^\ast(3)\circledS \mathbb{R}^{3*} \times \mathbb{R}^2
\times \mathbb{R}^{2*} \to \mathfrak{se}^\ast(3)\circledS
\mathbb{R}^{3*}$, the $R_p$-reduced
controlled underwater vehicle-rotor system is the 4-tuple
$(\mathcal{O}_{(\mu,a_1,a_2)} \times
\mathbb{R}^2 \times \mathbb{R}^{2*},\tilde{\omega}_{\mathcal{O}_{(\mu,a_1,a_2)} \times
\mathbb{R}^2 \times \mathbb{R}^{2*}}^{-},h_{(\mu,a_1,a_2)},u_{(\mu,a_1,a_2)}),
$ where $\mathcal{O}_{(\mu,a_1,a_2)}\\ \subset
\mathfrak{se}^\ast(3)\circledS \mathbb{R}^{3*}$ is the co-adjoint orbit,
$\tilde{\omega}_{\mathcal{O}_{(\mu,a_1,a_2)} \times
\mathbb{R}^2 \times \mathbb{R}^{2*}}^{-}$ is the induced symplectic form
on $\mathcal{O}_{(\mu,a_1,a_2)} \times
\mathbb{R}^2 \times \mathbb{R}^{2*}$,
$h_{(\mu,a_1,a_2)}(\Pi, \Gamma, P,\theta,l)\cdot \pi_{(\mu,a_1,a_2)}
=H(A,c,b,\Pi,\Gamma, P,\theta,l)|_{\mathcal{O}_{(\mu,a_1,a_2)} \times
\mathbb{R}^2 \times \mathbb{R}^{2*}},$ \\
$u_{(\mu,a_1,a_2)}(\Pi,\Gamma,P,\theta,l)\cdot \pi_{(\mu,a_1,a_2)}= u(A,c,b,\Pi,
\Gamma, P, \theta,l)|_{\mathcal{O}_{(\mu,a_1,a_2)} \times
\mathbb{R}^2 \times \mathbb{R}^{2*}}$, and
the dynamical vector field of the $R_p$-reduced controlled
underwater vehicle-rotor system satisfies
$(5.4)$ and $(5.5)$.
\hskip 0.3cm $\blacksquare$
\end{theo}

\subsection{Hamilton-Jacobi Theorems}

In the following we shall give precisely the geometric constraint conditions
of the $R_p$-reduced symplectic form
$\tilde{\omega}_{\mathcal{O}_{(\mu,a_1,a_2)} \times \mathbb{R}^2 \times
\mathbb{R}^{2*}}^{-}$ for the dynamical vector field of the regular point
reducible controlled underwater vehicle-rotor system with
non-coincident centers of buoyancy and gravity, that is, Type I and Type II of
Hamilton-Jacobi equation for the $R_p$-reduced controlled underwater vehicle-rotor system
$(\mathcal{O}_{(\mu,a_1,a_2)}\times \mathbb{R}^2\times \mathbb{R}^{2*},
\omega^{-}_{\mathcal{O}_{(\mu,a_1,a_2)}\times \mathbb{R}^2\times \mathbb{R}^{2*}},
h_{(\mu,a_1,a_2)}, u_{(\mu,a_1,a_2)}) .$\\

Assume that $\gamma: Q \rightarrow T^* Q $ is an one-form on
$Q=W\times \mathbb{R}^2=\textmd{SE}(3)\circledS \mathbb{R}^3 \times \mathbb{R}^2 $, and
for $(A,c) \in \textmd{SE}(3)$, $b \in \mathbb{R}^3 $ and
$\theta \in \mathbb{R}^2$, denote $\gamma(A,c,b,\theta)=(\gamma_1, \cdots, \gamma_{22})$,
and $\gamma$ is closed with respect to $T\pi_Q: TT^* Q \rightarrow TQ. $
For $(\mu,a_1,a_2) \in \mathfrak{se}^\ast(3)\circledS \mathbb{R}^3,$ the regular value of
$\mathbf{J}_Q: T^* Q \cong \textmd{SE}(3)\circledS \mathbb{R}^3\times
\mathfrak{se}^\ast(3)\circledS \mathbb{R}^{3*} \times \mathbb{R}^2
\times \mathbb{R}^{2*} \to \mathfrak{se}^\ast(3)\circledS
\mathbb{R}^{3*}$, assume that
$\textmd{Im}(\gamma)\subset \mathbf{J}_Q^{-1}(\mu,a_1,a_2), $ and it is
$(\textmd{SE}(3)\circledS \mathbb{R}^3)_{(\mu,a_1,a_2)}$-invariant,
and $\bar{\gamma}=\pi_{(\mu,a_1,a_2)}(\gamma):
Q=\textmd{SE}(3)\circledS \mathbb{R}^3\times \mathbb{R}^2 \rightarrow \mathcal{O}_{(\mu,a_1,a_2)} \times
\mathbb{R}^2 \times \mathbb{R}^{2*}$. Denote by
$\bar{\gamma}= (\bar{\gamma}_1, \cdots, \bar{\gamma}_{13}) \in \mathcal{O}_{(\mu,a_1,a_2)} \times
\mathbb{R}^2 \times \mathbb{R}^{2*}(\subset \mathfrak{se}^\ast(3)\circledS \mathbb{R}^3
\times \mathbb{R}^2 \times \mathbb{R}^{2*}), $ where
$\pi_{(\mu,a_1,a_2)}: \mathbf{J}_Q^{-1}(\mu,a_1,a_2) \rightarrow \mathcal{O}_{(\mu,a_1,a_2)} \times
\mathbb{R}^2 \times \mathbb{R}^{2*}. $ We choose that
$(\Pi,\Gamma, P,\theta,l)\in
\mathcal{O}_{(\mu,a_1,a_2)}\times \mathbb{R}^2 \times \mathbb{R}^{2*}, $ and
$\Pi=(\Pi_1,\Pi_2,\Pi_3)=(\bar{\gamma}_1,\bar{\gamma}_2,\bar{\gamma}_3)$,
$\Gamma=(\Gamma_1,\Gamma_2,\Gamma_3)=(\bar{\gamma}_4,\bar{\gamma}_5,\bar{\gamma}_6)$,
$P=(P_1,P_2,P_3)=(\bar{\gamma}_7,\bar{\gamma}_8,\bar{\gamma}_9)$,
$\theta=(\bar{\gamma}_{10}, \bar{\gamma}_{11})$ and $l=(\bar{\gamma}_{12},
\bar{\gamma}_{13}) $. Then $h_{(\mu,a_1,a_2)} \cdot \bar{\gamma}:
Q=\textmd{SE}(3)\circledS \mathbb{R}^3 \times \mathbb{R}^2 \rightarrow \mathbb{R} $ is given by
\begin{align*} & h_{(\mu,a_1,a_2)}(\Pi,\Gamma, P,\theta,l) \cdot \bar{\gamma}=
H (A,c,b,\Pi,\Gamma, P,\theta,l)|_{\mathcal{O}_{(\mu,a_1,a_2)} \times \mathbb{R}^2 \times \mathbb{R}^{2*}}
\cdot \bar{\gamma}\\ &
=\frac{1}{2}[\frac{(\bar{\gamma}_1-
\bar{\gamma}_{12})^2}{\bar{I}_1}+\frac{(\bar{\gamma}_2-
\bar{\gamma}_{13})^2}{\bar{I}_2}
+\frac{\bar{\gamma}_3^2}{\bar{I}_3}+
\frac{\bar{\gamma}_7^2}{m_1}+\frac{\bar{\gamma}_8^2}{m_2}+\frac{\bar{\gamma}_9^2}{m_3}
+\frac{\bar{\gamma}_{12}^2}{J_1}+\frac{\bar{\gamma}_{13}^2}{J_2}] \\
& \;\;\; \;\;\; + gh(\bar{\gamma}_4\cdot\chi_1+ \bar{\gamma}_5\cdot\chi_2+
\bar{\gamma}_6\cdot\chi_3),
\end{align*} and the vector field
\begin{align*}
& X_{h_{(\mu,a_1,a_2)}}(\Pi) \cdot \bar{\gamma}=\{\Pi,\;
h_{(\mu,a_1,a_2)}\}_{-}|_{\mathcal{O}_{(\mu,a_1,a_2)}
\times \mathbb{R}^2\times \mathbb{R}^{2*}}\cdot \bar{\gamma}\\
& = -\Pi\cdot(\nabla_\Pi\Pi\times\nabla_\Pi (h_{(\mu,a_1,a_2)})) \cdot
\bar{\gamma}- \Gamma\cdot(\nabla_\Pi\Pi\times\nabla_\Gamma
(h_{(\mu,a_1,a_2)})-\nabla_\Pi
(h_{(\mu,a_1,a_2)}) \times\nabla_\Gamma \Pi)\cdot \bar{\gamma}\\
& \;\;\;\; -P\cdot(\nabla_\Pi\Pi\times\nabla_P
(h_{(\mu,a_1,a_2)})-\nabla_\Pi
(h_{(\mu,a_1,a_2)}) \times\nabla_P \Pi)\cdot \bar{\gamma}+ \{\Pi,\;
h_{(\mu,a_1,a_2)}\}_{\mathbb{R}^2}|_{\mathcal{O}_{(\mu,a_1,a_2)} \times
\mathbb{R}^2\times \mathbb{R}^{2*}}\cdot \bar{\gamma}\\
& = -\nabla_\Pi\Pi \cdot(\nabla_\Pi (h_{(\mu,a_1,a_2)}) \times \Pi)\cdot
\bar{\gamma}- \nabla_\Pi\Pi \cdot (\nabla_\Gamma
(h_{(\mu,a_1,a_2)}) \times \Gamma) \cdot \bar{\gamma}\\
& \;\;\;\; - \nabla_\Pi\Pi \cdot (\nabla_P
(h_{(\mu,a_1,a_2)}) \times P) \cdot \bar{\gamma}
+ \sum_{k=1}^2(\frac{\partial \Pi}{\partial \theta_k}
\frac{\partial (h_{(\mu,a_1,a_2)})}{\partial l_k}- \frac{\partial
(h_{(\mu,a_1,a_2)})}{\partial
\theta_k}\frac{\partial \Pi}{\partial l_k})\cdot \bar{\gamma}\\
&=(\Pi_1,\Pi_2,\Pi_3)\times (\frac{(\Pi_1- l_1)}{ \bar{I}_1},
\frac{(\Pi_2- l_2)}{ \bar{I}_2}, \frac{\Pi_3}{\bar{I}_3})\cdot
\bar{\gamma}
+gh(\Gamma_1,\Gamma_2,\Gamma_3)\times (\chi_1,\chi_2,\chi_3)\cdot \bar{\gamma}\\
& \;\;\;\;\;\; +(P_1,P_2,P_3)\times (\frac{P_1}{m_1},\frac{P_2}{m_2},\frac{P_3}{m_3})\cdot \bar{\gamma}\\
&= ( \frac{(\bar{I}_2-\bar{I}_3)\bar{\gamma}_2\bar{\gamma}_3+
\bar{I}_3\bar{\gamma}_3\bar{\gamma}_{13}}{\bar{I}_2\bar{I}_3}
+  \frac{(m_2-m_3)\bar{\gamma}_8\bar{\gamma}_9}{m_2m_3}
+ gh(\bar{\gamma}_5\chi_3-\bar{\gamma}_6\chi_2),\\
& \;\;\;\;\;\; \frac{(\bar{I}_3-\bar{I}_1)\bar{\gamma}_3
\bar{\gamma}_1-
\bar{I}_3\bar{\gamma}_3\bar{\gamma}_{12}}{\bar{I}_3\bar{I}_1}
+  \frac{(m_3-m_1)\bar{\gamma}_9\bar{\gamma}_7}{m_3m_1}
+ gh(\bar{\gamma}_6\chi_1-\bar{\gamma}_4\chi_3),\\
& \;\;\;\;\;\;
\frac{(\bar{I}_1-\bar{I}_2)\bar{\gamma}_1\bar{\gamma}_2-\bar{I}_1\bar{\gamma}_1
\bar{\gamma}_{13}+
\bar{I}_2\bar{\gamma}_2\bar{\gamma}_{12}}{\bar{I}_1\bar{I}_2}
+  \frac{(m_1-m_2)\bar{\gamma}_7\bar{\gamma}_8}{m_1m_2}
+ gh(\bar{\gamma}_4\chi_2-\bar{\gamma}_5\chi_1) ),
\end{align*}
since $\nabla_{\Pi_i}\Pi_i=1,\; \nabla_{\Pi_i}\Pi_j=0, \; i\neq j, \;
\nabla_{\Pi_i}\Gamma_j=\nabla_{\Gamma_i}\Pi_j=\nabla_{\Pi_i}P_j=\nabla_{P_i}\Pi_j=0$
and $\chi=(\chi_1,\chi_2,\chi_3), \; \nabla_{\Gamma_j}
(h_{(\mu,a_1,a_2)})= mgh\chi_j, \; \nabla_{\Pi_k} (h_{(\mu,a_1,a_2)}
)= (\Pi_k-l_k)/\bar{I}_k ,\; \nabla_{\Pi_3}
(h_{(\mu,a_1,a_2)})= \Pi_3/\bar{I}_3 ,$
$ \nabla_{P_i} (h_{(\mu,a)})= \frac{P_i}{m_i}, \;
\frac{\partial \Pi}{\partial \theta_k}= \frac{\partial (h_{(\mu,a_1,a_2)})}{\partial
\theta_k}=0, \; i, j=1,2,3, \; k=1,2. $

\begin{align*}
& X_{h_{(\mu,a_1,a_2)}}(\Gamma) \cdot \bar{\gamma}=\{\Gamma,\;
h_{(\mu,a_1,a_2)}\}_{-}|_{\mathcal{O}_{(\mu,a_1,a_2)} \times \mathbb{R}^2\times
\mathbb{R}^{2*}}\cdot \bar{\gamma}\\
& =-\Pi\cdot(\nabla_\Pi\Gamma\times\nabla_\Pi (h_{(\mu,a_1,a_2)})) \cdot
\bar{\gamma}-\Gamma\cdot(\nabla_\Pi\Gamma\times\nabla_\Gamma
(h_{(\mu,a_1,a_2)})-\nabla_\Pi
(h_{(\mu,a_1,a_2)}) \times\nabla_\Gamma\Gamma)\cdot \bar{\gamma}\\
& \;\;\;\; -P\cdot(\nabla_\Pi \Gamma \times\nabla_P
(h_{(\mu,a_1,a_2)})-\nabla_\Pi
(h_{(\mu,a_1,a_2)}) \times\nabla_P \Gamma )\cdot \bar{\gamma}+ \{\Gamma,\;
h_{(\mu,a_1,a_2)}\}_{\mathbb{R}^2}|_{\mathcal{O}_{(\mu,a_1,a_2)} \times
\mathbb{R}^2\times \mathbb{R}^{2*}}\cdot \bar{\gamma}\\
& =\nabla_\Gamma\Gamma\cdot(\Gamma\times\nabla_\Pi
(h_{(\mu,a_1,a_2)}))\cdot \bar{\gamma}+ \sum_{k=1}^2(\frac{\partial
\Gamma}{\partial \theta_k} \frac{\partial (h_{(\mu,a_1,a_2)})}{\partial
l_k}- \frac{\partial (h_{(\mu,a_1,a_2)})}{\partial
\theta_k}\frac{\partial \Gamma}{\partial l_k})\cdot \bar{\gamma}\\
&=(\Gamma_1,\Gamma_2,\Gamma_3)\times (\frac{(\Pi_1- l_1)}{
\bar{I}_1}, \frac{(\Pi_2- l_2)}{ \bar{I}_2},
\frac{\Pi_3}{\bar{I}_3})\cdot \bar{\gamma}\\
&= ( \frac{\bar{I}_2\bar{\gamma}_5\bar{\gamma}_3-\bar{I}_3\bar{\gamma}_6\bar{\gamma}_2+
\bar{I}_3\bar{\gamma}_3\bar{\gamma}_{13}}{\bar{I}_2\bar{I}_3}, \;\;
\frac{\bar{I}_3\bar{\gamma}_6\bar{\gamma}_1-\bar{I}_1\bar{\gamma}_4\bar{\gamma}_3-
\bar{I}_3\bar{\gamma}_3\bar{\gamma}_{12}}{\bar{I}_3\bar{I}_1}, \\
& \;\;\;\;\;\;
\frac{\bar{I}_1\bar{\gamma}_4\bar{\gamma}_2-\bar{I}_2\bar{\gamma}_5\bar{\gamma}_1-
\bar{I}_1\bar{\gamma}_1\bar{\gamma}_{13} + \bar{I}_2\bar{\gamma}_2\bar{\gamma}_{12}}{\bar{I}_1\bar{I}_2} ),
\end{align*}
since $\nabla_{\Gamma_i}\Gamma_i =1, \; \nabla_{\Gamma_i}\Gamma_j=0, \; i\neq j, \;
\nabla_{\Pi_i}\Gamma_j =\nabla_{P_i}\Gamma_j =0, $ and
$\nabla_{\Pi_k} (h_{(\mu,a_1,a_2)})= (\Pi_k-l_k)/\bar{I}_k ,$
$ \nabla_{\Pi_3} (h_{(\mu,a_1,a_2)})= \Pi_3/\bar{I}_3 , \; \frac{\partial
\Gamma_j}{\partial \theta_k}= \frac{\partial (h_{(\mu,a_1,a_2)})}{\partial
\theta_k}=0, \; i, j= 1,2,3, \; k=1,2.$

\begin{align*}
& X_{h_{(\mu,a_1,a_2)}}(P) \cdot \bar{\gamma}=\{P,\;
h_{(\mu,a_1,a_2)}\}_{-}|_{\mathcal{O}_{(\mu,a_1,a_2)} \times \mathbb{R}^2\times
\mathbb{R}^{2*}}\cdot \bar{\gamma}\\
& =-\Pi\cdot(\nabla_\Pi P \times\nabla_\Pi (h_{(\mu,a_1,a_2)})) \cdot
\bar{\gamma} -\Gamma \cdot(\nabla_\Pi P \times\nabla_\Gamma
(h_{(\mu,a_1,a_2)})-\nabla_\Pi
(h_{(\mu,a_1,a_2)}) \times\nabla_\Gamma P)\cdot \bar{\gamma}\\
& \;\;\;\; -P\cdot(\nabla_\Pi P\times\nabla_P
(h_{(\mu,a_1,a_2)})-\nabla_\Pi
(h_{(\mu,a_1,a_2)}) \times\nabla_P P)\cdot \bar{\gamma}
+ \{P,\; h_{(\mu,a_1,a_2)}\}_{\mathbb{R}^2}|_{\mathcal{O}_{(\mu,a_1,a_2)} \times
\mathbb{R}^2\times \mathbb{R}^{2*}}\cdot \bar{\gamma}\\
& =\nabla_P P \cdot(P \times\nabla_\Pi
(h_{(\mu,a_1,a_2)}))\cdot \bar{\gamma}+ \sum_{k=1}^2(\frac{\partial
P}{\partial \theta_k} \frac{\partial (h_{(\mu,a_1,a_2)})}{\partial
l_k}- \frac{\partial (h_{(\mu,a_1,a_2)})}{\partial
\theta_k}\frac{\partial P}{\partial l_k})\cdot \bar{\gamma}\\
&=(P_1,P_2,P_3)\times (\frac{(\Pi_1- l_1)}{
\bar{I}_1}, \frac{(\Pi_2- l_2)}{ \bar{I}_2},
\frac{\Pi_3}{\bar{I}_3})\cdot \bar{\gamma}\\
&= ( \frac{\bar{I}_2\bar{\gamma}_8\bar{\gamma}_3-\bar{I}_3\bar{\gamma}_9\bar{\gamma}_2+
\bar{I}_3\bar{\gamma}_3\bar{\gamma}_{13}}{\bar{I}_2\bar{I}_3}, \;\;
\frac{\bar{I}_3\bar{\gamma}_9\bar{\gamma}_1-\bar{I}_1\bar{\gamma}_7\bar{\gamma}_3-
\bar{I}_3\bar{\gamma}_3\bar{\gamma}_{12}}{\bar{I}_3\bar{I}_1}, \\
& \;\;\;\;\;\;
\frac{\bar{I}_1\bar{\gamma}_7\bar{\gamma}_2-\bar{I}_2\bar{\gamma}_8\bar{\gamma}_1-
\bar{I}_1\bar{\gamma}_1\bar{\gamma}_{13} + \bar{I}_2\bar{\gamma}_2\bar{\gamma}_{12}}{\bar{I}_1\bar{I}_2} ),
\end{align*}
since $\nabla_{P_i}P_i =1, \; \nabla_{P_i}P_j=0, \; i\neq j, \; \nabla_{\Pi_i}P_j =\nabla_{\Gamma_i}P_j =0, $ and
$\nabla_{\Pi_k} (h_{(\mu,a_1,a_2)})= (\Pi_k-l_k)/\bar{I}_k ,$
$ \nabla_{\Pi_3} (h_{(\mu,a_1,a_2)})= \Pi_3/\bar{I}_3 , \; \frac{\partial
P_j}{\partial \theta_k}= \frac{\partial (h_{(\mu,a_1,a_2)})}{\partial
\theta_k}=0, \; i, j= 1,2,3, \; k=1,2.$

\begin{align*}
& X_{h_{(\mu,a_1,a_2)}}(\theta) \cdot \bar{\gamma}=\{\theta,\;
h_{(\mu,a_1,a_2)}\}_{-}|_{\mathcal{O}_{(\mu,a_1,a_2)} \times \mathbb{R}^2\times
\mathbb{R}^{2*}}\cdot \bar{\gamma}\\
& =-\Pi\cdot(\nabla_\Pi \theta \times\nabla_\Pi (h_{(\mu,a_1,a_2)}))\cdot
\bar{\gamma}-\Gamma \cdot(\nabla_\Pi \theta \times\nabla_\Gamma
(h_{(\mu,a_1,a_2)})-\nabla_\Pi
(h_{(\mu,a_1,a_2)}) \times\nabla_\Gamma \theta)\cdot \bar{\gamma}\\
& \;\;\;\; -P \cdot(\nabla_\Pi \theta \times\nabla_P
(h_{(\mu,a_1,a_2)})-\nabla_\Pi
(h_{(\mu,a_1,a_2)}) \times\nabla_P \theta)\cdot \bar{\gamma}+ \{\theta,\;
h_{(\mu,a_1,a_2)}\}_{\mathbb{R}^2}|_{\mathcal{O}_{(\mu,a_1,a_2)} \times
\mathbb{R}^2\times \mathbb{R}^{2*}}\cdot \bar{\gamma}\\
& = \sum_{k=1}^2(\frac{\partial \theta}{\partial \theta_k}
\frac{\partial (h_{(\mu,a_1,a_2)})}{\partial l_k}- \frac{\partial
(h_{(\mu,a_1,a_2)})}{\partial
\theta_k}\frac{\partial \theta}{\partial l_k})\cdot \bar{\gamma}
 = (-\frac{(\bar{\gamma}_1- \bar{\gamma}_{12})}{ \bar{I}_1}+
\frac{\bar{\gamma}_{12}}{J_1}, \; -\frac{(\bar{\gamma}_2-
\bar{\gamma}_{13})}{ \bar{I}_2}+ \frac{\bar{\gamma}_{13}}{J_2}),
\end{align*}
since $\nabla_{\Pi_i}\theta=\nabla_{\Gamma_i}\theta=\nabla_{P_i}\theta=0, $ $\frac{\partial \theta_k}{\partial
\theta_k}= 1, \; \frac{\partial \theta_n}{\partial \theta_k}=
0, \; n\neq k, \;\; \frac{\partial (h_{(\mu,a_1,a_2)})}{\partial \theta_k}=0, $
and $\frac{\partial (h_{(\mu,a_1,a_2)})}{\partial l_k}= -(\Pi_k-l_k)/\bar{I}_k
+\frac{l_k}{J_k}, \; i= 1,2,3, \; k,n=1,2. $

\begin{align*}
& X_{h_{(\mu,a_1,a_2)}}(l) \cdot \bar{\gamma}=\{l,\;
h_{(\mu,a_1,a_2)}\}_{-}|_{\mathcal{O}_{(\mu,a_1,a_2)} \times \mathbb{R}^2\times
\mathbb{R}^{2*}}\cdot \bar{\gamma}\\
& =-\Pi\cdot(\nabla_\Pi l \times\nabla_\Pi (h_{(\mu,a_1,a_2)})) \cdot
\bar{\gamma}-\Gamma \cdot(\nabla_\Pi l \times\nabla_\Gamma
(h_{(\mu,a_1,a_2)})-\nabla_\Pi
(h_{(\mu,a_1,a_2)}) \times\nabla_\Gamma l)\cdot \bar{\gamma}\\
& \;\;\;\; -P \cdot(\nabla_\Pi l \times\nabla_P
(h_{(\mu,a_1,a_2)})-\nabla_\Pi
(h_{(\mu,a_1,a_2)}) \times\nabla_P l)\cdot \bar{\gamma}+ \{l,\;
h_{(\mu,a_1,a_2)}\}_{\mathbb{R}^2}|_{\mathcal{O}_{(\mu,a_1,a_2)} \times
\mathbb{R}^2\times \mathbb{R}^{2*}}\cdot \bar{\gamma}\\
& = \sum_{k=1}^2(\frac{\partial l}{\partial \theta_k} \frac{\partial
(h_{(\mu,a_1,a_2)})}{\partial l_k}- \frac{\partial (h_{(\mu,a_1,a_2)})}{\partial
\theta_k}\frac{\partial l}{\partial l_k})\cdot \bar{\gamma}=(0,0),
\end{align*}
since $\nabla_{\Pi_i} l =\nabla_{\Gamma_i} l=\nabla_{P_i} l=0,$ and $\frac{\partial l}{\partial \theta_k}=
\frac{\partial (h_{(\mu,a_1,a_2)})}{\partial \theta_k}=0, \; i=1,2,3, \; k=1,2. $\\

On the other hand, from the expressions of the dynamical vector field $\tilde{X}$
and Hamiltonian vector field $X_H$, we have that
\begin{align*}
\tilde{X}(\Pi, \Gamma, P,\theta, l)^\gamma &
=T\pi_Q \cdot \tilde{X}\cdot\gamma(\Pi, \Gamma, P,\theta, l)\\
& =T\pi_Q \cdot (X_H+ \textnormal{vlift}(u))\cdot\gamma (\Pi, \Gamma, P,\theta, l)\\
&=T\pi_Q \cdot X_H \cdot\gamma (\Pi, \Gamma, P,\theta, l) = X_H\cdot\gamma(\Pi,  \Gamma, P,\theta, l),
\end{align*}
that is,
\begin{align*}
\tilde{X}(\Pi)^\gamma & = X_H(\Pi)\cdot\gamma \\
& = ( \frac{(\bar{I}_2-\bar{I}_3)\gamma_{11}\gamma_{12}+
\bar{I}_3\gamma_{12}\gamma_{22}}{\bar{I}_2\bar{I}_3}
+ \frac{(m_2-m_3)\gamma_{17}\gamma_{18}}{m_2m_3}+
gh(\gamma_{14}\chi_3-\gamma_{15}\chi_2),\\
& \;\;\;\;\;\;
\frac{(\bar{I}_3-\bar{I}_1)\gamma_{12}\gamma_{10}-
\bar{I}_3\gamma_{12}\gamma_{21}}{\bar{I}_3\bar{I}_1}
+ \frac{(m_3-m_1)\gamma_{18}\gamma_{16}}{m_3m_1}+gh(\gamma_{15}\chi_1-\gamma_{13}\chi_3), \\
& \;\;\;\;\;\;
\frac{(\bar{I}_1-\bar{I}_2)\gamma_{10}\gamma_{11}-
\bar{I}_1\gamma_{10}\gamma_{22} + \bar{I}_2\gamma_{11}\gamma_{21}}{\bar{I}_1\bar{I}_2}
 + \frac{(m_1-m_2)\gamma_{16}\gamma_{17}}{m_1m_2}+ gh(\gamma_{13}\chi_2-\gamma_{14}\chi_1) ),
\end{align*}

\begin{align*}
\tilde{X}(\Gamma)^\gamma & = X_H(\Gamma)\cdot\gamma \\
&= ( \frac{\bar{I}_2\gamma_{14}\gamma_{12}-\bar{I}_3\gamma_{15}\gamma_{11}+
\bar{I}_3\gamma_{12}\gamma_{22}}{\bar{I}_2\bar{I}_3}, \;\;
\frac{\bar{I}_3\gamma_{15}\gamma_{10}-\bar{I}_1\gamma_{13}\gamma_{12}-
\bar{I}_3\gamma_{12}\gamma_{21}}{\bar{I}_3\bar{I}_1}, \\
& \;\;\;\;\;\;
\frac{\bar{I}_1\gamma_{13}\gamma_{11}-\bar{I}_2\gamma_{14}\gamma_{10}-
\bar{I}_1\gamma_{10}\gamma_{22} + \bar{I}_2\gamma_{11}\gamma_{21}}{\bar{I}_1\bar{I}_2} ),
\end{align*}

\begin{align*}
\tilde{X}(P)^\gamma & = X_H(P)\cdot\gamma \\
&= ( \frac{\bar{I}_2\gamma_{17}\gamma_{12}-\bar{I}_3\gamma_{18}\gamma_{11}+
\bar{I}_3\gamma_{12}\gamma_{22}}{\bar{I}_2\bar{I}_3}, \;\;
\frac{\bar{I}_3\gamma_{18}\gamma_{10}-\bar{I}_1\gamma_{16}\gamma_{12}-
\bar{I}_3\gamma_{12}\gamma_{21}}{\bar{I}_3\bar{I}_1}, \\
& \;\;\;\;\;\;
\frac{\bar{I}_1\gamma_{16}\gamma_{11}-\bar{I}_2\gamma_{17}\gamma_{10}-
\bar{I}_1\gamma_{10}\gamma_{22} + \bar{I}_2\gamma_{11}\gamma_{21}}{\bar{I}_1\bar{I}_2} ),
\end{align*}

\begin{align*}
\tilde{X}(\theta)^\gamma & = X_H(\theta)\cdot\gamma
 = ( -\frac{(\gamma_{10}- \gamma_{21})}{ \bar{I}_1}
+\frac{\gamma_{21}}{J_1},\;\; -\frac{(\gamma_{11}- \gamma_{22})}{ \bar{I}_2}
+\frac{\gamma_{22}}{J_2} ),
\end{align*}
\begin{align*}
\tilde{X}(l)^\gamma & = X_H(l)\cdot\gamma=(0,0).
\end{align*}
Since $\gamma$ is closed with respect to
$T\pi_{(\textmd{SE}(3) \times \mathbb{R}^2)}: TT^* (\textmd{SE}(3) \times \mathbb{R}^2)
\rightarrow T(\textmd{SE}(3) \times \mathbb{R}^2), $
then $\pi_{(\textmd{SE}(3) \times \mathbb{R}^2)}^*(\mathbf{d}\gamma)=0.$ We choose that
$(\gamma_{10},\gamma_{11},\gamma_{12})=\Pi=(\Pi_1,\Pi_2,\Pi_3)=
(\bar{\gamma}_1,\bar{\gamma}_2,\bar{\gamma}_3), \;
(\gamma_{13},\gamma_{14},\gamma_{15})=\Gamma=(\Gamma_1,\Gamma_2,\Gamma_3)=
(\bar{\gamma}_4,\bar{\gamma}_5,\bar{\gamma}_6), $
$ (\gamma_{16},\gamma_{17},\gamma_{18})=P=(P_1,P_2,P_3)=
(\bar{\gamma}_7,\bar{\gamma}_8,\bar{\gamma}_9), $
and $(\gamma_{19},\gamma_{20})= \theta= (\theta_1,\theta_2)
=(\bar{\gamma}_{10},\bar{\gamma}_{11}),
\; (\gamma_{21},\gamma_{22})= l=(l_1,l_2)
=(\bar{\gamma}_{12},\bar{\gamma}_{13}). $
Hence
\begin{align*}
& T\bar{\gamma}\cdot \tilde{X}(\Pi)^\gamma
= X_{h_{(\mu,a_1,a_2)}}(\Pi) \cdot \bar{\gamma}, \;\;\;
T\bar{\gamma}\cdot \tilde{X}(\Gamma)^\gamma
= X_{h_{(\mu,a_1,a_2)}}(\Gamma) \cdot \bar{\gamma}, \;\;\;
 T\bar{\gamma}\cdot \tilde{X}(P)^\gamma
= X_{h_{(\mu,a_1,a_2)}}(P) \cdot \bar{\gamma}, \\
& T\bar{\gamma}\cdot \tilde{X}(\theta)^\gamma
= X_{h_{(\mu,a_1,a_2)}}(\theta) \cdot \bar{\gamma}, \;\;\;\;\;
T\bar{\gamma}\cdot \tilde{X}(l)^\gamma
= X_{h_{(\mu,a_1,a_2)}}(l) \cdot \bar{\gamma}.
\end{align*}
Thus, the Type I of Hamilton-Jacobi equation for the
$R_p$-reduced controlled underwater vehicle-rotor system
$(\mathcal{O}_{((\mu,a_1,a_2)}\times \mathbb{R}^2\times \mathbb{R}^{2*},
\omega^{-}_{\mathcal{O}_{(\mu,a_1,a_2)}\times \mathbb{R}^2\times \mathbb{R}^{2*}},
h_{(\mu,a_1,a_2)}, u_{(\mu,a_1,a_2)})$ holds.\\

Next, for $(\mu,a_1,a_2) \in \mathfrak{se}^\ast(3)\circledS \mathbb{R}^3,$
the regular value of $\mathbf{J}_Q$,
and a $(\textmd{SE}(3)\circledS \mathbb{R}^3)_{(\mu,a_1,a_2)}$-invariant symplectic map
$\varepsilon: T^* Q \rightarrow T^* Q,$
assume that $\varepsilon(A,c,b,\Pi,\Gamma, P,\theta, l) =(\varepsilon_1,\cdots,
\varepsilon_{22}),$ and
$\varepsilon(\mathbf{J}^{-1}(\mu,a_1,a_2)\subset
\mathbf{J}^{-1}(\mu,a_1,a_2). $ Denote by
$\bar{\varepsilon}=\pi_{(\mu,a_1,a_2)}(\varepsilon):
\mathbf{J}^{-1}(\mu,a_1,a_2) \rightarrow
\mathcal{O}_{(\mu,a_1,a_2)}\times \mathbb{R}^2\times \mathbb{R}^{2*}, $
and $\bar{\varepsilon}=(\bar{\varepsilon}_1,\cdots,\bar{\varepsilon}_{13})
\in \mathcal{O}_{(\mu,a_1,a_2)}\times \mathbb{R}^2\times \mathbb{R}^{2*}
(\subset \mathfrak{se}^\ast(3)\times \mathbb{R}^2\times \mathbb{R}^{2*}), $
and $\lambda= \gamma \cdot \pi_Q: T^* Q \rightarrow T^* Q,$
and $\lambda(A,c,b,\Pi,\Gamma, P, \theta, l) =(\lambda_1,\cdots, \lambda_{22}),$ and
$\bar{\lambda}=\pi_{(\mu,a_1,a_2)}(\lambda): T^* Q
\rightarrow \mathcal{O}_{(\mu,a_1,a_2)}\times \mathbb{R}^2\times \mathbb{R}^{2*}, $ and
$\bar{\lambda}=
(\bar{\lambda}_1,\cdots,\bar{\lambda}_{13}) \in
\mathcal{O}_{(\mu,a_1,a_2)}\times \mathbb{R}^2\times \mathbb{R}^{2*}
(\subset \mathfrak{se}^\ast(3)\circledS \mathbb{R}^3\times \mathbb{R}^2\times \mathbb{R}^{2*}). $
We choose that $(\Pi,\Gamma,P,\theta,l)\in
\mathcal{O}_{(\mu,a_1,a_2)}\times \mathbb{R}^2 \times \mathbb{R}^{2*}, $ and
$\Pi=(\Pi_1,\Pi_2,\Pi_3)
=(\bar{\varepsilon}_1,\bar{\varepsilon}_2,\bar{\varepsilon}_3),
\; \Gamma= (\Gamma_1,\Gamma_2,\Gamma_3)=
(\bar{\varepsilon}_4,\bar{\varepsilon}_5,\bar{\varepsilon}_6), $
$P= (P_1,P_2,P_3)=
(\bar{\varepsilon}_7,\bar{\varepsilon}_8,\bar{\varepsilon}_9), $
$\theta= (\theta_1,\theta_2)
=(\bar{\varepsilon}_{10}, \bar{\varepsilon}_{11})$ and
$ l=(l_1,l_2)=(\bar{\varepsilon}_{12}, \bar{\varepsilon}_{13}) $.
Then $ h_{(\mu,a_1,a_2)} \cdot \bar{\varepsilon}: T^* Q \rightarrow \mathbb{R} $ is given by
\begin{align*}
& h_{(\mu,a_1,a_2)}(\Pi,\Gamma, P,\theta,l) \cdot \bar{\varepsilon}=
H(A,c,b,\Pi,\Gamma, P, \theta, l)|_{\mathcal{O}_{(\mu,a_1,a_2)}
\times \mathbb{R}^2\times \mathbb{R}^{2*}} \cdot \bar{\varepsilon}\\
& =\frac{1}{2}[\frac{(\bar{\varepsilon}_1-
\bar{\varepsilon}_{12})^2}{\bar{I}_1}+\frac{(\bar{\varepsilon}_2-
\bar{\varepsilon}_{13})^2}{\bar{I}_2}
+\frac{\bar{\varepsilon}_3^2}{\bar{I}_3}+
\frac{\bar{\varepsilon}_7^2}{m_1}+\frac{\bar{\varepsilon}_8^2}{m_2}+\frac{\bar{\varepsilon}_9^2}{m_3}
+\frac{\bar{\varepsilon}_{12}^2}{J_1}+\frac{\bar{\varepsilon}_{13}^2}{J_2}] \\
& \;\;\;\;\;\; + gh(\bar{\varepsilon}_4\cdot\chi_1+ \bar{\varepsilon}_5\cdot\chi_2+
\bar{\varepsilon}_6\cdot\chi_3),
\end{align*} and the vector field
\begin{align*}
& X_{h_{(\mu,a_1,a_2)}}(\Pi) \cdot \bar{\varepsilon}
= \{\Pi,h_{(\mu,a_1,a_2)}\}_{-}|_{\mathcal{O}_{(\mu,a_1,a_2)}\times \mathbb{R}^2\times \mathbb{R}^{2*}}\cdot \bar{\varepsilon}\\
&=(\Pi_1,\Pi_2,\Pi_3)\times (\frac{(\Pi_1- l_1)}{ \bar{I}_1},
\frac{(\Pi_2- l_2)}{ \bar{I}_2}, \frac{\Pi_3}{\bar{I}_3})\cdot
\bar{\varepsilon}
+gh(\Gamma_1,\Gamma_2,\Gamma_3)\times (\chi_1,\chi_2,\chi_3)\cdot \bar{\varepsilon}\\
& \;\;\;\;\;\; +(P_1,P_2,P_3)\times (\frac{P_1}{m_1},\frac{P_2}{m_2},\frac{P_3}{m_3})\cdot \bar{\varepsilon}\\
&= ( \frac{(\bar{I}_2-\bar{I}_3)\bar{\varepsilon}_2\bar{\varepsilon}_3+
\bar{I}_3\bar{\varepsilon}_3\bar{\varepsilon}_{13}}{\bar{I}_2\bar{I}_3}
+  \frac{(m_2-m_3)\bar{\varepsilon}_8\bar{\varepsilon}_9}{m_2m_3}+
gh(\bar{\varepsilon}_5\chi_3-\bar{\varepsilon}_6\chi_2), \\
& \;\;\;\;\;\; \frac{(\bar{I}_3-\bar{I}_1)\bar{\varepsilon}_3
\bar{\varepsilon}_1-
\bar{I}_3\bar{\varepsilon}_3\bar{\varepsilon}_{12}}{\bar{I}_3\bar{I}_1}
+  \frac{(m_3-m_1)\bar{\varepsilon}_9\bar{\varepsilon}_7}{m_3m_1}+
gh(\bar{\varepsilon}_6\chi_1-\bar{\varepsilon}_4\chi_3),\\
& \;\;\;\;\;\;
\frac{(\bar{I}_1-\bar{I}_2)\bar{\varepsilon}_1\bar{\varepsilon}_2-\bar{I}_1\bar{\varepsilon}_1
\bar{\varepsilon}_{13}+
\bar{I}_2\bar{\varepsilon}_2\bar{\varepsilon}_{12}}{\bar{I}_1\bar{I}_2}
+  \frac{(m_1-m_2)\bar{\varepsilon}_7\bar{\varepsilon}_8}{m_1m_2}+
gh(\bar{\varepsilon}_4\chi_2-\bar{\varepsilon}_5\chi_1) ),
\end{align*}

\begin{align*}
& X_{h_{(\mu,a_1,a_2)}}(\Gamma) \cdot
\bar{\varepsilon}
= \{\Gamma, h_{(\mu,a_1,a_2)}\}_{-}|_{\mathcal{O}_{(\mu,a_1,a_2)}
\times \mathbb{R}^2\times \mathbb{R}^{2*}}\cdot \bar{\varepsilon}\\
& =(\Gamma_1, \Gamma_2, \Gamma_3)\times (\frac{(\Pi_1- l_1)}{
\bar{I}_1}, \frac{(\Pi_2- l_2)}{ \bar{I}_2},
\frac{\Pi_3}{\bar{I}_3})\cdot \bar{\varepsilon}\\
&= ( \frac{\bar{I}_2\bar{\varepsilon}_5\bar{\varepsilon}_3
-\bar{I}_3\bar{\varepsilon}_6\bar{\varepsilon}_2+
\bar{I}_3\bar{\varepsilon}_3\bar{\varepsilon}_{13}}{\bar{I}_2\bar{I}_3}, \;\;
\frac{\bar{I}_3\bar{\varepsilon}_6\bar{\varepsilon}_1-\bar{I}_1\bar{\varepsilon}_4\bar{\varepsilon}_3-
\bar{I}_3\bar{\varepsilon}_3\bar{\varepsilon}_{12}}{\bar{I}_3\bar{I}_1}, \\
& \;\;\;\;\;\;
\frac{\bar{I}_1\bar{\varepsilon}_4\bar{\varepsilon}_2-\bar{I}_2\bar{\varepsilon}_5\bar{\varepsilon}_1-
\bar{I}_1\bar{\varepsilon}_1\bar{\varepsilon}_{13} + \bar{I}_2\bar{\varepsilon}_2\bar{\varepsilon}_{12}}{\bar{I}_1\bar{I}_2} ),
\end{align*}

\begin{align*}
& X_{h_{(\mu,a_1,a_2)}}(P) \cdot
\bar{\varepsilon}
= \{P,h_{(\mu,a_1,a_2)}\}_{-}|_{\mathcal{O}_{(\mu,a_1,a_2)}\times \mathbb{R}^2\times \mathbb{R}^{2*}}\cdot \bar{\varepsilon}\\
& =(P_1,P_2,P_3)\times (\frac{(\Pi_1- l_1)}{
\bar{I}_1}, \frac{(\Pi_2- l_2)}{ \bar{I}_2},
\frac{\Pi_3}{\bar{I}_3})\cdot \bar{\varepsilon}\\
&= ( \frac{\bar{I}_2\bar{\varepsilon}_8\bar{\varepsilon}_3-\bar{I}_3\bar{\varepsilon}_9\bar{\varepsilon}_2+
\bar{I}_3\bar{\varepsilon}_3\bar{\varepsilon}_{13}}{\bar{I}_2\bar{I}_3}, \;\;
\frac{\bar{I}_3\bar{\varepsilon}_9\bar{\varepsilon}_1-\bar{I}_1\bar{\varepsilon}_7\bar{\varepsilon}_3-
\bar{I}_3\bar{\varepsilon}_3\bar{\varepsilon}_{12}}{\bar{I}_3\bar{I}_1}, \\
& \;\;\;\;\;\;
\frac{\bar{I}_1\bar{\varepsilon}_7\bar{\varepsilon}_2-\bar{I}_2\bar{\varepsilon}_8\bar{\varepsilon}_1-
\bar{I}_1\bar{\varepsilon}_1\bar{\varepsilon}_{13} + \bar{I}_2\bar{\varepsilon}_2\bar{\varepsilon}_{12}}{\bar{I}_1\bar{I}_2} ),
\end{align*}

\begin{align*}
& X_{h_{(\mu,a_1,a_2)}}(\theta) \cdot \bar{\varepsilon}=\{\theta,\;
h_{(\mu,a_1,a_2)}\}_{-}|_{\mathcal{O}_{(\mu,a_1,a_2)} \times \mathbb{R}^2\times
\mathbb{R}^2}\cdot \bar{\varepsilon}\\
& = \sum_{k=1}^2(\frac{\partial \theta}{\partial \theta_k}
\frac{\partial (h_{(\mu,a_1,a_2)})}{\partial l_k}- \frac{\partial
(h_{(\mu,a_1,a_2)})}{\partial
\theta_k}\frac{\partial \theta}{\partial l_k})\cdot \bar{\varepsilon}
= (-\frac{(\bar{\varepsilon}_1- \bar{\varepsilon}_{12})}{ \bar{I}_1}+
\frac{\bar{\varepsilon}_{12}}{J_1}, \; -\frac{(\bar{\varepsilon}_2-
\bar{\varepsilon}_{13})}{ \bar{I}_2}+ \frac{\bar{\varepsilon}_{13}}{J_2}),
\end{align*}

\begin{align*}
& X_{h_{(\mu,a_1,a_2)}}(l) \cdot \bar{\varepsilon}=\{l,\;
h_{(\mu,a_1,a_2)}\}_{-}|_{\mathcal{O}_{(\mu,a_1,a_2)} \times \mathbb{R}^2\times
\mathbb{R}^2}\cdot \bar{\varepsilon}\\
& = \sum_{k=1}^2(\frac{\partial l}{\partial \theta_k} \frac{\partial
(h_{(\mu,a_1,a_2)})}{\partial l_k}- \frac{\partial (h_{(\mu,a_1,a_2)})}{\partial
\theta_k}\frac{\partial l}{\partial l_k})\cdot \bar{\varepsilon}=(0,0).
\end{align*}
On the other hand, from the expressions of the dynamical vector field $\tilde{X}$
and Hamiltonian vector field $X_H$, we have that
\begin{align*}
\tilde{X}(\Pi, \Gamma, P, \theta, l)^\varepsilon
& =T\pi_Q \cdot \tilde{X}\cdot\varepsilon(\Pi, \Gamma,P,\theta, l)\\
& =T\pi_Q \cdot (X_H+ \textnormal{vlift}(u))\cdot\varepsilon (\Pi, \Gamma, P,\theta, l)\\
& =T\pi_Q \cdot X_H \cdot\varepsilon (\Pi, \Gamma, P,\theta, l)
= X_H\cdot\varepsilon(\Pi,\Gamma, P,\theta, l),
\end{align*}
that is,
\begin{align*}
\tilde{X}(\Pi)^\varepsilon & = X_H(\Pi)\cdot\varepsilon \\
& = ( \frac{(\bar{I}_2-\bar{I}_3)\varepsilon_{11}\varepsilon_{12}+
\bar{I}_3\varepsilon_{12}\varepsilon_{22}}{\bar{I}_2\bar{I}_3}
+ \frac{(m_2-m_3)\varepsilon_{17}\varepsilon_{18}}{m_2m_3}+
gh(\varepsilon_{14}\chi_3-\varepsilon_{15}\chi_2),\\
& \;\;\;\;\;\;
\frac{(\bar{I}_3-\bar{I}_1)\varepsilon_{12}\varepsilon_{10}-
\bar{I}_3\varepsilon_{12}\varepsilon_{21}}{\bar{I}_3\bar{I}_1}
+ \frac{(m_3-m_1)\varepsilon_{18}\varepsilon_{16}}{m_3m_1}
+ gh(\varepsilon_{15}\chi_1-\varepsilon_{13}\chi_3), \\
& \;\;\;\;\;\;
\frac{(\bar{I}_1-\bar{I}_2)\varepsilon_{10}\varepsilon_{11}-
\bar{I}_1\varepsilon_{10}\varepsilon_{22}
 + \bar{I}_2\varepsilon_{11}\varepsilon_{21}}{\bar{I}_1\bar{I}_2}
+ \frac{(m_1-m_2)\varepsilon_{16}\varepsilon_{17}}{m_1m_2}
+ gh(\varepsilon_{13}\chi_2-\varepsilon_{14}\chi_1) ),
\end{align*}

\begin{align*}
\tilde{X}(\Gamma)^\varepsilon & = X_H(\Gamma)\cdot\varepsilon \\
&= ( \frac{\bar{I}_2\varepsilon_{14}\varepsilon_{12}-\bar{I}_3\varepsilon_{15}\varepsilon_{11}+
\bar{I}_3\varepsilon_{12}\varepsilon_{22}}{\bar{I}_2\bar{I}_3}, \;\;
\frac{\bar{I}_3\varepsilon_{15}\varepsilon_{10}-\bar{I}_1\varepsilon_{13}\varepsilon_{12}-
\bar{I}_3\varepsilon_{12}\varepsilon_{21}}{\bar{I}_3\bar{I}_1}, \\
& \;\;\;\;\;\;
\frac{\bar{I}_1\varepsilon_{13}\varepsilon_{11}-\bar{I}_2\varepsilon_{14}\varepsilon_{10}-
\bar{I}_1\varepsilon_{10}\varepsilon_{22} + \bar{I}_2\varepsilon_{11}\varepsilon_{21}}{\bar{I}_1\bar{I}_2} ),
\end{align*}

\begin{align*}
\tilde{X}(P)^\varepsilon & = X_H(P)\cdot\varepsilon \\
&= ( \frac{\bar{I}_2\varepsilon_{17}\varepsilon_{12}-\bar{I}_3\varepsilon_{18}\varepsilon_{11}+
\bar{I}_3\varepsilon_{12}\varepsilon_{22}}{\bar{I}_2\bar{I}_3}, \;\;
\frac{\bar{I}_3\varepsilon_{18}\varepsilon_{10}-\bar{I}_1\varepsilon_{16}\varepsilon_{12}-
\bar{I}_3\varepsilon_{12}\varepsilon_{21}}{\bar{I}_3\bar{I}_1}, \\
& \;\;\;\;\;\;
\frac{\bar{I}_1\varepsilon_{16}\varepsilon_{11}-\bar{I}_2\varepsilon_{17}\varepsilon_{10}-
\bar{I}_1\varepsilon_{10}\varepsilon_{22} + \bar{I}_2\varepsilon_{11}\varepsilon_{21}}{\bar{I}_1\bar{I}_2} ),
\end{align*}
\begin{align*}
\tilde{X}(\theta)^\varepsilon & = X_H(\theta)\cdot\varepsilon
= ( -\frac{(\varepsilon_{10}- \varepsilon_{21})}{ \bar{I}_1}
+\frac{\varepsilon_{21}}{J_1},\;\; -\frac{(\varepsilon_{11}- \varepsilon_{22})}{ \bar{I}_2}
+\frac{\varepsilon_{22}}{J_2} ),
\end{align*}
\begin{align*}
\tilde{X}(l)^\varepsilon & = X_H(l)\cdot\varepsilon =(0,0),
\end{align*}
then we have that
\begin{align*}
T\bar{\gamma}\cdot \tilde{X}(\Pi)^\varepsilon
&= ( \frac{(\bar{I}_2-\bar{I}_3)\bar{\gamma}_2\bar{\gamma}_3+
\bar{I}_3\bar{\gamma}_3\bar{\gamma}_{13}}{\bar{I}_2\bar{I}_3}
+  \frac{(m_2-m_3)\bar{\gamma}_8\bar{\gamma}_9}{m_2m_3}+
gh(\bar{\gamma}_5\chi_3-\bar{\gamma}_6\chi_2),\\
& \;\;\;\;\;\; \frac{(\bar{I}_3-\bar{I}_1)\bar{\gamma}_3
\bar{\gamma}_1-
\bar{I}_3\bar{\gamma}_3\bar{\gamma}_{12}}{\bar{I}_3\bar{I}_1}
+  \frac{(m_3-m_1)\bar{\gamma}_9\bar{\gamma}_7}{m_3m_1}+
gh(\bar{\gamma}_6\chi_1-\bar{\gamma}_4\chi_3),\\
& \;\;\;\;\;\;
\frac{(\bar{I}_1-\bar{I}_2)\bar{\gamma}_1\bar{\gamma}_2-\bar{I}_1\bar{\gamma}_1
\bar{\gamma}_{13}+
\bar{I}_2\bar{\gamma}_2\bar{\gamma}_{12}}{\bar{I}_1\bar{I}_2}
+  \frac{(m_1-m_2)\bar{\gamma}_7\bar{\gamma}_8}{m_1m_2}+
gh(\bar{\gamma}_4\chi_2-\bar{\gamma}_5\chi_1) ),
\end{align*}

\begin{align*}
T\bar{\gamma}\cdot \tilde{X}(\Gamma)^\varepsilon
&= ( \frac{\bar{I}_2\bar{\gamma}_5\bar{\gamma}_3-\bar{I}_3\bar{\gamma}_6\bar{\gamma}_2+
\bar{I}_3\bar{\gamma}_3\bar{\gamma}_{13}}{\bar{I}_2\bar{I}_3}, \;\;
\frac{\bar{I}_3\bar{\gamma}_6\bar{\gamma}_1-\bar{I}_1\bar{\gamma}_4\bar{\gamma}_3-
\bar{I}_3\bar{\gamma}_3\bar{\gamma}_{12}}{\bar{I}_3\bar{I}_1}, \\
& \;\;\;\;\;\;
\frac{\bar{I}_1\bar{\gamma}_4\bar{\gamma}_2-\bar{I}_2\bar{\gamma}_5\bar{\gamma}_1-
\bar{I}_1\bar{\gamma}_1\bar{\gamma}_{13}
+ \bar{I}_2\bar{\gamma}_2\bar{\gamma}_{12}}{\bar{I}_1\bar{I}_2} ),
\end{align*}

\begin{align*}
T\bar{\gamma}\cdot \tilde{X}(P)^\varepsilon
&= ( \frac{\bar{I}_2\bar{\gamma}_8\bar{\gamma}_3-\bar{I}_3\bar{\gamma}_9\bar{\gamma}_2+
\bar{I}_3\bar{\gamma}_3\bar{\gamma}_{13}}{\bar{I}_2\bar{I}_3}, \;\;
\frac{\bar{I}_3\bar{\gamma}_9\bar{\gamma}_1-\bar{I}_1\bar{\gamma}_7\bar{\gamma}_3-
\bar{I}_3\bar{\gamma}_3\bar{\gamma}_{12}}{\bar{I}_3\bar{I}_1}, \\
& \;\;\;\;\;\;
\frac{\bar{I}_1\bar{\gamma}_7\bar{\gamma}_2-\bar{I}_2\bar{\gamma}_8\bar{\gamma}_1-
\bar{I}_1\bar{\gamma}_1\bar{\gamma}_{13}
+ \bar{I}_2\bar{\gamma}_2\bar{\gamma}_{12}}{\bar{I}_1\bar{I}_2} ),
\end{align*}
\begin{align*}
T\bar{\gamma}\cdot \tilde{X}(\theta)^\varepsilon
& = (-\frac{(\bar{\gamma}_1- \bar{\gamma}_{12})}{ \bar{I}_1}+
\frac{\bar{\gamma}_{12}}{J_1}, \; -\frac{(\bar{\gamma}_2-
\bar{\gamma}_{13})}{ \bar{I}_2}+ \frac{\bar{\gamma}_{13}}{J_2}),
\end{align*}
\begin{align*}
T\bar{\gamma}\cdot \tilde{X}(l)^\varepsilon=(0,0).
\end{align*}
Note that $$T\bar{\lambda}\cdot \tilde{X} \cdot \varepsilon=T\pi_{(\mu,a_1,a_2)}\cdot T\lambda \cdot (X_H+ \textnormal{vlift}(u))\cdot\varepsilon
=T\pi_{(\mu,a_1,a_2)}\cdot T\gamma \cdot T\pi_Q \cdot (X_H+ \textnormal{vlift}(u))\cdot\varepsilon
=T\bar{\lambda}\cdot X_H \cdot \varepsilon,$$ that is,
\begin{align*}
T\bar{\lambda}\cdot \tilde{X}(\Pi) \cdot \varepsilon & =T\bar{\lambda}\cdot X_H(\Pi) \cdot \varepsilon\\
&= ( \frac{(\bar{I}_2-\bar{I}_3)\bar{\lambda}_2\bar{\lambda}_3+
\bar{I}_3\bar{\lambda}_3\bar{\lambda}_{13}}{\bar{I}_2\bar{I}_3}
+  \frac{(m_2-m_3)\bar{\lambda}_8\bar{\lambda}_9}{m_2m_3}+
gh(\bar{\lambda}_5\chi_3-\bar{\lambda}_6\chi_2),\\
& \;\;\;\;\;\; \frac{(\bar{I}_3-\bar{I}_1)\bar{\lambda}_3
\bar{\lambda}_1-
\bar{I}_3\bar{\lambda}_3\bar{\lambda}_{12}}{\bar{I}_3\bar{I}_1}
+  \frac{(m_3-m_1)\bar{\lambda}_9\bar{\lambda}_7}{m_3m_1}+
gh(\bar{\lambda}_6\chi_1-\bar{\lambda}_4\chi_3),\\
& \;\;\;\;\;\;
\frac{(\bar{I}_1-\bar{I}_2)\bar{\lambda}_1\bar{\lambda}_2-\bar{I}_1\bar{\lambda}_1
\bar{\lambda}_{13}+
\bar{I}_2\bar{\lambda}_2\bar{\lambda}_{12}}{\bar{I}_1\bar{I}_2}
+  \frac{(m_1-m_2)\bar{\lambda}_7\bar{\lambda}_8}{m_1m_2}+
gh(\bar{\lambda}_4\chi_2-\bar{\lambda}_5\chi_1) ),
\end{align*}

\begin{align*}
T\bar{\lambda}\cdot \tilde{X}(\Gamma) \cdot \varepsilon &
=T\bar{\lambda}\cdot X_H(\Gamma) \cdot \varepsilon\\
&= ( \frac{\bar{I}_2\bar{\lambda}_5\bar{\lambda}_3-\bar{I}_3\bar{\lambda}_6\bar{\lambda}_2+
\bar{I}_3\bar{\lambda}_3\bar{\lambda}_{13}}{\bar{I}_2\bar{I}_3}, \;\;
\frac{\bar{I}_3\bar{\lambda}_6\bar{\lambda}_1-\bar{I}_1\bar{\lambda}_4\bar{\lambda}_3-
\bar{I}_3\bar{\lambda}_3\bar{\lambda}_{12}}{\bar{I}_3\bar{I}_1}, \\
& \;\;\;\;\;\;
\frac{\bar{I}_1\bar{\lambda}_4\bar{\lambda}_2-\bar{I}_2\bar{\lambda}_5\bar{\lambda}_1-
\bar{I}_1\bar{\lambda}_1\bar{\lambda}_{13}
+ \bar{I}_2\bar{\lambda}_2\bar{\lambda}_{12}}{\bar{I}_1\bar{I}_2} ),
\end{align*}

\begin{align*}
T\bar{\lambda}\cdot \tilde{X}(P) \cdot \varepsilon & =T\bar{\lambda}\cdot X_H(P) \cdot \varepsilon\\
&= ( \frac{\bar{I}_2\bar{\lambda}_8\bar{\lambda}_3-\bar{I}_3\bar{\lambda}_9\bar{\lambda}_2+
\bar{I}_3\bar{\lambda}_3\bar{\lambda}_{13}}{\bar{I}_2\bar{I}_3}, \;\;
\frac{\bar{I}_3\bar{\lambda}_9\bar{\lambda}_1-\bar{I}_1\bar{\lambda}_7\bar{\lambda}_3-
\bar{I}_3\bar{\lambda}_3\bar{\lambda}_{12}}{\bar{I}_3\bar{I}_1}, \\
& \;\;\;\;\;\;
\frac{\bar{I}_1\bar{\lambda}_7\bar{\lambda}_2-\bar{I}_2\bar{\lambda}_8\bar{\lambda}_1-
\bar{I}_1\bar{\lambda}_1\bar{\lambda}_{13}
+ \bar{I}_2\bar{\lambda}_2\bar{\lambda}_{12}}{\bar{I}_1\bar{I}_2} ),
\end{align*}

\begin{align*}
T\bar{\lambda}\cdot \tilde{X}(\theta) \cdot \varepsilon=T\bar{\lambda}\cdot X_H(\theta) \cdot \varepsilon
= (-\frac{(\bar{\lambda}_1- \bar{\lambda}_{12})}{ \bar{I}_1}+
\frac{\bar{\lambda}_{12}}{J_1}, \; -\frac{(\bar{\lambda}_2-
\bar{\lambda}_{13})}{ \bar{I}_2}+ \frac{\bar{\lambda}_{13}}{J_2}),
\end{align*}
\begin{align*}
T\bar{\lambda}\cdot \tilde{X}(l) \cdot \varepsilon=T\bar{\lambda}\cdot X_H(l) \cdot \varepsilon=(0,0).
\end{align*}
Thus, when we choose that
$(\Pi,\Gamma, P,\theta,l)\in
\mathcal{O}_{(\mu,a_1,a_2)}\times \mathbb{R}^2 \times \mathbb{R}^{2*}, $ and
 $(\varepsilon_{10},\varepsilon_{11},\varepsilon_{12})=\Pi=(\Pi_1,\Pi_2,\Pi_3)
 =(\bar{\gamma}_1,\bar{\gamma}_2,\bar{\gamma}_3)=
(\bar{\varepsilon}_1,\bar{\varepsilon}_2,\bar{\varepsilon}_3)=
(\bar{\lambda}_1,\bar{\lambda}_2,\bar{\lambda}_3), $ and
$(\varepsilon_{13},\varepsilon_{14},\varepsilon_{15})=\Gamma
=(\Gamma_1,\Gamma_2,\Gamma_3)=(\bar{\gamma}_4,\bar{\gamma}_5,\bar{\gamma}_6)=
(\bar{\varepsilon}_4,\bar{\varepsilon}_5,\bar{\varepsilon}_6)=
(\bar{\lambda}_4,\bar{\lambda}_5,\bar{\lambda}_6), $
$(\varepsilon_{16},\varepsilon_{17},\varepsilon_{18})=P=(P_1,P_2,P_3)
=(\bar{\gamma}_7,\bar{\gamma}_8,\bar{\gamma}_9)=
(\bar{\varepsilon}_7,\bar{\varepsilon}_8,\bar{\varepsilon}_9)=
(\bar{\lambda}_7,\bar{\lambda}_8,\bar{\lambda}_9), $
and $(\varepsilon_{19},\varepsilon_{20})=\theta= (\theta_1,\theta_2)=(\bar{\gamma}_{10},\bar{\gamma}_{11})
=(\bar{\varepsilon}_{10},\bar{\varepsilon}_{11})=(\bar{\lambda}_{10},\bar{\lambda}_{11}),$
$ (\varepsilon_{21},\varepsilon_{22})=l=(l_1,l_2)=(\bar{\gamma}_{12},\bar{\gamma}_{13})
=(\bar{\varepsilon}_{12},\bar{\varepsilon}_{13})=(\bar{\lambda}_{12},\bar{\lambda}_{13}). $ we must have that
\begin{align*}
& T\bar{\gamma}\cdot \tilde{X}(\Pi)^\varepsilon =X_{h_{(\mu,a_1,a_2)}}(\Pi) \cdot \bar{\varepsilon}
=T\bar{\lambda}\cdot \tilde{X}(\Pi) \cdot \varepsilon, \\
& T\bar{\gamma}\cdot \tilde{X}(\Gamma)^\varepsilon=X_{h_{(\mu,a_1,a_2)}}(\Gamma) \cdot \bar{\varepsilon}
=T\bar{\lambda}\cdot \tilde{X}(\Gamma) \cdot \varepsilon,\\
& T\bar{\gamma}\cdot \tilde{X}(P)^\varepsilon=X_{h_{(\mu,a_1,a_2)}}(P) \cdot \bar{\varepsilon}
=T\bar{\lambda}\cdot \tilde{X}(P) \cdot \varepsilon,\\
& T\bar{\gamma}\cdot \tilde{X}(\theta)^\varepsilon =X_{h_{(\mu,a_1,a_2)}}(\theta) \cdot \bar{\varepsilon}
=T\bar{\lambda}\cdot \tilde{X}(\theta) \cdot \varepsilon, \\
& T\bar{\gamma}\cdot \tilde{X}(l)^\varepsilon =X_{h_{(\mu,a_1,a_2)}}(l) \cdot \bar{\varepsilon}
=T\bar{\lambda}\cdot \tilde{X}(l) \cdot \varepsilon.
\end{align*}
Since the map $\varepsilon: T^* Q \rightarrow T^* Q $ is symplectic, then
$T\bar{\varepsilon}\cdot X_{h_{(\mu,a_1,a_2)} \cdot \bar{\varepsilon}}
=X_{h_{(\mu,a_1,a_2)}} \cdot \bar{\varepsilon}. $
Thus, in this case, we must have that
$\varepsilon$ and $\bar{\varepsilon} $ are the solution of the Type II of
Hamilton-Jacobi equation
$T\bar{\gamma}\cdot \tilde{X}^\varepsilon= X_{h_{(\mu,a_1,a_2)}}\cdot \bar{\varepsilon}, $
for the $R_p$-reduced controlled underwater vehicle-rotor system
$(\mathcal{O}_{(\mu,a_1,a_2)}\times \mathbb{R}^2\times \mathbb{R}^{2*},
\omega^{-}_{\mathcal{O}_{(\mu,a_1,a_2)}\times \mathbb{R}^2\times \mathbb{R}^{2*}},
h_{(\mu,a_1,a_2)}, u_{(\mu,a_1,a_2)})$, if and only if they satisfy
the equation $T\bar{\varepsilon}\cdot(X_{h_{(\mu,a_1,a_2)} \cdot \bar{\varepsilon}})
= T\bar{\lambda}\cdot \tilde{X}\cdot\varepsilon. $\\

To sum up the above discussion, we have the following Theorem 5.2.
For convenience, the maps involved in
the following theorem are shown in Diagram-4.

\begin{center}
\hskip 0cm \xymatrix{ \mathbf{J}_Q^{-1}(\mu,a_1,a_2) \ar[r]^{i_{(\mu,a_1,a_2)}} & T^* Q
\ar[d]_{X_{H\cdot \varepsilon}} \ar[dr]^{\tilde{X}^\varepsilon} \ar[r]^{\pi_Q}
& Q \ar[d]^{\tilde{X}^\gamma} \ar[r]^{\gamma}
& T^*Q \ar[d]_{\tilde{X}} \ar[dr]^{X_{h_{(\mu,a_1,a_2)} \cdot\bar{\varepsilon}}} \ar[r]^{\pi_{(\mu,a_1,a_2)}}
& \;\;\;\;\;\; \mathcal{O}_{(\mu,a_1,a_2)}\times \mathbb{R}^2\times \mathbb{R}^{2*} \ar[d]^{X_{h_{(\mu,a_1,a_2)}}} \\
& T(T^*Q)  & TQ \ar[l]^{T\gamma}
& T(T^*Q) \ar[l]^{T\pi_Q} \ar[r]_{T\pi_{(\mu,a_1,a_2)}}
& \;\;\;\;\;\; T(\mathcal{O}_{(\mu,a_1,a_2)}\times \mathbb{R}^2\times \mathbb{R}^{2*})}
\end{center}
$$\mbox{Diagram-4}$$

\begin{theo}
In the case of non-coincident centers of buoyancy and gravity,
if the 5-tuple \\ $(T^\ast Q,\textmd{SE}(3)\circledS \mathbb{R}^3,\omega_Q,H,u), $ where $Q=
\textmd{SE}(3)\circledS \mathbb{R}^3\times \mathbb{R}^2, $ is a regular point reducible
underwater vehicle-rotor system with the control torque $u$ acting on the rotors,
then for a point $(\mu,a_1,a_2)\in \mathfrak{se}^\ast(3)\circledS
\mathbb{R}^{3*}$, the regular value of the momentum map $\mathbf{J}_Q:
T^* Q \cong \textmd{SE}(3)\circledS \mathbb{R}^3 \times
\mathfrak{se}^\ast(3)\circledS \mathbb{R}^{3*} \times
\mathbb{R}^2\times \mathbb{R}^{2*} \to \mathfrak{se}^\ast(3)\circledS
\mathbb{R}^{3*}$, the $R_p$-reduced controlled underwater vehicle-rotor system is the 4-tuple
$(\mathcal{O}_{(\mu,a_1,a_2)} \times \mathbb{R}^2 \times \mathbb{R}^{2*},
\tilde{\omega}^{-}_{\mathcal{O}_{(\mu,a_1,a_2)} \times
\mathbb{R}^2 \times \mathbb{R}^{2*}},h_{(\mu,a_1,a_2)},u_{(\mu,a_1,a_2)})$,
where $ \mathcal{O}_{(\mu,a_1,a_2)}
\subset \mathfrak{se}^\ast(3)\circledS \mathbb{R}^{3*}$ is the co-adjoint orbit
of the semidirect product Lie group $\textmd{SE}(3)\circledS \mathbb{R}^3$.
Assume that
$\gamma: Q \rightarrow T^* Q $ is an one-form on $Q=\textmd{SE}(3)\circledS
\mathbb{R}^3 \times \mathbb{R}^2 $, and $\lambda=\gamma \cdot \pi_Q:
T^* Q \rightarrow T^* Q, $ and $\varepsilon: T^* Q \rightarrow T^* Q $ is a
$\textmd{(SE}(3)\circledS \mathbb{R}^3)_{(\mu,a_1,a_2)}$-invariant symplectic map,
where $(\textmd{SE}(3)\circledS \mathbb{R}^3)_{(\mu,a_1,a_2)}$ is
the isotropy subgroup of co-adjoint $\textmd{SE}(3)\circledS
\mathbb{R}^3$-action at the point $(\mu,a_1,a_2)$.
Denote
$\tilde{X}^\gamma = T\pi_Q \cdot \tilde{X}\cdot \gamma$, and
$\tilde{X}^\varepsilon = T\pi_Q \cdot \tilde{X}\cdot \varepsilon$,
where $\tilde{X}=X_{(T^\ast Q,\textmd{SE}(3)\circledS
\mathbb{R}^3,\omega_Q,H,u)}$ is the dynamical vector field of
the controlled underwater vehicle-rotor system $(T^\ast Q,\textmd{SE}(3)\circledS
\mathbb{R}^3,\omega_Q,H,u)$.
Moreover, assume that
and $\textmd{Im}(\gamma)\subset
\mathbf{J}_Q^{-1}((\mu,a_1,a_2)), $ and it is
$(\textmd{SE}(3)\circledS \mathbb{R}^3)_{(\mu,a_1,a_2)}$-invariant,
and $\varepsilon (\mathbf{J}_Q^{-1}(\mu,a_1,a_2))\subset \mathbf{J}_Q^{-1}(\mu,a_1,a_2). $
 Denote by
$\bar{\gamma}=\pi_{(\mu,a_1,a_2)}(\gamma): Q=\textmd{SE}(3)\circledS
\mathbb{R}^3 \times \mathbb{R}^2 \rightarrow \tilde{\mathcal{O}}, $
and $\bar{\lambda}=\pi_{(\mu,a_1,a_2)}(\lambda):
T^* Q \rightarrow \tilde{\mathcal{O}} $,
 $\bar{\varepsilon}=\pi_{(\mu,a_1,a_2)}(\varepsilon): \mathbf{J}_Q^{-1}(\mu,a_1,a_2)\rightarrow
\tilde{\mathcal{O}}. $
Then the following two assertions hold:\\
\noindent $(\mathbf{i})$
If the one-form $\gamma: Q \rightarrow T^* Q $ is closed with respect to
$T\pi_Q: TT^* Q \rightarrow TQ, $
then $\bar{\gamma}$ is a solution of the Type I of Hamilton-Jacobi equation
$T\bar{\gamma}\cdot \tilde{X}^\gamma= X_{h_{(\mu,a_1,a_2)}}\cdot \bar{\gamma}; $\\
\noindent $(\mathbf{ii})$
The $\varepsilon$ and $\bar{\varepsilon} $ satisfy the Type II of Hamilton-Jacobi equation
$T\bar{\gamma}\cdot \tilde{X}^\varepsilon= X_{h_{(\mu,a_1,a_2)}}\cdot \bar{\varepsilon}, $
if and only if they satisfy
the equation $T\bar{\varepsilon}\cdot(X_{h_{(\mu,a_1,a_2)} \cdot \bar{\varepsilon}})
= T\bar{\lambda}\cdot \tilde{X}\cdot\varepsilon. $ \hskip 0.3cm $\blacksquare$
\end{theo}

It is worthy of noting that, for the regular point reducible controlled
underwater vehicle-rotor system
$(T^\ast Q,\textmd{SE}(3)\circledS \mathbb{R}^3,\omega_Q, H,u)$
with the $R_p$-reduced controlled underwater vehicle-rotor system
$(\mathcal{O}_{(\mu,a_1,a_2)} \times \mathbb{R}^2 \times
\mathbb{R}^{2*},\tilde{\omega}_{\mathcal{O}_{(\mu,a_1,a_2)} \times \mathbb{R}^2
\times \mathbb{R}^{2*}}^{-},
h_{(\mu,a_1,a_2)}, u_{(\mu,a_1,a_2)}) $, we know that the Hamiltonian vector fields
$X_{H}$ and $X_{h_{(\mu,a_1,a_2)}}$ for the corresponding
Hamiltonian system $(T^*Q,\textmd{SE}(3)\circledS \mathbb{R}^3,\omega_Q, H)$
and its $R_p$-reduced system $(\mathcal{O}_{(\mu,a_1,a_2)} \times \mathbb{R}^2 \times
\mathbb{R}^{2*},\tilde{\omega}_{\mathcal{O}_{(\mu,a_1,a_2)}\times \mathbb{R}^2
\times \mathbb{R}^{2*}}^{-}, h_{(\mu,a_1,a_2)} )$, are $\pi_{(\mu,a_1,a_2)}$-related, that is,
$X_{h_{(\mu,a_1,a_2)}}\cdot \pi_{(\mu,a_1,a_2)}
=T\pi_{(\mu,a_1,a_2)}\cdot X_{H}\cdot i_{(\mu,a_1,a_2)}.$
From the above Theorem 3.8 we can obtain the following Theorem 5.3, which states the relationship
between the solutions of Type II of Hamilton-Jacobi equations and the
regular point reduction.

\begin{theo}
In the case of non-coincident centers of buoyancy and gravity, for the regular point reducible controlled
underwater vehicle-rotor system $(T^\ast Q,\textmd{SE}(3)\circledS \mathbb{R}^3,\omega_Q,H,u)$
with the $R_p$-reduced controlled underwater vehicle-rotor system
$(\mathcal{O}_{(\mu,a_1,a_2)} \times \mathbb{R}^2 \times
\mathbb{R}^{2*},\tilde{\omega}_{\mathcal{O}_{(\mu,a_1,a_2)}\times \mathbb{R}^2
\times \mathbb{R}^{2*}}^{-}, h_{(\mu,a_1,a_2)}, \\ u_{(\mu,a_1,a_2)}) $,
assume that $\gamma: Q \rightarrow T^* Q $ is an one-form on
$Q=\textmd{SE}(3)\circledS \mathbb{R}^3\times \mathbb{R}^2 $, and $\varepsilon:
T^* Q \rightarrow T^* Q $ is a
$(\textmd{SE}(3)\circledS \mathbb{R}^3)_{(\mu,a_1,a_2)}$-invariant symplectic map,
$\bar{\varepsilon}=\pi_{(\mu,a_1,a_2)}(\varepsilon): \mathbf{J}_Q^{-1}(\mu, a_1,a_2 )\rightarrow
\mathcal{O}_{(\mu,a_1,a_2)}\times \mathbb{R}^2\times \mathbb{R}^{2*}. $
Under the hypotheses and notations of Theorem 5.2, then we have that
$\varepsilon$ is a solution of the Type II of Hamilton-Jacobi equation
$T\gamma\cdot \tilde{X}^\varepsilon= X_H\cdot \varepsilon, $ for the
regular point reducible controlled underwater vehicle-rotor system
$(T^\ast Q,\textmd{SE}(3)\circledS \mathbb{R}^3,\omega_Q,H,u), $ if and only if
$\varepsilon$ and $\bar{\varepsilon} $ satisfy the Type II of Hamilton-Jacobi equation
$T\bar{\gamma}\cdot \tilde{X}^\varepsilon= X_{h_{(\mu,a_1,a_2)}}\cdot \bar{\varepsilon}, $ for the
$R_p$-reduced controlled underwater vehicle-rotor system
$(\mathcal{O}_{(\mu,a_1,a_2)} \times \mathbb{R}^2 \times
\mathbb{R}^{2*},\tilde{\omega}_{\mathcal{O}_{(\mu,a_1,a_2)} \times \mathbb{R}^2
\times \mathbb{R}^{2*}}^{-}, h_{(\mu,a_1,a_2)}, u_{(\mu,a_1,a_2)}) $.
 \hskip 0.3cm $\blacksquare$
\end{theo}

\begin{rema}
When the underwater vehicle does not carry any internal rotor, in
this case the configuration space is $Q=W=\textmd{SE}(3)\circledS
\mathbb{R}^3, $ the motion of underwater vehicle is just the
rotation with drift and translation motion of a rigid body, the
above $R_p$-reduced controlled underwater vehicle-rotors system is
just the Marsden-Weinstein reduced underwater vehicle system,
that is, 3-tuple $(\mathcal{O}_{(\mu,a_1,a_2)},
\omega_{\mathcal{O}_{(\mu,a_1,a_2)}}^{-}, h_{(\mu,a_1,a_2)} )$,
where $\mathcal{O}_{(\mu,a_1,a_2)}\subset \mathfrak{se}^\ast(3)\circledS \mathbb{R}^3$
is the co-adjoint orbit of the semidirect product
Lie group $\textmd{SE}(3)\circledS \mathbb{R}^3$,
$\omega_{\mathcal{O}_{(\mu,a_1,a_2)}}^{-}$
is the orbit symplectic form on $\mathcal{O}_{(\mu,a_1,a_2)}$,
which is induced by the (-)-semidirect product Lie-Poisson brackets on
$\mathfrak{se}^\ast(3)\circledS \mathbb{R}^3$,
see Marsden et al. \cite{mamiorpera07, marawe84a, marawe84b},
$h_{(\mu,a_1,a_2)}(\Pi, \Gamma, P)\cdot \pi_{(\mu,a_1,a_2)}
=H(A,c,b,\Pi,\Gamma, P)|_{\mathcal{O}_{(\mu,a_1,a_2)}}.$
From the above Theorem 5.2 we can obtain the two
types of Hamilton-Jacobi equation for the Marsden-Weinstein
reduced underwater vehicle system $(\mathcal{O}_{(\mu,a_1,a_2)},
\omega_{\mathcal{O}_{(\mu,a_1,a_2)}}^{-}, h_{(\mu,a_1,a_2)} )$,
see Wang \cite{wa17}.
Moreover, from the above Theorem 5.3 we can also
state the relationship between the solutions of Type II of Hamilton-Jacobi equations
and the Marsden-Weinstein reduction.
\end{rema}

The theory of controlled mechanical systems is a very important
subject, following the theoretical and applied development of
geometric mechanics, a lot of important problems about this subject
are being explored and studied. In this paper, as an application of
the regular point symplectic reduction and Hamilton-Jacobi theory for
an RCH system with symmetry and momentum map,
we give precisely the regular point reduction and the geometric constraint
conditions of the reduced symplectic forms for the
dynamical vector fields of the regular point reducible controlled underwater vehicle-rotor system,
that is, the two types of Hamilton-Jacobi equations for the reduced
controlled underwater vehicle-rotor system by calculation in
detail, in the cases of coincident and
non-coincident centers of buoyancy and gravity. However,
note that in the cases of coincident and
non-coincident centers of buoyancy and gravity,
since the motions of the controlled underwater vehicle-rotor system are different, and
the configuration spaces, the Hamiltonian functions, the actions of Lie groups,
the $R_p$-reduced symplectic forms and the $R_p$-reduced systems of
the controlled underwater vehicle-rotor system are also different. But,
the two types of Hamilton-Jacobi equations
given by calculation in detail are same,
that is, the internal rules are same by comparing Theorem 4.2 and Theorem 5.2.
In consequence, we reveal the deeply internal
relationships of the geometrical structures of phase spaces, the dynamical
vector fields and controls of the controlled underwater vehicle-rotor system.
It is the key thought of the researches of geometrical mechanics
of the professor Jerrold E. Marsden to explore and study the deeply internal
relationship between the geometrical structure of phase space and the dynamical
vector field of a mechanical system.
It is also our goal of pursuing and inheriting.
In addition, we note also that there have been a lot of
beautiful results of reduction theory of Hamiltonian systems in
celestial mechanics, hydrodynamics and plasma physics. Thus, it is
an important topic to study the application of reduction theory of the
RCH systems in celestial mechanics, hydrodynamics
and plasma physics. These are our goals in future research.\\

\end{document}